\documentclass[11pt]{article} 
\usepackage{latexsym} 
\usepackage{amscd}
\usepackage{amsmath,amsfonts,amssymb,amsthm}
\usepackage{graphicx}

\usepackage[T1]{fontenc}
\usepackage{hyperref}
\usepackage{xcolor}

\usepackage[letterpaper,top=3.5cm,bottom=3.5cm,left=3.5cm,right=3.5cm,marginparwidth=1.75cm]{geometry}

\title{Inverse Spectral Problems for Collapsing Manifolds I: Uniqueness and Stability} 

\author{Yaroslav Kurylev,\, Matti Lassas,\, Jinpeng Lu,\, Takao Yamaguchi}

\date{}

\newcommand{\la}{{\lambda}} 
\newcommand{\R}{{\mathbb R}} 
\newcommand{\Z}{{\mathbb Z}} 
 
\newcommand{\C}{{\mathbb C}} 
\newcommand{\N}{{\mathbb N}}
 
\newcommand{\x}{{\bf x}}

\def\hat{\widehat}
\def\normalized{\bar}
\def\tilde{\widetilde}
\def \bfo {\begin {eqnarray*} }
\def \efo {\end {eqnarray*} }
\def \ba {\begin {eqnarray*} }
\def \ea {\end {eqnarray*} }
\def \beq {\begin {eqnarray}}
\def \eeq {\end {eqnarray}}
\def \bequ {\begin {equation}}
\def \eequ {\end {equation}}
\def \supp {\hbox{supp }}
\def \diam {\hbox{diam}}

\def \det {\hbox{det}}
\def\bra{\langle}
\def\cet{\rangle}
\def \e {\varepsilon}
\def \p {\partial}
\def\M{{\mathcal M}}

\def\a{\alpha}
\def\g{\gamma}

\def\h{H}
\def \d {\delta}
\def\y {{\bf y}}

\newtheorem{definition}{Definition}[section] 
\newtheorem{theorem}[definition]{Theorem} 
\newtheorem{lemma}[definition]{Lemma} 
\newtheorem{problem}[definition]{Problem} 
\newtheorem{proposition}[definition]{Proposition} 
\newtheorem{corollary}[definition]{Corollary} 
\newtheorem{remark}[definition]{Remark} 
\newtheorem{example}[definition]{Example}

\numberwithin{equation}{section}

\begin{document}

\AtEndDocument{\bigskip{\footnotesize%
  
\textsc{Yaroslav Kurylev, Department of Mathematics, University College London, UK} \par  
  
\addvspace{\medskipamount}

  

  \addvspace{\medskipamount}
  \textsc{Matti Lassas and Jinpeng Lu: Department of Mathematics and Statistics, University of Helsinki, Finland.
  Email: matti.lassas@helsinki.fi, jinpeng.lu@helsinki.fi}

\addvspace{\medskipamount}
 \textsc{Takao Yamaguchi: Institute of Mathematics, University of Tsukuba,  Japan. Email: takao@math.tsukuba.ac.jp}

  
}}

\maketitle 

\begin{abstract}
We consider the geometric inverse problem of determining a closed Riemannian manifold from measurements of the heat kernel in an open subset of the manifold. In this paper we analyze the stability of this problem in the class of $n$-dimensional Riemannian manifolds with bounded diameter and sectional curvature.
It is well-known that a sequence in this class of manifolds can collapse to a lower dimensional stratified space when the injectivity radius of the sequence of manifolds goes to zero.
We prove the uniqueness of the inverse problem on the limiting spaces of the collapsing manifolds.
As a result, we obtain stability results for the inverse problem in the class of manifolds with bounded diameter and sectional curvature.
\end{abstract}


 \section{Introduction}\label{introduction 1}

\subsection{Inverse problems for collapsing manifolds}\label{introduction 2}

Let $(M,h)$ be a connected, closed, smooth Riemannian manifold of dimension $n$
with Riemannian metric $h$. Given $p\in M$, we say that the triple $(M,p,h)$ is a pointed Riemannian manifold.  
We denote by $dV_h$ the Riemannian volume element of $(M,h)$, and by
$d\mu_M$ the normalized measure,
\bequ \label{normalised_measure}
d\mu_M= \frac1 {\hbox{Vol}(M)} dV_h,
\eequ
where $\hbox{Vol}(M)$ is the Riemannian volume of $(M, h)$.
We consider the heat kernel $H(x,y,t)$ on $(M,h)$ associated to the Laplacian operator,
\bequ\label{heat kernel1}
\big(\frac{\p}{\p t}+\Delta_M \big) H(\cdotp,y,t)=0\quad\hbox{on }M\times\R_+,\quad  \h(\cdotp,y,0)=\delta_y,
\eequ
where $\delta_y$ is the normalized Dirac delta-distribution, i.e., $\int_M \delta_y(x)\varphi(x)\,d\mu_M(x)=\varphi(y)$ for all $\varphi\in C^\infty(M)$, and
$\Delta_{M}$  is the (nonnegative definite)
Laplace-Beltrami operator  on $(M, h)$, which has the following form in local coordinates,
\bequ \label{000.2}
\Delta_{M} u=-|h|^{-\frac12} \sum_{j,k=1}^n \frac{\p}{\p x_j}\Big( |h|^{\frac12}h^{jk}\frac{\p}{\p x_k} u \Big), \quad |h|=\det\left(h_{jk} \right).
\eequ

In the following, we assume that we are given the values of the heat kernel at
points $\{z_\a: \a=1,2,\dots\}$, which form a dense set in a ball  $B=B_M(p,r)$ of $(M, p,h)$ having center at $p$ and radius $r>0$.
We define pointwise heat data, $PHD$, that is the ordered sequence
\bequ \label{intro-PHD}
PHD=(H_{\a, \beta, \ell})_{\a, \beta, \ell=1}^\infty,\quad
H_{\a, \beta, \ell}:= \h(z_\a,z_\beta,t_\ell),
\eequ
where $\{t_\ell: \ell=1,2,\dots \}$ is a dense set in $\R_+$.
Let us emphasize that the mutual relations of the measurement points $z_\a$, e.g.
the distances between these points, and the topology of the ball $B$ are not
{\it a priori} known.
%

We consider the following generalization of Gel'fand's inverse problem \cite{Gel}.
\begin{problem} \label{Gelfand:1}
{\rm Suppose that the pointwise heat data \eqref{intro-PHD} of two connected, closed, smooth, pointed Riemannian manifolds $(M, p)$ and $(M', p')$ coincide, i.e., for some $r>0$,
$$H(z_{\alpha},z_{\beta},t_{\ell})= H'(z'_{\alpha},z'_{\beta},t_{\ell}), \;\textrm{ for all }z_{\alpha}\in B_M(p,r),\, z'_{\alpha}\in B_{M'}(p',r),\, t_{\ell}\in \R_+.$$
Are the two manifolds isometric?}
\end{problem}

\begin{figure}[h]
\centering
\includegraphics[height=3.5cm]{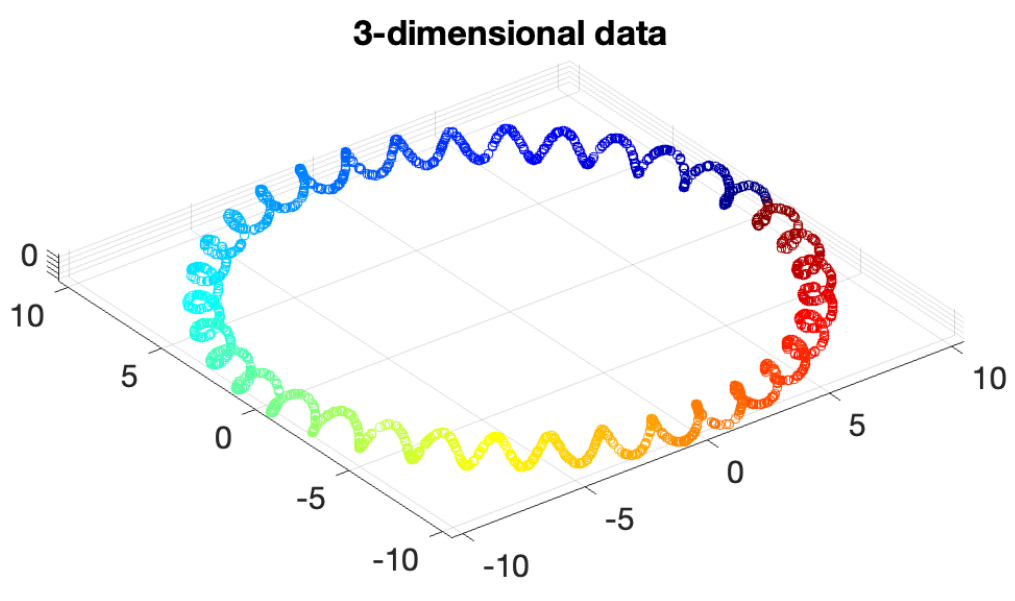}
\includegraphics[height=3.5cm]{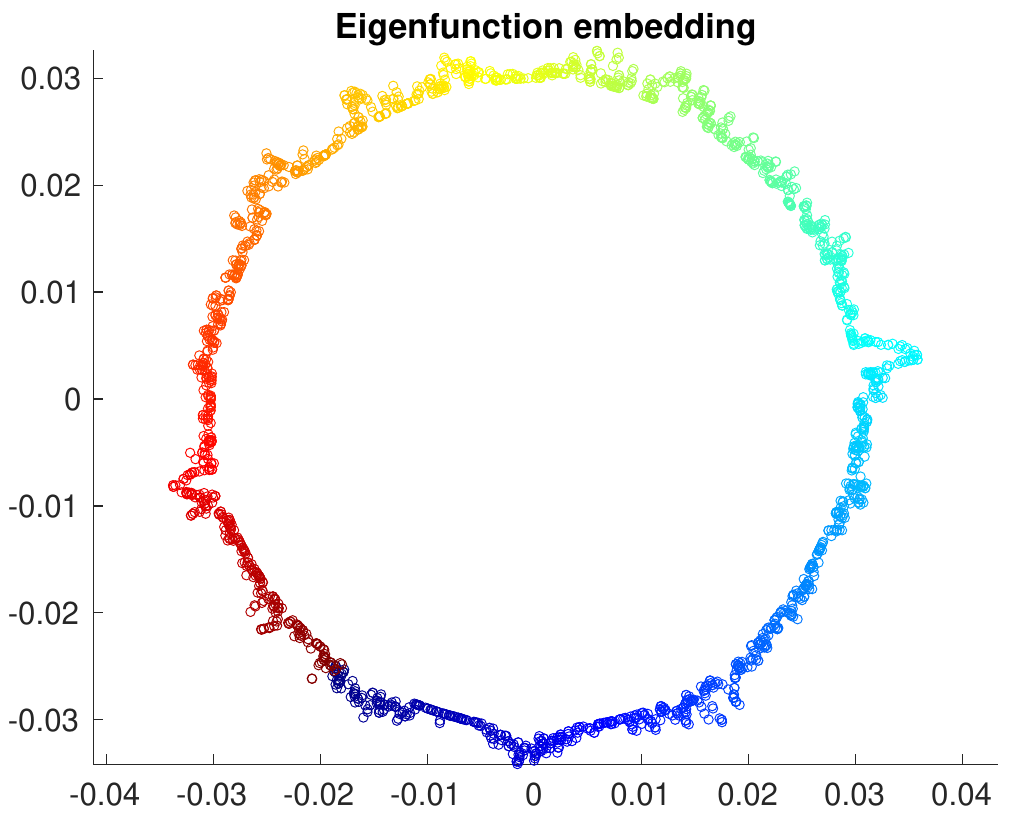}
\caption{Applications in manifold learning for almost collapsed manifolds. {\bf Left:}  
Data points $X=\{x_j\in \R^3:\ j=1,2,\dots, 2048\}$ sampled  from
a 2-dimensional torus embedded in $\R^3$. The points
are sampled along a helical curve.
The distances between the points $x_j$ are defined using the flat metric of the torus $M_{R,r}=\mathbb S^1_R\times \mathbb S^1_r$, with the
larger radius $R=10$ and the smaller radius $r=\frac 12$.
When $r$ is small, the torus $M_{R,r}$
is `almost collapsed' to the circle $\mathbb S^1_R$. 
Note that the points could also be sampled randomly, but we have
used points sampled randomly from a helical curve on $M_{R,r}$ to make the visualization clearer.
{\bf Right:} The image of the 2-dimensional eigenfunction map
$\Phi^{(2,2)}(X)$ associated to the eigenfunctions for an approximation of the heat kernel of $M_{R,r}$.
The image $\Phi^{(2,2)}(X)\subset \R^2$
is close to a circle, that is, it is an approximation
of the limiting space $\mathbb{S}^1_R$, see \eqref{eigenmap}. We study the local version of
Gel'fand's  problem \ref{Gelfand:1}, where the heat kernel $H(x_j,x_{j'},t_\ell)$ are  given at points $x_j$ that do not fill
the whole manifold $M=M_{R,r}$ but only fill a possibly small metric ball $B\subset M$, for example, only the red part of the helical curve on the left picture.
The data missing from the points $x\in M\setminus B$ are compensated by measuring
the heat kernel at several times $t_\ell>0$, see Theorem \ref{main uni thm}.
For details of the figure, see the part II of this paper \cite[Section 7]{LLY}.}
\label{fig_embedding}
\end{figure}

Gel'fand's problem \ref{Gelfand:1} is encountered in imaging sciences and inverse problems for partial differential equations, and in dimensionality reduction, especially in the study of diffusion maps \cite{Coifman2,Coifman1} and manifold learning \cite{FMN,ISOMAP}. 
In this paper we are particularly concerned with the situation where the spaces have collapsing structures. 
One such example is
manifold learning,
also called dimensionality reduction in data science \cite{ISOMAP},
where the data sampled from a submanifold of a high-dimensional space need to be 
approximated by a  submanifold of a low-dimensional space.
In manifold learning one is given a point cloud, that is, a finite subset $X=\{x_j\}_{j=1}^N\subset \R^n$ of 
 $N$ points which lie close to a $d$-dimensional submanifold $M_0$ in an $n$-dimensional Euclidean  space, where $n$ is
 much larger than $d$ \cite{FMN}. Then, the goal in learning is to find a manifold $M\subset \R^m$ 
 that is diffeomorphic to $M_0$ so that $d<m<n$. In a more general problem, called the geometric Whitney
 problem \cite{FIKLN}, one is given a possibly discrete metric space $(Y,d_Y)$ and the task is to construct a smooth Riemannian manifold $M$
that is close to $Y$ in the Gromov-Hausdorff sense.
An extensively studied method in manifold learning is the diffusion maps (or the spectral embedding) introduced
 by Coifman and Lafon in \cite{Coifman2,Coifman1}.
In this method, one considers the data points  $X=\{x_j\}_{j=1}^N$, which are sampled from some  manifold $M_0$ and are given with distances $d_X(x_j,x_{j'})$,
and computes a kernel function $h_t:X\times X\to \R$ that approximates the heat kernel of the manifold $M_0$ at a time $t>0$. Then, one
continues the algorithm by computing the first $J$ eigenfunctions $\phi_j(x)$
of the integral operator defined by  the kernel $h_t(x,y),$ $x,y\in X$, and uses the eigenfunction map,
definded as 
\begin{eqnarray}\label{eigenmap}
 & &\Phi^{(K,J)}:X\to \R^J,\\ \nonumber
& &\Phi^{(K,J)}(x)=\big(\phi_K(x),\phi_{K+1}(x),\dots,\phi_{K+J-1}(x) \big),
\end{eqnarray}
to construct the set $\Phi^{(K,J)}(X)\subset \R^J$ that approximates the Riemannian manifold $M_0$, see Figure \ref{fig_embedding} for a toy example. In this paper we consider a similar problem of reconstructing an approximation
of the Riemannian manifold $M_0$ in the case when we have restricted data, that is, one has only local data
$h_{t_\ell}(x_k,x_{k'})=H(x_k,x_{k'},{t_\ell}),$ where the sample points $\{x_k:\ k=1,2,\dots,K\}$ are an $\e$-dense set of points in a possibly small metric ball $B(x_0,r)$ of $M_0$, not an $\e$-dense set of points in the whole manifold $M_0$.
We also define and analyze a local version of the eigenfunction map \eqref{eigenmap}, see Definition \ref{local_data} later. 
The collapsing of manifolds studied in this paper are closely related to multi-scale models in manifold learning where the intrinsic dimension of the data set is modelled by a function that depends on the scale, see \cite{WM}.
For example, this occurs when the cryogenic electron-microscopy \cite{Singer} is applied to image the structure of a large molecule which is connected to a small part of the molecule with a connection that allows a rotation (roughly
speaking, the molecule has a moving `tail'). In this case, the problem of imaging the molecule is to find a manifold diffeomorphic to, e.g.  $SO(3)\times S^1(\epsilon)$ which almost collapses to $SO(3)$ of lower dimension, where $S^1(\epsilon)$ is a circle of small radius $\epsilon$.
Problem \ref{Gelfand:1} for collapsing manifolds is also encountered in mathematical physics, for example in a predecessor of the string theory, the Kaluza-Klein theory, in which Einstein equations are considered on $M_\epsilon=\mathbb{R}^4\times S^1(\epsilon)$, and when
$\epsilon\to 0$, the Einstein equations on $M_\epsilon$ converge to equations on $\mathbb{R}^4$ containing both the standard
Einstein equations and Maxwell's equations. In this case, the observations analogous to the data \eqref{intro-PHD}
give information of the fine structure of the almost collapsed metric.
These applications for collapsing manifolds are considered in the part II of this paper \cite[Section 7]{LLY}.

When the positions of $z_{\alpha}\in B$ are known, the positive answer to Problem \ref{Gelfand:1} was essentially given in \cite{BK,KKL,KrKuLa}, using the boundary control method pioneered by Belishev \cite{Bel} on domains of $\mathbb{R}^n$ and Tataru's unique continuation theorem \cite{Tataru}, see also 
\cite{Alexakis,BKL3,LL}. Generalizations and alternative methods to solve the problem have been studied in \cite{AKKLT,Caday,KKLM,KOP,LO}, and the related determination of smooth structure was recently studied in \cite{FIKLN,FILLN}. On a given domain of the Euclidean space, the problem can be reduced to inverse coefficient problems for elliptic equations which were solved in \cite{SyU}.
We are concerned with the stability of the inverse problem.

\begin{problem}
{\rm Does the solution of the Inverse Problem \ref{Gelfand:1} depend continuously on the given data?}
\end{problem}

Gel'fand's inverse problem is ill-posed in the sense of Hadamard, as one can make large changes to the geometry without affecting the local measurements much. 
Moreover, in the study of inverse coefficient problems in Euclidean spaces, allowing non-smooth coefficients in operators can lead to counterexamples where even the change of topology is not observed in the measurements \cite{GLU}. This phenomenon, in turn, has led to engineering applications in the emerging field of transformation optics and invisibility cloaking \cite{GKLU1} in Calder{\'o}n's inverse problem,
see \cite{AP,KSU,Na2}.
This also makes clear the necessity of studying the question of uniqueness in the limiting non-smooth case, in order to understand the stability of the solution of Gel'fand's inverse problem for smooth manifolds.  

To stabilize the inverse problem, it is natural to impose \emph{a priori} bounds in an invariant form on geometric parameters such as dimension, diameter, sectional curvature and injectivity radius. Within such a class of manifolds, an abstract continuity result for the stability of the inverse problem and a $\log$-$\log$ type of stability estimate were proved in \cite{AKKLT,BKL3,BILL}. 
In this class of manifolds, imposing a lower bound on the injectivity radius is essential as the limiting space is always a Riemannian manifold of the same dimension.
Stronger H\"older type of stability estimates can be obtained in \cite{AG,BD,SU,SU2},
see also \cite{Gui,SUV1,SUV2,UW} for results on closely related tomography problems, if additional geometric assumptions are assumed, e.g., if the metric is close to being simple. 

However, the situation gets much more complicated in general, as it is well-known that a sequence of smooth $n$-dimensional manifolds without a lower bound on the injectivity radius can collapse
to a space of lower dimension and the limiting space is not necessarily a manifold, see e.g. 
\cite{GLP}.
It might be possible that a similar invisibility phenomenon could occur in geometric inverse problems:  the information on microstructures can vanish in the collapse of dimension and one can determine only
some effective properties of the metric. In this paper we study the question of
what global information on collapsing manifolds can be determined from the local measurements, and prove the stability of the Inverse Problem \ref{Gelfand:1} in the class of manifolds with bounded diameter and sectional curvature.

\smallskip
In this paper, we  work with the class ${\frak M}{\frak M}_p={\frak M}{\frak M}_p(n,\Lambda,D)$ 
of connected, closed, smooth, pointed Riemannian manifolds $(M, p, \mu_M)$,  $p\in M$, satisfying
\bequ\label{basic cond 1}
\dim(M)= n, \quad |R(M)|\leq \Lambda^2, \quad \diam(M)\leq D,
\eequ
where $R(M),\diam(M)$ are the sectional curvature and diameter of $M$, and $\mu_M$ is the normalized Riemannian measure \eqref{normalised_measure}.
We also use the notation ${\frak M}{\frak M}(n,\Lambda,D)$ for the same class of manifolds when there is no need to specify a point $p$.

When $(X,d)$  is a metric space,
 we denote by $d_X(x,y)$ the distance between $x$ and $y$ on $X$. When
$A\subset B\subset  X$, we say that $A$ is 
a $\delta$-net in $B$
if for any $y\in B$, there is $x\in A$ such that $d_X(x,y)<\delta$. 
Denote by $A^{\e}$ the $\e$-neighborhood of $A$, that is,
$A^{\e}:= \{ x\in X:\, d_X(x, A) < \e\}$. 
We say that $\psi:X\to X'$ is 
an $\e$-Gromov-Hausdorff approximation (or $\e$-approximation
in short)  if
\beq \label{GH-approximation}
    |d_{X'}(\psi(x), \psi(y)) - d_X( x,y)|<\e, \quad\hbox{and}\quad      \psi(X)\ \hbox{is an $\e$-net in } X'.
\eeq

A metric-measure space $(X,d,\mu)$ is a metric space $(X,d)$
endowed with a Borel measure $\mu$. Below we often
denote  $(X,d,\mu)$ and  $(X',d',\mu')$ just by   $(X,\mu)$ and  $(X',\mu')$, respectively. 
We use the following 
 measured Gromov-Hausdorff\ ``distance''  (which is equivalent to the Prokhorov metric on the space 
of probability measures on a fixed space, see \cite{Bill:conv}).

\begin{definition} \label{mm-distance}
{\rm Let $(X,\mu)$ and  $(X',\mu')$ be metric-measure spaces such that
$X$  and $X'$  are compact. 
The measured Gromov-Hausdorff distance
$d_{mGH}((X,\mu),(X',\mu'))$ 
is defined as the infimum of those 
$\e >0$ such that there are measurable $\e$-approximations
$\psi:X\to X'$ and $\psi':X'\to X$ satisfying
\bequ \label{eq:meas-approx}
      \mu(\psi^{-1}(A')) <  \mu'((A')^{\e}) + \e, \quad
      \mu'((\psi')^{-1}(A)) < \mu(A^{\e}) + \e, 
\eequ
for all Borel sets $A\subset X$ and $A'\subset X'$, 
where $A^{\e}$ is the $\e$-neighborhood of $A$.}

{\rm For pointed compact metric-measure spaces $(X,p,\mu)$ and  $(X',p',\mu')$, the pointed measured Gromov-Hausdorff distance $d_{pmGH}$ is defined as the infimum of those $\e >0$ such that there are measurable $\e$-approximations
$\psi:X\to X'$ and $\psi':X'\to X$ satisfying \eqref{eq:meas-approx} and
\bequ \label{def-pointed}
d_{X'}(\psi(p),p')<\e, \quad d_{X}(\psi'(p'),p)<\e.
\eequ}
\end{definition}

%
%

\smallskip
Let $\overline{{\frak M\frak M}}_p(n,\Lambda, D)$ denote the completion of 
the set ${\frak M\frak M}_p(n,\Lambda, D)$, that is, the space of equivalence classes of 
the Cauchy sequences, with respect to the pointed measured Gromov-Hausdorff  distance
$d_{pmGH}.$ Usually we say that  $\overline{{\frak M\frak M}}_p(n,\Lambda, D)$ is
the closure of  ${{\frak M\frak M}}_p(n,\Lambda, D)$. 

It is well-known that a sequence of $n$-dimensional manifolds in the class ${\frak M\frak M}_p(n,\Lambda, D)$ can collapse to a lower dimensional space when the injectivity radius of the sequence of manifolds goes to zero.
The structure of the metric-measure spaces
 $(X,\mu_X) \in \overline{{\frak M}\frak M}_p(n,\Lambda,D)$ was studied in \cite {Fuk_JDG2}--\cite {Fuk:haus}.
The space $X$ has the stratification 
\bequ
     X = S_0(X) \supset S_1(X) \supset \cdots \supset S_d(X),
       \label{eq:stratif}
\eequ
with the following property: if $S_j(X)\setminus S_{j+1}(X)$ is non-empty, then it is a 
$(d-j)$-dimensional Riemannian manifold of class $C^{1,\alpha}$ for any $0<\alpha<1$, where $d=\dim(X)\leq n$. The regular part of $X$ is the set
$X^{reg}:=S_0(X) \setminus S_1(X)$ which is an open $d$-dimensional manifold of class $C^{1,\alpha}$, and the singular set is its complement, $X^{sing}:=X\setminus X^{reg}$. In particular, the singular set $X^{sing}$ has dimension at most $d-1$. In this paper we improve on the regularity of the limiting spaces. We show that on $X^{reg}$, the metric is
locally determined by $C^2_*(X^{reg})$-smooth Riemannian tensor $h_X$ and the measure $\mu_X$
is absolutely continuous with respect to the Riemannian measure $dV_{h_X}$
with Radon-Nikodym derivative $\rho_X\in C^2_*(X^{reg})$. Here,
$C^2_*(X^{reg})$ is the Zygmund space having the relation $C^{1,1}(X^{reg})\subset C^2_*(X^{reg})\subset
C^{1,\alpha}(X^{reg})$ for any $0<\alpha<1$. Moreover,  $X^{reg}$ is convex and $\mu_X(X^{sing})=0$. 
We remark that locally on the regular part $X^{reg}$, the $C^2_*$-regularity provides the local uniqueness and stability for the geodesic equation in the Riemannian sense (e.g. \cite{AKKLT},\cite[Prop. A.2]{PDE-Taylor}), which does not hold in general for $C^{1,\alpha}$ metrics (\cite{H50}).
More details on the structure of collapsing are reviewed in Section \ref{sec:review}.

As shown in \cite {Fuk_inv}, the Dirichlet's quadratic form $A[u,v]=\bra du,dv\cet _{L^2(X^{reg},d\mu_X)}^2$,
$u,v\in C^{0,1}(X)$,  defines
a self-adjoint operator $\Delta_X$ on $(X,\mu_X)$, which we call the weighted Laplacian, denoted by $\Delta_X$. In local coordinates
  on $X^{reg}$, it has the form
\bequ
 \Delta_X u= -\frac{1}{\rho_X |h_X|^{\frac12}} \sum_{j,k=1}^d \frac \p{\p x^j} \Big(\rho_X |h_X|^{\frac12} h^{jk}_X \frac \p{\p x^k}  u  \Big), \quad |h_X|=\det\big((h_X)_{jk} \big).
\eequ
The associated semigroup $e^{-t\Delta_X}$  has the Schwartz kernel $H_X(x,y,t)$ that we call
the heat kernel associated to $(X,\mu_X)$. Denote by $\lambda_j$ the $j$-th eigenvalue of the weighted Laplacian $\Delta_X$ and by $\phi_j$ the corresponding orthonormalized eigenfunction in $L^2(X,\mu_X)$.
Fukaya proved in \cite{Fuk_inv} that the $j$-th eigenvalue, for any $j$, is a continuous function on  $\overline{{\frak M}\frak M}_p(n,\Lambda,D)$ with respect to the measured Gromov-Hausdorff topology.

\subsection{Main results}

Our first main result is the following uniqueness theorem for the inverse problem for collapsing manifolds.

\begin{theorem} \label{main uni thm}
Let $r,\Lambda,D >0, \,n \in \Z_+$.
Denote by $\,\overline{\frak M \frak M}_p(n,\Lambda, D)$ the closure of the class of connected closed smooth pointed Riemannian manifolds
defined by the conditions (\ref{basic cond 1}) in the pointed measured Gromov-Hausdorff topology.
Let $(X,p,\mu),(X',p',\mu')\in \overline{{\frak M\frak M}}_p(n,\Lambda, D)$.
Let $\{z_\a \}_{\a=0}^{\infty} \subset X,\,\,
\{z'_\a \}_{\a=0}^{\infty} \subset X'$
be dense sequences in the ball $B_X(p, r)$, $B_{X'}(p', r)$ with $z_0=p$, $z'_0=p'$, 
and $\{t_{\ell}\}_{\ell=1}^\infty$ be a dense set in $\mathbb{R}_+$.
Suppose that
\bequ \label{eq: uniqueness 1}
H(z_\a, z_\beta, t_{\ell})=H'(z'_\a, z'_\beta, t_{\ell}), \quad \hbox{for all } \a, \beta\in \N, \; \ell\in \Z_+,
\eequ
where $H, \, H'$ are the heat kernels on $X,\, X'$.
Then there exists an isometry $F:X\to X'$ satisfying $F(p)=p'$ and $\mu=F^*\mu'$. 
\end{theorem} 

In particular, the spaces $X$ and $X'$ necessarily have the same dimension.
Note that the assumption on the positions of $p,p'$ relative to the dense sequences $\{z_{\alpha}\},\{z'_{\alpha}\}$ is necessary: otherwise, it is not possible to determine the location of $p,p'$ when $r>{\rm diam}(X)$.
Theorem \ref{main uni thm} implies our second main result: the stability theorem for the Inverse Problem \ref{Gelfand:1}.

\begin{theorem} \label{main stability thm}
Let $r,\Lambda,D >0, \,n \in \Z_+$.
Then, there exists an increasing, continuous function $\omega(s)=\omega_{(r, n,\Lambda, D)}(s)$ that depends on $r, n,\Lambda, D$,
which maps $\omega: [0, 1) \to
[0, \infty)$ with $\omega(0)=0$,
such that the following holds.

\smallskip
\noindent Let $(X, p, \mu),\, (X', p', \mu') \in \overline{{\frak M \frak M}}_p (n,\Lambda, D)$.
Let
$\{z_\a \}_{\a=0}^{N} \subset X,\,\,
\{z'_\a \}_{\a=0}^{N} \subset X'$
be $\delta$-nets in the ball $B_X(p, r), \,B_{X'}(p', r)$ with $z_0=p$, $z'_0=p'$, and 
$\{ t_{\ell} \}_{\ell=1}^{L}$ be a 
$\delta$-net in $(\delta,\, \delta^{-1})\subset \R_+ $, where $N,L\in \Z_+$. 
Suppose that
\bequ \label{000.3}
|H(z_\a, z_\beta, t_{\ell})-H'(z'_\a, z'_\beta, t_{\ell})| < \delta, \quad \hbox{for }\, 0 \leq \a, \beta \leq N,\,\, 
1 \leq \ell \leq L,
\eequ
where $H, \, H'$ are the heat kernels on $X,\, X'$.
Then
\bequ \label{2.28.07.10}
d_{pmGH}\left((X, p, \mu),\, (X', p', \mu')  \right) < \omega(\delta).
\eequ
\end{theorem} 

In particular, in the case of the spaces $X, X'$ being $n$-dimensional manifolds, Theorem \ref{main stability thm} removes the injectivity radius assumption for the stability result in \cite{AKKLT}.

\smallskip
If $(X,p,\mu)\in \overline{{\frak M}{\frak M}}_p(n,\Lambda,D)$ is \emph{a priori} an orbifold, i.e., $\textrm{dim}(X)=n-1$, it is possible to obtain an explicit form of the modulus of continuity $\omega$, given the interior spectral data $\{\lambda_j,\phi_j|_{B_X(p,r)}\}$ for the weighted Laplacian $\Delta_X$ on a ball $B_X(p,r) \subset X^{reg}$, provided that the metric tensor on the regular part $X^{reg}$ is of class $C^4$.
To this end, let us denote by ${\frak M}{\frak M}_p(n,\Lambda,D,\Lambda_3)$ the class of connected closed smooth pointed Riemannian manifolds $(M,p,\mu_M)\in {\frak M}{\frak M}_p(n,\Lambda,D)$ satisfying additionally $\|\nabla^i R(M)\|\leq \Lambda_3$ for $i=1,2,3$, and denote its closure by $\overline{{\frak M}{\frak M}}_p(n,\Lambda,D,\Lambda_3)$.

\begin{theorem} \label{thm-orbifold}
Let $(X,p,\mu) \in \overline{{\frak M}{\frak M}}_p(n,\Lambda,D,\Lambda_3)$ with $p\in X^{reg}$. Suppose $\dim(X)=n-1$ and ${\rm Vol}_{n-1}(X) \geq v_0$.
Let $\sigma\in (0,1)$, and $r>0$ such that the ball $B_X(p,r)\subset X^{reg}$.
Then there exists $\widehat{\delta}=\widehat{\delta}(X,r,\sigma)>0$, such that
the finite interior spectral data $\{\lambda_j,\phi_j|_{B_X(p,r)}\}_{j=0}^{\delta^{-1}}$ for $\delta<\widehat{\delta}$ determine a finite metric space $\widehat{X}$ such that
$$d_{GH}(X\setminus S_{\sigma,\delta},\widehat{X}) < C_1(X,\sigma) \Big(\log\log|\log \delta| \Big)^{-C_2},$$
where $S_{\sigma,\delta}$ is a subset 
of the $C_3\sigma^{1/4}$-neighborhood of the singular set of $X$, and $X\setminus S_{\sigma,\delta}$ is equipped with the restriction of the metric of $X$.
The constant $C_1(X,\sigma)$ depends on $X,\sigma,r,n,\Lambda,\Lambda_3,D,v_0$, and $C_1(X,\sigma)\to \infty$ as $\sigma\to 0$. The constant $C_2$ depends on $n$, and $C_3$ depends on $n,\Lambda,D$.
\end{theorem}

Theorem \ref{thm-orbifold} will be proved in the part II of this paper \cite{LLY}. 

\subsection{Plan of the exposition} \label{Plan}

This paper is organized as follows.
Section \ref{sec:review} is of an expository nature, where we review, in a slightly modified
form appropriate for our purposes, Fukaya's results on the measured Gromov-Hausdorff convergence of Riemannian manifolds and provide some further results in this direction. 
In Section \ref{smooth}, we show that the density
function $\rho_X$ for $(X,\mu_X) \in \overline{\frak M \frak M}_p(n,\Lambda,D)$ is of $C^2_*(X^{reg})$. This improves on the earlier results in \cite{Fuk_inv,Kas:measII,K93}. Note that our proof differs from that in the above papers as it is based on
the analysis of smoothness of the transformation groups 
into the scale of Zygmund-type functions.
In turn, this requires an extension of the classical
Montgomery-Zippin results, which we will discuss in the part II of this paper \cite{LLY}.
Section \ref{section: Aux results} is of an auxiliary nature. Here we prove various results
concerning the behavior of the spectrum, eigenfunctions and heat kernels on 
${\overline {\frak M \frak M}}_p(n,\Lambda,D)$, and analyze the spectral information contained
in the heat kernel data. In Section \ref{Continuity of the direct map}, we continue to study the spectral
behavior on ${\overline {\frak M \frak M}}_p(n,\Lambda,D)$, obtaining uniform estimates
for the eigenfunctions and heat kernels, and prove the pointwise convergence
of the heat kernels with respect to the pointed measured
Gromov-Hausdorff convergence. Some auxiliary results dealing with the relations between
the Laplacian on manifolds and their orthonormal frame bundles, as well as the corresponding
structures on ${\overline {\frak M \frak M}}_p (n,\Lambda,D)$, are discussed in Appendix \ref{OTA}.
In Section \ref{metric}, by extending the geometric boundary control method to the collapsing manifolds, we show that the heat kernel data for any 
$(X,\mu_X) \in {\overline {\frak M \frak M}}_p (n,\Lambda,D)$
uniquely determine its
metric-measure structure (Theorem \ref{main uni thm}). At last, Section \ref{Stability of inverse problem} is devoted to
the proof of the stability result, Theorem \ref{main stability thm}.

\medskip
\noindent {\bf Acknowledgements.}
 M.L. and J.L. were supported by PDE-Inverse project of the European Research Council of the European Union,
 the FAME-flagship and the grant 336786 of the Research Council of Finland. Views and opinions expressed are
those of the authors only and do not necessarily reflect those of the European
Union or the other funding organizations.  
T.Y. was supported by JSPS KAKENHI Grant Number 21H00977.

\section{Basic results in the theory of  collapsing} \label{sec:review}

Given $n\in \mathbb{N}_+$ and $\Lambda, D\in \mathbb{R}_+$, the notation
${\frak M}(n,\Lambda, D)$ stands for the class  of $n$-dimensional connected, closed, smooth Riemannian manifolds $(M,h)$ satisfying 
\[
      \dim(M)= n, \quad |R(M)|\leq \Lambda^2, \quad \diam(M)\leq D,
\]
where  $R(M)$ and $\hbox{diam}(M)$ are the sectional 
curvature and diameter of $M$.
The structure of collapsing in the moduli space ${\frak M}(n,\Lambda, D)$, with respect to the Gromov-Hausdorff distance,
was extensively studied by Fukaya (\cite{Fuk_inv, Fuk_JDG, Fuk:haus}).   
In this section we review some of the results.

\subsection{Properties of the limit spaces} \label{subsection-basic-property}

\noindent\textbf{2.1.1 \, Basic example.} It is known (e.g. \cite{GLP}) that ${\frak M}(n, \Lambda, D)$ is precompact in the Gromov-Hausdorff topology.
Therefore we assume that a sequence $(M_i,h_i) \in {\frak M}(n, \Lambda, D)$, $i=1,2, \dots$,
converges to some compact metric space $X$ with respect to the 
Gromov-Hausdorff distance. 
The space $X\in \overline{{\frak M}}(n, \Lambda, D)$ can be locally described as orbit spaces.
It has the stratification 
\bequ
     X = S_0(X) \supset S_1(X) \supset \cdots \supset S_d(X),
       \label{eq:stratif}
\eequ
such that, if $S_j(X)\setminus S_{j+1}(X)$ is non-empty, it is a 
$(d-j)$-dimensional Riemannian manifold of class $C^{1,\alpha}$ for any $0<\alpha<1$,
where $d:= \dim (X) \leq n$. 
Actually $S_j(X)\setminus S_{j+1}(X)$ is defined as the set of all point $p\in X$ such
that the tangent cone $K_p(X)$ of $X$ at $p$ (see e.g. \cite{BBI,BGP}) is isometric to a product of the form
$\Bbb R^{d-j}\times Y^{j}$, where $Y$ has no such a
nontrivial $\Bbb R$-factor.

We denote by 
\bequ \label{def-sin-reg}
 X^{sing}:= S_1(X),\quad X^{reg}:= X\setminus S(X)
\eequ
the singular and the regular set of $X$, respectively.
By definition, the regular set $X^{reg}$ is an open $d$-dimensional Riemannian manifold with $d=\textrm{dim}(X)$, and the singular set $X^{sing}$ has dimension at most $d-1$.
When $\dim (X) \le n-1$, 
we say that 
$M_i$ collapses to $X$,
and that the collapsing is $m$-dimensional for  $m=n-\dim (X)$.

\smallskip
Let us briefly explain the structure of $1$-dimensional collapse in the following example.

\noindent
\begin{minipage}[t]{0.73\textwidth}

\begin{example} \label{orbif_collapse}{\rm
 For $\epsilon>0$, let $M_\epsilon$ be a $3$-dimensional Riemannian manifold of the following form.
 We start with a product manifold $\mathbb{S}^2 \times [0, \epsilon]$, where $\mathbb{S}^2\subset \R^3$ is the $2$-dimensional unit sphere 
 with the canonical metric.  Consider the action of the finite group $\Z_m, \,
 m \geq 2$, on $\mathbb{S}^2$ by a rotation of angle $2 \pi/m$ along the $z$-axis. Then we define $M_{\epsilon}$ by identifying points 
 $({\bf x}, 0) \in \mathbb{S}^2 \times \{0\}$ with $(e^{2\pi i/m} \cdot {\bf x}, \epsilon) \in \mathbb{S}^2 \times \{\epsilon\}$ in the manifold $\mathbb{S}^2 \times [0, \epsilon]$,
 where $e^{2\pi i/m} \cdot {\bf x}$ stands for the rotation by the angle $2 \pi/m$.
 Observe that this gives rise to closed vertical geodesics of length $m \epsilon$ except for the points corresponding to the north and south poles.
 When $\epsilon \to 0$, $M_\epsilon$ collapses to a $2$-dimensional space $X=\mathbb{S}^2/\Z_m$, a Riemannian orbifold as illustrated on the right.
 This orbifold $X$ has conic singular points only at the north and south poles. For analytically-oriented readers, we cite to \cite[Section 7.2]{LLY} for the details and additional pictures on this example.}
 \end{example} 
\end{minipage}
\begin{minipage}[t]{0.27\textwidth}

\centering
\raisebox{\dimexpr \topskip-\height}{%
  \includegraphics[width=0.9\textwidth]{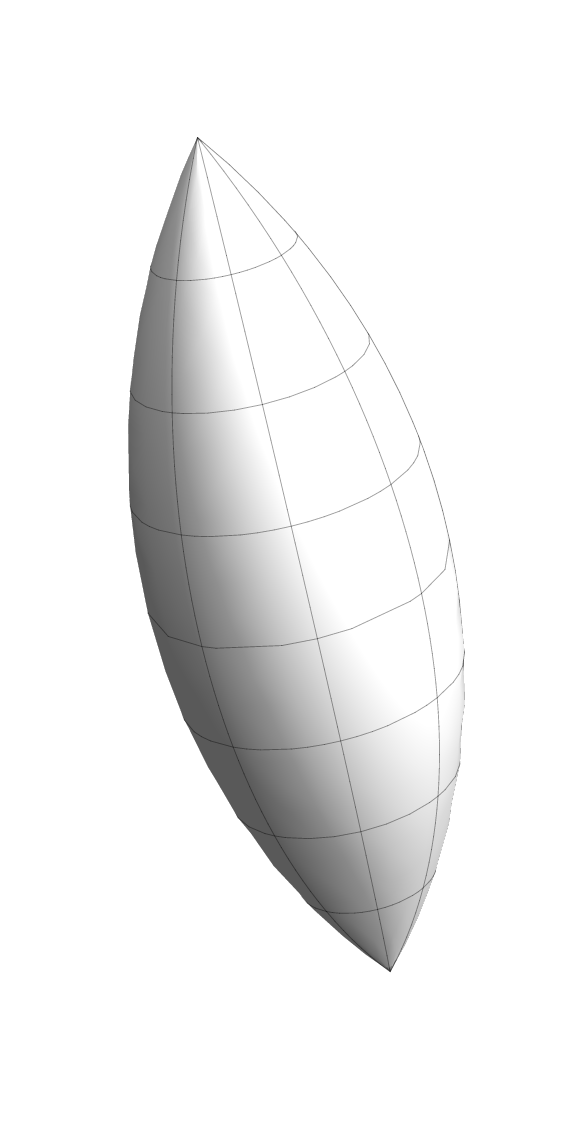}}
\label{fig_football}

\end{minipage}
 \smallskip

 An example of the limit spaces above is so-called good orbifolds (e.g. \cite{Thurston}), that is, 
$$X=N / G,$$
where $N$ is a Riemannian manifold and $G$ is a discrete group acting properly and isometrically on $N$.
For another example, one could consider $N=\mathbb{S}^2\times \mathbb{S}^1$, and $\mathbb{Z}_m$ is the finite group acting on $\mathbb{S}^2\subset \mathbb{R}^3$ by a $2\pi/m$-rotation along the $z$-axis as in Example \ref{orbif_collapse}.
Then the quotient space $N/\mathbb{Z}_m=(\mathbb{S}^2/\Z_m) \times \mathbb{S}^1$ is a good orbifold, where the singular set of the orbifold is $\{\textrm{north pole, south pole}\} \times \mathbb{S}^1$ and the complement is the regular set.

\medskip
\noindent\textbf{2.1.2 \, Local pseudogroup construction.} 
In this subsection, we consider the geometric properties of the collapsing manifolds and the properties of the limit space $X$.
An analytically-oriented reader may consider the orbifold space described in Example \ref{orbif_collapse} as a basic example of the limit of collapsing manifolds, and skip in the first reading the general local construction in this subsection.

Let $(M_i,h_i)\in {\frak M}(n,\Lambda, D)$ be a sequence of closed Riemannian manifolds converging to a compact metric space $X$ in the Gromov-Hausdorff topology.
Fix any point $p\in X$ and put  $p_i:= \psi_i(p)$, where $\psi_i:X\to M_i$ is
an $\e_i$-approximation as defined in \eqref{GH-approximation} with $\lim\limits_{i\to\infty} \e_i = 0$.
Let $B \subset \R^n$ be the open ball around the origin $\it O$ in $\R^n$ of radius $\pi/\Lambda$, and
let $\exp_i: B\to M_i$ be the composition of the exponential map $\exp_{p_i}:
 T_{p_i}(M_i)\to M_i$ and a li isometric embedding 
$B\to B({\it O},\pi/\Lambda)\subset T_{p_i}(M_i)$.
Since $|R(M_i)| \leq \Lambda^2$, the exponential map on $B(O,\pi/\Lambda)$ has maximal rank, and 
we have the pull-back metric $\tilde h_i := \exp_{i}^*(h_{i})$ on $B\subset \R^n$.
Moreover, the injectivity radius of $(B,\tilde h_i)$ is uniformly bounded from below, e.g. \cite[Lemma 8.19]{GLP}.
Therefore, we may assume that $(B, \tilde h_i)$ converges to a
$C^2_{*}$-metric $(B, \tilde h_0)$ in the $C^{1,\alpha}$-topology, 
for any $0<\alpha < 1$ (see \cite{A90,AKKLT,GW:lipschitz,Pets:conv}).

Let $G_i$ denote the set of all  isometric 
embeddings $\gamma:(B', \tilde h_i) \to (B, \tilde h_i)$ such that
$\exp_i \circ \gamma = \exp_i$ on $B'$, where $B':= B({\it O}, \frac{\pi}{2\Lambda})\subset\R^n$. 
Then $G_i$ forms a local pseudogroup, see e.g. \cite[Section 7]{Fuk:haus}.
Passing to a subsequence, we may assume that
$G_i$ converges to a local pseudogroup   $G$ consisting of isometric 
embeddings $\gamma:(B', \tilde h_0) \to (B, \tilde h_0)$, and that
$(B,\tilde h_i, G_i)$ converges to $(B, \tilde h_0,G)$ in the 
equivariant Gromov-Hausdorff topology.
 This means that there exist $\e_i$-approximations 
 $ \tilde \phi_i: (B, \tilde h_i) \to (B, \tilde h_0)$,  
$\tilde \psi_i:(B, \tilde h_0) \to (B, \tilde h_i)$
with $\lim\limits_{i\to\infty} \e_i = 0$, and  maps
$ \rho_i: G_i \to G$, $\sigma_i:  G \to G_i$
such that for every $x, y \in B$ and $\gamma_i \in G_i, \,\gamma \in
 G,$
the following hold
 \bequ \label{29.06.3}
d_0\big(\tilde\phi_i(\gamma_i (x)), \rho_i(\gamma_i)(\tilde\phi_i(x)) \big) < \e_i,   \quad
  d_i \big(\tilde\psi_i(\gamma(x)), \sigma_i(\gamma)(\tilde\psi_i(x))\big) <  \e_i,
 \eequ
whenever they make sense, where $d_0,\, d_i$ denote the distance function on $(B,\tilde h_0),\, (B,\tilde h_i)$, respectively.
 Roughly speaking, this means that the pseudogroup action of $G_i$ on $(B,\tilde h_i)$ 
is close to that of $G$ on $(B, \tilde h_0)$.
In particular, the quotient space $(B',\tilde h_i)/G_i=B(p_i, \frac{\pi}{2\Lambda})\subset M_i$ 
converges to $(B',\tilde h_0)/G$,
which implies that 
 \bequ \label{eq:limit-char}
     (B', \tilde h_0)/G=B(p, \frac{\pi}{2\Lambda}) \subset X.
 \eequ
 (See \cite{FY:fundgp} for further details on
basic properties of the equivariant Gromov-Hausdorff convergence.)
%

\medskip
\noindent\textbf{2.1.3 \, Geometric structure of the limit space.}
The following lemma is known, where (1), (2) are due to Fukaya
\cite{Fuk_JDG}
(and \cite{AKKLT}), and (3) is the consequence of \cite{Pet}.

\begin{lemma} \label{injectivity radius}
Let $X$ be any element of $\overline{{\frak M}}(n,\Lambda, D)$ which is not a point.
Then
\begin{enumerate}
  \item[$(1)$] $X^{reg}$ is a $C^2_*$-Riemannian manifold, that is, there exist coordinate charts on $X^{reg}$ for which the transition maps are of class $C^3_*$, and the metric tensor in these charts is of class $C^2_*$.
  \item[$(2)$] For any compact subset $K \subset X^{reg}$,  there exists a positive number 
    $i_K>0$ such that ${\rm inj}(p) \geq i_K$ for all $p \in K$, 
     where ${\rm inj}(p)$ denotes the injectivity radius at $p$.
  \item[$(3)$] $X^{reg}$ is convex in $X$, that is, every minimizing geodesic joining two points
    in $X^{reg}$ is contained in $X^{reg}$.
\end{enumerate}
\end{lemma}

For reader's convenience, we shall sketch the outline of the proof of
Lemma \ref{injectivity radius}(1) and (2) below.
Let $X$ be the limit of $(M_i, h_i) \in {\frak M}(n,\Lambda, D)$.
By \cite{Fuk_JDG}, there exist Riemannian metrics  $h_i^{\e}$
on $M_i$ such that,  as $\e \to 0$, $h_i^{\e} \to h_i$ in the $C^{1,\alpha}$-topology, for any 
$0<\alpha<1$,  and 

(i) $(M_i, h_i^{\e})$ converges  to $X^\e$, as $i\to\infty$,
    with respect to the Gromov-Hausdorff distance;
    
(ii) the regular part $(X^{\e})^{reg}$ of $X^{\e}$ is a
   Riemannian manifold of class $C^{\infty}$;
   
(iii) $X^{\e}$ is $\e$-isometric to $X$. Namely, there
   exists a  bi-Lipschitz map  $f^{\e}:X \to X^{\e}$ satisfying
      \[
          \left| \frac{d(f^{\e}(x), f^{\e}(y) )}{d(x,y)} -1\right| < \e.
       \]
More precisely, the norm of the $k$-th covariant derivaitives 
of the curvature tensor $R_{h_i^{\e}}$ of 
$(M_i,h_i^{\e})$  is uniformly bounded $ \| \nabla^k R_{h_i^{\e}} \| \le
C(n,\Lambda,k,\e)$ for any fixed $k$ and $\e$ (see \cite{BMR}).
The space $X^{\e}$ here is 
called a {\it smooth element} in \cite{Fuk_JDG}.

Since $X$ and $X^{\e}$ have orbit-type singularities,  the above (iii) implies that 
$f^{\e}(S_j(X))= S_j(X^{\e})$ for small $\e$, where 
$S_j(X)$ are the stratification of $X$ in \eqref{eq:stratif}.
In particular, $f^{\e}(X^{reg})=  (X^{\e})^{reg}$.

For a small $\delta>0$, let $V$  be the $\delta$-neighborhood of $K$  
in $X$ such that $\overline V\subset X^{reg}$. Set $V^{\e}:= f^{\e}(V)$.
Then,
\bequ \label{04.11.01}
-\Lambda^2 \leq R\left(X^\e\right)|_{V^{\e}}  \leq C(n,\Lambda),
\eequ
see \cite[Theorem 0.9]{Fuk_JDG}. 
Note that, by (iii), $\hbox{diam}(V^\e) \leq C$, $\hbox{Vol}(V^\e) \geq C$ for
some 
uniform constant $C=C(K)>0$, which is independent of $\e$. Thus, 
(\ref{04.11.01}) together with  Cheeger's theorem \cite{Ch:finite}
implies that  there is a positive number $i_K<\delta$ independent of $\e$
such that 
\bequ
       \hbox{inj}_{X^\e}(p) \ge  i_K,   \label{eq:inj1}
\eequ
for all $p \in K^{\e}:=f^{\e}(K)$. 
Now a standard argument using the Cheeger-Gromov compactness applied to
the convegence $V^{\e}\to V$  together with Alexandrov geometry
proves (2).
It follows from (\ref{04.11.01}), \eqref{eq:inj1} and \cite{AKKLT} that
the metric of $V$ and hence of $X^{reg}$ is of class  $C^2_*$.
Then the $C^3_*$-smoothness of the transition maps between harmonic coordinate charts is due to \cite{A90,AKKLT}.

Let us recall here the definition of the Zygmund spaces, e.g. \cite[Chapter 2.7]{Tri}. The Zygmund space $C^s_*(B)$ on a ball $B\subset \mathbb{R}^n$ coincides with the H\"older space $C^s(B)$ for $s\in \mathbb{R}_+\setminus \Z_+$, and forms a complex interpolation scale (for an overview on the interpolation theory of function spaces, see \cite{BL}).
For $s\in \Z_+$, $C^s_*(B)$ is defined as the space of functions in $C^{s-1}(B)$ satisfying that the following norm is finite,
\begin{equation}\label{def-Zygmund}
\|f\|_{C^s_*(B)}:=\|f\|_{C^{s-1}(B)}+\sum_{|\alpha|=s-1}\sup_{\substack{x\neq y \\ x,y\in B}} \frac{\big|D^{\alpha}f(x)-2D^{\alpha}f(\frac{x+y}{2})+D^{\alpha}f(y) \big|}{|x-y|}\, , \;\; s\in \Z_+.
\end{equation}
We point out that when $k$ is a positive integer and $0<\alpha<1$, the H\"older spaces $C^{k,\alpha }(B)$ and
the Zygmund spaces $ C^k_*(B)$ have the relation 
$C^{k-1,1}(B)\subset C^k_*(B)\subset C^{k-1,\alpha }(B).$


\smallskip
The geometric structure of the limit space X can be described in the following ways.

\begin{theorem}[\cite{Fuk_JDG}]  \label{thm:limitstr}
For every $p\in X$, let $G$ be local pseudogroup defined as above,
and set    $\ell := n - \dim (G\cdot {\it O})$, where $G\cdot {\it O}$ denotes the 
orbit $G({\it O})$.
\begin{enumerate}
 \item[$(1)$] There exist a neighborhood $U$ of $p$,
    a compact Lie group $G_p$ and a faithful representation 
    of $G_p$
    into the orthogonal group $O(\ell)$, 
    a $G_p$-invariant smooth metric on a neighborhood $V$ of
    ${\it O}$ in $\R^{\ell}$ such that $U$ is bi-Lipschitz homeomorphic to  $V/G_p;$
 \item[$(2)$] There exists a $C^2_*$-Riemannian manifold $Y$ with
    $\dim (Y) = \dim (X)+ \dim (O(n))$  on which 
    $O(n)$ acts as isometries in such a way that
   \begin{enumerate}
    \item $X$ is isometric to $Y/O(n)$. Let $\pi:Y\to X$ be the 
     projection$;$
    \item For every $p\in X$ and $\bar p\in\pi^{-1}(p)$, the isotropy  group 
       \begin{equation*}
               H_{\bar p} := \{ g\in O(n)\,|\, g(\bar p) = \bar p \}
        \end{equation*}
     is isomorphic to $G_p$, where $G_p$ is as
       in $(1)$.
   \end{enumerate}
\end{enumerate}
\end{theorem}

We will make an extensive use of Theorem \ref{thm:limitstr}(2) in this paper, and the fact of $Y$ being a Riemannian manifold is crucial for our purposes.
Here we briefly explain the construction of $Y$.
Let $X$ be the limit of $M_i \in {\frak M}(n,\Lambda, D)$.
Let $FM_i$ denote the orthonormal frame bundle of $M_i$ endowed with
the natural Riemannian metric, which  has uniformly bounded  sectional 
curvature and diameter.
Note  that $O(n)$ isometrically acts  on $FM_i$, and $M_i=FM_i/O(n)$.
Passing to a subsequence, we may assume that 
$(FM_i, O(n))$ converges to $(Y, O(n))$ in the
equivariant  GH-topology. 
Then it follows that $X=Y/O(n)$. 

To show that $Y$ is a Riemannian manifold,
let $B'\subset B\subset \R^n$, $(B, \tilde h_i,G_i)$ and  $(B, \tilde h_0,G)$ be
as described at the beginning of this section, so that
$(B', \tilde h_i)/G_i=B(p_i,\frac{\pi}{2\Lambda})$ and $(B', \tilde h_0)/G=B(p, \frac{\pi}{2\Lambda})$.
The pseudogroup action of $G_i$ on $(B', \tilde h_i)$ induces an isometric pseudogroup
action, denoted by $\hat G_i$, on the frame bundle $F(B', \tilde h_i)$ of
$(B', \tilde h_i)$ defined by differential. Therefore, 
$F(B', \tilde h_i)/\hat G_i = FB(p_i, \frac{\pi}{2\Lambda})$. Passing to a subsequence,
we may assume that  $(F(B', \tilde h_i), \hat G_i)$ 
converges to $(F(B', \tilde h_0), \hat G)$ in the equivariant
GH-topology, where $\hat G$ denotes the
isometric pseudogroup action on $F(B',\tilde h_0)$ induced from that of 
$G$ on $(B', \tilde h_0)$. 
Since the action of $G$ on $(B',\tilde h_0)$ is isometric, the action of $\hat G$ on  $F(B', \tilde h_0)$
is  free. 
Therefore, $F(B',\tilde h_0)/\hat G$ is a Riemannian manifold, and
so is $Y$.


\begin{theorem}[\cite{Fuk_JDG}, Theorem 10.1]\label{thm:fib}
Suppose  a sequence $X_i$ in ${\overline{\frak M}}(n,\Lambda, D)$ converges to 
$X$ with respect to the Gromov-Hausdorff distance.
Then there are $O(n)$-Riemannian manifolds  $Y_i$ and $Y$ of class $C^{2}_*$ and 
$O(n)$-maps $\tilde f_i:Y_i \to Y$ and maps $f_i:X_i\to X$ such that
\begin{enumerate}
 \item[$(1)$] $X_i=Y_i/O(n)$, $X=Y/O(n)$. Let $\pi_i:Y_i\to X_i$, $\pi:Y\to X$ 
   be the projections;
  \item[$(2)$] $\tilde f_i$ are $\e_i$-Riemannian submersions as well as 
   $\e_i$-approximations, where $\lim\limits_{i\to\infty} \e_i = 0$. Namely, $\tilde f_i$ satisfies
   \bequ \label{1.22.01.2012}
       e^{-\e_i} < \frac{|d\tilde f_i(\xi)|}{|\xi|} < e^{\e_i},
   \eequ
    for all tangent vectors $\xi$ orthogonal to fibers of $\tilde f_i;$
  \item[$(3)$] $f_i\circ \pi_i = \pi\circ \tilde f_i;$
  \item[$(4)$] for every $y\in Y$, the isotropy subgroup 
      $\{ g\in O(n)\,|\, g(y)=y\}$ is isomorphic to $G_{\pi(y)}$, where
     $G_{\pi(y)}$ is as in Theorem \ref{thm:limitstr}.
\end{enumerate}
\end{theorem}

We shall call the maps $f_i:X_i\to X$ {\it regular $\e_i$-approximations}
for simplicity.

\smallskip
Note that $p$ is an orbifold point if and only if 
$G_p\simeq H_{\bar p}$ is finite.
This actually occurs for every $p\in X$ in the case of collapsing being
one dimensional.

\begin{corollary} [\cite{Fuk:haus}, Proposition 11.5] \label{orbifold_collapse}
If $\dim X= n-1$, then 
$X$ is an orbifold.
\end{corollary}


\subsection{Measured Gromov-Hausdorff distance}
Let $\mu_i$ and $\mu$ be probability Borel measures on 
compact metric spaces $X_i$ and $X$.
Fukaya defined the notion of measured Gromov-Hausdorff convergence
$(X_i,\mu_i)\to (X,\mu)$, or weak convergence in short.
By definition, this is the case when  there are measurable $\e_i$-approximations
$\psi_i:X_i\to X$ with $\lim\limits_{i\to\infty} \e_i = 0$
such that the pushfoward measure
$(\psi_i)_{*}\mu_i$ weakly converges to $\mu$ in the usual sense:
namely, 
\bequ \label{measured-def}
\int_{X_i} f\circ\psi_i\, d\mu_i \to \int_X f\,d\mu,\quad \textrm{as }\, i\to\infty,
\eequ
for any $f\in C(X)$, 
where $C(X)$ denotes the space of 
continuous functions on $X$.
In this subsection, we provide basic properties of the 
 measured Gromov-Hausdorff distance 
$d_{mGH}((X,\mu),(X',\mu'))$ defined in Definition \ref{mm-distance}
between compact metric measure spaces.



\begin{lemma} \label{triangle}
The measured Gromov-Hausdorff ``distance'' satisfies an almost triangle inequality:
\ba
    && d_{mGH}((X_1,\mu_1),  (X_2, \mu_2))  \\
       && \hspace{5mm} \le 
       2 \left( d_{mGH}((X_1,\mu_1), (X_3, \mu_3)) + d_{mGH}((X_3,\mu_3), (X_2, \mu_2))\right).
\ea
\end{lemma}

\begin{proof}
Let $d_{ij}:= d_{mGH}((X_i,\mu_i), (X_j, \mu_j))$, and $d := d_{13}+d_{32}$.
By definition, for any $\e>0$ and for any $i,j\in \{ 1, 3\}$ or $i,j\in \{ 3,2\}$,
there are $(d_{ij}+\e/2)$-approximation
$\psi_{ij} : X_i\to X_j$, satisfying 
\[
 \mu_i(\psi_{ij}^{-1}(A_j)) <  \mu_j(A_j^{d_{ij}+\e/2}) + d_{ij}+\e/2, 
\]
for any Borel set $A_j\subset X_j$.
Define $\psi_{12}:X_1\to X_2$ by $\psi_{12}:=\psi_{32}\circ\psi_{13}$, which is a 
$2(d +\e)$-approximation.
Then for any Borel set $A_2\subset X_2$, 
\ba
    \mu_1(\psi_{12}^{-1} (A_2)) &=&\mu_1(\psi_{13}^{-1}\psi_{32}^{-1}(A_2))  
            < \mu_3((\psi_{32}^{-1}(A_2))^{d_{13}+\e/2}) +d_{13}+\e/2 \\
             & < & \mu_3(\psi_{32}^{-1}((A_2)^{d_{13}+d_{32}+\e/2})) +d_{13}+\e/2 \\
            &<& \mu_2((A_2)^{d+d_{32}+\e}) +d+\e.
\ea
Similarly, for  $\psi_{21}:=\psi_{31}\circ\psi_{23}$ we have
$\mu_2(\psi_{21}^{-1} (A_1)) < \mu_1((A_1)^{d+d_{31}+\e}) +d+\e$, 
for any Borel set $A_1\subset X_2$, 
and, therefore, the almost triangle inequality follows.
\end{proof}

\begin{lemma}
$d_{mGH}((X,\mu), (X', \mu'))=0$ if and only if there exists an isometry $\psi:X\to X'$ 
such that $\psi_*(\mu) = \mu'$.
\end{lemma}

\begin{proof}
Suppose $d_{mGH}((X,\mu), (X', \mu'))=0$. By definition, 
there are $\e_i$-approximations 
$\psi_i : X \to X'$ with $\lim \e_i =0$ such that 
\[
        \mu(\psi_i^{-1}(A')) <  \mu'((A')^{\e_i}) + \e_i,
\]
for every closed subset $A'\subset X'$. As $X, X'$ are compact,
we may assume, using (\ref{GH-approximation}),
that $\psi_i$ uniformly converges to an isometry
$\psi:X \to X'$. Since $\psi^{-1}(A') \subset (\psi_i^{-1}(A'))^{\d_i}$
for some $\d_i \to 0$, it follows that
\ba
  && \mu(\psi^{-1}(A'))  \le  \mu( (\psi_i^{-1}(A'))^{\d_i}) \\
          &&  \hspace{20mm} \le \mu(\psi_i^{-1}((A')^{\d_i+\e_i})) 
             \le \mu'((A')^{\d_i+2\e_i}) + \e_i.
\ea
Letting $i\to\infty$, we obtain $ \mu(\psi^{-1}(A')) \le  \mu'(A')$.
Taking complement, we have  $ \mu(\psi^{-1}(U')) \ge  \mu'(U')$
for any open set $U'\subset X'$.
It follows that $ \mu(\psi^{-1}((A')^{\e}))  \ge  \mu'((A')^{\e})$.
Letting $\e\to 0$, we obtain $ \mu(\psi^{-1}(A')) \ge  \mu'(A')$.
Thus, we have  $ \mu(\psi^{-1}(A')) =  \mu'(A')$
for every closed set $A'$ and hence for every Borel subset  $A'$.
This completes the proof of the lemma.
\end{proof}

\begin{proposition}
Let $\mu_i$ and $\mu$ be probability Borel measures on 
compact metric spaces $X_i$ and $X$.
A sequence $(X_i,\mu_i)$ weakly converges to $(X, \mu)$ 
if and only if  \\
$\lim\limits_{i\to \infty} d_{mGH}((X_i,\mu_i), (X, \mu))=0$.
\end{proposition}

\begin{proof}
Take $\e_i$-approximations $\psi_i : X_i \to X$ with $\lim\limits_{i\to\infty} \e_i = 0$
such that $\int_{X_i} f\circ \psi_i\, d\mu_i \to  \int_{X} f d\mu$
for every $f\in C(X)$.
First, using the weak convergence $(\psi_i)_{*} \mu_i \to \mu$,  we show
by contradiction that, for any Borel set $A\subset X$,
\beq 
            ((\psi_i)_{*}\mu_i)(A) & < & \mu(A^{\e_i'}) + \e_i', \label{eq:semi1}\\
        \mu(A) & <&  ((\psi_i)_{*}\mu_i)(A^{\e_i'})  + \e_i', \label{eq:semi2}
\eeq
for some $\e_i'\to 0$.
Suppose \eqref{eq:semi1} does not hold.
Then there are Borel sets $A_i$ of $X$ such that
\bequ 
    ((\psi_i)_{*}\mu_i)(A_i)  \ge  \mu(A_i^{c}) + c, \label{eq:semi3}
\eequ
for some constant $c>0$ independent of $i$.
We may assume that $A_i$ converges to a closed set $A$ with 
respect to the Hausdorff distance in $X$.
Take $\e_2>\e_1>0$ with 
$A_i\subset A^{\e_1}\subset A^{\e_2}\subset A_i^c$
for sufficiently large $i$. Choose $f\in C(X)$ such that
$0\le f\le 1$, $f=1$ on $A^{\e_1}$, and
 $\supp (f)\subset A^{\e_2}$.
Then 
\ba
   \mu(A_i^c) \ge \mu(A^{\e_2}) & \ge & \int_X f\,d\mu = \lim_{i\to \infty} \int_{X_i}f\circ\psi_i\,d\mu_i\\
     &  \ge &\limsup_{i\to \infty} ((\psi_i)_{*}\mu_i)(A^{\e_1}) \ge \limsup_{i\to \infty}  ((\psi_i)_{*}\mu_i)(A_i).
\ea
This is a contradiction to \eqref{eq:semi3}.

Next suppose \eqref{eq:semi2} does not hold. 
Then for some Borel sets $A_i$ of $X$ we have
\bequ
   \mu(A_i) \ge  ((\psi_i)_{*}\mu_i)(A_i^{c})  + c, \label{eq:semi4}
\eequ
for some constant $c>0$ independent of $i$.
Let $A_i\subset A^{\e_1}\subset A^{\e_2}\subset A_i^c$ and $f\in C(X)$ be 
given as above. Then, 
\ba 
   \liminf_{i\to \infty}  ((\psi_i)_{*}\mu_i)(A_i^c) & \ge & \liminf_{i\to \infty} ((\psi_i)_{*}\mu_i)(A^{\e_2})\\
   & \ge &\lim_{i\to \infty} \int_{X_i} f\circ\psi_i\,d\mu_i = \int_{X}f\,d\mu \ge \mu(A^{\e_1})\ge \mu(A_i).
\ea
This is a contradiction to \eqref{eq:semi4}.

Let $\psi_i':X\to X_i$ be any measurable $\e_i$-approximation such that
$d_i(\psi_i'\circ \psi_i(x_i),x_i)<\e_i$ for every $x_i\in X_i$, and $\,d(\psi_i\circ \psi'_i(x),x)<\e_i$ for every $x\in X$, see e.g. \cite[Lemma 2.5]{Fuk_inv}.
Using \eqref{eq:semi2}, we obtain, for any Borel $A_i\subset X_i$,
\ba
   \mu((\psi_i')^{-1}(A_i)) & < &\mu_i(\psi_i^{-1}(((\psi_i')^{-1}(A_i))^{\e_i'}))  + \e_i'\\
         & < &\mu_i(A_i^{2\e_i+\e_i'}) +\e_i'.
\ea
Together with \eqref{eq:semi1}, we have $d_{mGH}((X_i,\mu_i), (X,\mu)) < 2\e_i+\e_i'$.

Finally we shall prove the converse. Since $\mu_i, \mu$ are probability measures,
shifting $f\in C(X),\, f \mapsto f- \min(f)$, and normalising it,
$f \mapsto f/ \max(f),$
we may assume $0 \leq f \leq 1$.  Take a large positive integer $k$, and set
$A_j := \{ x\in X\,|\, f(x) \ge j/k\}$ for $0\le j \le k$.
It is straightfoward to see that, for any Borel measure $\mu$ on $X$,
\bequ
   \frac{1}{k}\sum_{j=1}^k \mu(A_j) \le \int_X f\,d\mu \le \frac{1}{k} + \frac{1}{k}\sum_{j=1}^k \mu(A_j).
                     \label{eq:straight}
\eequ
Now take $\e_i$-approximation $\psi_i :X_i \to X$ with $\e_i\to 0$ such that 
$((\psi_i)_{*}\mu_i)(A) < \mu(A^{\e_i}) +\e_i$ for every closed set $A\subset X$, where $\e_i\to 0$.
Letting $i\to\infty$,  we have $\limsup ((\psi_i)_{*})(A) \le \mu(A)$.
Therefore, in the above situation, we obtain
$\limsup\limits_{i\to \infty} ((\psi_i)_{*}\mu_i)(A_j) \le \mu(A_j)$ for each $0\le j\le k$. 
It follows from \eqref{eq:straight} that
$\limsup\limits_{i\to \infty} \int_X f\,d((\psi_i)_{*}\mu_i) \le  \frac{1}{k} + \int_X f\,d\mu$, and 
letting $k\to\infty$, 
\[
      \limsup_{i\to \infty} \int_X f\,d((\psi_i)_{*}\mu_i) \le \int_X f\,d\mu.
\]
Replacing $f$ by $1-f$, we get 
$\liminf\limits_{i\to \infty} \int_X f\,d((\psi_i)_{*}\mu_i) \ge \int_X f\,d\mu$,
and 
\[
   \lim_{i\to \infty}\int_X f\,d((\psi_i)_{*}\mu_i) =\int_X f\,d\mu,
\]
 as required.
\end{proof}


\section{Smoothness of the density functions}\label{smooth}

In this section, we consider the class ${\frak M}{\frak M}_p(n, \Lambda, D)$
of pointed closed Riemannian manifolds $(M, p, \mu_M)$ with $M\in {\frak M}(n, \Lambda, D)$ and the normalized 
measure $\mu_M = dV_{M}/\hbox{Vol}(M)$.
By Fukaya \cite{Fuk_inv}, we have the precompactness of 
${\frak  M}{\frak M}_p(n, \Lambda, D)$ with respect to the pointed 
measured Gromov-Hausdorff topology.

Let us consider a sequence $M_i\in {\frak M}(n, \Lambda, D)$
converging to $X \in \overline{{\frak M}}(n, \Lambda, D) $.
Let $\varphi_i:M_i\to X$ be a measurable $\e_i$-approximation with $\lim\limits_{i\to \infty} \e_i = 0$.
Passing to a subsequence, we may assume that 
$(M_i, \mu_{M_i})$ converges to some $(X,\mu)$ in the 
pointed measured  Gromov-Hausdorff  topology. 
Here $\mu$ is some probability measure of $X$, see \cite[Section 3]{Fuk_inv}. More precisely,
the pushfoward measure $(\varphi_i)_{*}(\mu_i)$ 
weak$^*$ sub-converges to $\mu$.

The following lemma is known in  \cite{Fuk_inv}.
Recall the notations \eqref{def-sin-reg} that $S(X)=X^{sing}$ denotes the singular set of $X$ and $X^{reg}$ denotes the regular part of $X$.
We say a point $x\in X$ is an orbifold point if the group $G_x$ stated in Theorem \ref{thm:limitstr}(1) is finite.

\begin{lemma}[\cite{Fuk_inv}] \label{lem:fuk-den}
Let  $\hat X^{sing}$ be the set of all points of $X^{sing}$
which are not orbifold points. Then
 \begin{enumerate}
 \item[$(1)$] $\mu(X^{sing})=0;$
 \item[$(2)$] there exists a continuous density function $\rho_X$ on $X$
   with $\mu = \rho_X \, \mu_X$, where $\mu_X$ denotes the normalized 
   Riemannian volume element of $X^{reg}$, such that 
   $\hat X^{sing}= \{ x\in X^{sing}\,|\, \rho_X(x)= 0 \}$.
 \end{enumerate}
\end{lemma}

In this section, we discuss some properties of $\rho_X$ concerning 
the smoothness, and prove the following:

\begin{lemma} \label{lem:chi_X}
   $\rho_X$ is of class $C^2_*$ on the regular part $X^{reg}$.
\end{lemma}

Concerning Lemma \ref{lem:chi_X},  Kasue \cite{Kas:measII,K93} proved   that
$\rho_X$ is of class $C^{1,\alpha}$ for any $0<\alpha < 1$. The method
in \cite{Kas:measII,K93} used smooth approximations of the metric of 
$M_i$ as in the proof of Lemma \ref{injectivity radius}.
Our method discussed below  is more direct  and contains an extension of 
Montgomery and Zippin's result on the smoothness of 
isometric group actions, see \cite{LLY},
although we follow a basic line in \cite{Fuk_inv}. 

First we consider the case when $X^{sing}$ is empty, namely the case when 
$X$ is a Riemannian manifold.
In this case, we can approximate $\varphi_i$ by an
almost Rimennian submersion $f_i:M_i \to X$ such that 
(see \cite[Section 3]{Fuk_inv})
 \begin{enumerate}
      \item[$(1)$] the pushfoward measure $(f_i)_*(\mu_i)$ weak* converges
      to $\mu$,
     \item[$(2)$] $\mbox{Vol}(f_i^{-1}(q))/\mbox{Vol} (M_i)$ converges 
      to $\rho_X(q)/\mbox{Vol} (X)$ in the $C^0$-topology.
  \end{enumerate}

Fix $q_0\in X$ and put  $q_i:= \psi_i(q_0)$, where $\psi_i:X\to M_i$ is
an $\e_i$-approximation such that 
$d_X(f_i\circ \psi_i(x), x) <\e_i$.
Let $B'\subset B\subset \Bbb R^n$, $(B, \tilde h_i,G_i)$ and  $(B, \tilde h_0,G)$ be
as in Section \ref{sec:review} so that
$(B',\tilde h_i)/G_i=B(q_i,\frac{\pi}{2\Lambda})$ and $(B',\tilde h_0)/G=B(q_0, \frac{\pi}{2\Lambda})$.
Let 
\[
      \pi_i:B' \to B(q_i, \frac{\pi}{2\Lambda})\subset M_i, \,\,
        \pi:B' \to B(q_0, \frac{\pi}{2\Lambda})\subset X,
\]
be the natural projections.

We now need the following result on the smoothness of isometric
group actions.

\begin{theorem}[\cite{LLY}] \label{thm:MZ}
Let $G$ be a Lie group, and $M$ be a Riemannian manifold of class $C^{k}_*$ with $k\ge 1$.
Suppose that the action of $G$ on $M$ is isometric. Then the $G$-action on $M$, 
\[
       G\times M \to M, \,\, (g, x)\to gx,
\]
is of class $C^{k+1}_*$, where we consider the analytic structure on $G$.
\end{theorem}

It is proved in  Calabi-Hartman \cite{CH} and  Shefel \cite{Shf} 
that the transformation $g:M \to M$ defined by 
each $g\in  G$ is of class $C^{k,\alpha}$.
Theorem \ref{thm:MZ} is a generalization of Montgomery-Zippin \cite[p. 212]{MZ},
where it is stated that the $G$-action on $M$ is of class $C^k$.
The proof of Theorem \ref{thm:MZ} can be found in \cite[Appendix A]{LLY}.

Note that in our situation, the pseudo-group $G$ can be extended to a
nilpotent Lie group (\cite{Fuk_JDG}).
It follows from Theorem \ref{thm:MZ} that 
the pseudo-group action of $G$ on $B$ is of class $C^3_*$.

Let $d=\dim(X)$, $m:=\dim(G)$ with $n=d+m$, and 
take a $d$-dimensional $C^{\infty}$-submanifold $Q$ of $B$ 
which transversally meets the orbit $G\cdot {\it O}$
at the origin ${\it O}$.
Let $s:U_0\to Q$ be a smooth coordinate chart  of $Q$ around $0$,
where $U_0$ is an open subset of $\R^d$.
From Theorem \ref{thm:MZ} together with the inverse function theorem, 
taking $Q$ smaller if necessary, we may assume that,
for some neighborhood $U$ of ${\it O}$ in $B$ and for a neighborhood 
$G^*$ of the identity in  $G$, the mapping
\bequ \label{000.9}
      G^*\times U_0 \to U,\,\, (g,x)\to g(s(x)), \quad g \in G^*,\,\,x \in U_0,
\eequ
gives  $C^3_*$-coordinates in $U$.

We can consider every element  $V$ of the Lie algebra $\frak g$  of $G$
as a Killing field on $B'$
by setting 
\[
     V(\x) := \frac{d}{dt}(\exp tV\cdot \x) \big|_{t=0},
\]
where $\x$ denotes points in $U$. Then, 
for any $\x\in U$, there is a unique $g\in G^*$ and $x\in U_0$ with
$\x=g(s(x))$. We define 
\[
      \bar \rho (\x) = \big|Ad_g(V_1)(\x)\wedge\cdots\wedge Ad_g(V_m)(\x) \big|,
\]
where $V_1,\ldots, V_m$ be a basis of $\frak g$, and the norm is taken with 
respect to $\tilde h_0$.
Here the adjoint representation $Ad:G \to GL(\frak g)$ is defined by
\[
          Ad_g(V) := \frac{d}{dt}(g\cdot \exp\, tV\cdot g^{-1}) \big|_{t=0}.
\]

Obviously  $\bar\rho$ is $G$-invariant. We show that $\bar\rho$ is of class $C^2_*$.
Recall that the correspondence $\x=g(s(x)) \to (g,x)$ is of class
$C^3_*$ from the inverse function theorem. 
It follows that  the map $(t, \x)\to g(\exp tV_i)g^{-1}(\x)$ is of
 class $C^3_*$. Thus $\x\to Ad_g(V_i)(\x)$ is of class $C^2_*$, and so is 
$\bar\rho(\x)$.

Since  $\bar\rho$ is $G$-invariant, there is a function 
$\rho$ defined on a neighborhood $\pi(U)$ of the point $q_0$ such 
that $\bar\rho = \rho\circ\pi$.
Then $\rho$ is of class $C^2_*$.

We shall prove that  $\rho_X$ is of class $C^2_*$ by 
showing that
$\rho_X(x) /\rho(x)$ is constant for $x \in \pi(U)$.
Basically we follow the argument in \cite{Fuk:II}.
Let 
\bfo
& &     G_i' :=\{ g\in G_i\,|\, d_{(B, \tilde h_i)}(g({\it O}), {\it O})<1/2\},\\ \nonumber
& &     G' :=\{ g\in G\,|\, d_{(B, \tilde h_0)}(g({\it O}), {\it O})<1/2\}.
\efo
Let us consider a left-invariant Riemannian metric on $G$.
First we need to show  that 
\[
       \frac{\mbox{Vol} (G'(s(x)))}{\rho(x)} = \textstyle{const},
\]
on a neighborhood of $q_0$. Define $F^x:G'\to G'(s(x))$  by
$F^x (g)=g(s(x))$. Let $V_1,\ldots, V_m$ be an orthonormal basis 
of $\frak g$. For any $g \in G'$,
we have
\[
    F^x_{*}(V_i(g)) =\frac{d}{dt} (g\exp tV_i (s(x))) \big|_{t=0}
                  = Ad_g(V_i)(F^x(g)).
\]
Therefore,
\ba
  \mbox{Vol} (G'(s(x))) &
       =& \int_{G'} |Ad_g(V_1)\wedge\cdots\wedge Ad_g(V_m)|(g(s(x)))\\
     & = &\int_{G'} \bar \rho (g(s(x))
      =  \rho (x) \mbox{Vol}(G'). 
\ea
For the rest of the argument, we can go through along the same line as
in \cite{Fuk:II}, which we outline below for reader's convenience.

Set
\ba
 E_i(x,\delta) & := & \{ \y \in B\,|\,  \hbox{there \,exists}\,\, g\in
   G_i'\, \,   \hbox{such\, that}\,\, d_{\tilde h_i}(\y, g(s(x))) < \delta  \}, \\
 E_0(x,\delta) & := & \{ \y \in B\,|\,  \hbox{there\, exists}\,\, g\in
   G' \,\,   \hbox{such\, that}\,\, d_{\tilde h_0}(\y, g(s(x))) < \delta  \}.
 \ea 
Then one can check 
\beq
  &&  \lim_{i\to\infty} \sup_{x\in \pi(U)} \frac{\mbox{Vol} (E_i(x,\delta))}{\mbox{Vol}
  (E_0(x,\delta))} = 1, \\
  && \lim_{\delta\to 0} \frac{\mbox{Vol}( E_0(x,\delta))}{\delta^d} = \omega_d\, \mbox{Vol}
   (G'(s(x))),
\eeq
where $\omega_d$ denotes the volume of unit ball in $\R^d$.
This implies 
\[
  \lim_{\delta\to 0} \frac{\mbox{Vol}( E_0(x,\delta))}{\mbox{Vol}(
  E_0(x^\prime,\delta))}\frac{\rho(x')}{\rho(x)} = 1,
\]
for all $x, x'\in W$. One can prove that there exists $c>0$ independent of
$x$
such that 
\bequ
 \lim_{\delta\to 0}\lim_{i\to\infty} \frac{\mbox{Vol}( E_i(x,\delta))}
       {\hbox{Vol}( G_i' )\,\delta^n \, \mbox{Vol} (f_i^{-1}(x))} = c.
\eequ
These yield \[
     \frac{\mbox{Vol}(f_i^{-1}(x))}{\mbox{Vol}(f_i^{-1}(x'))}
        \frac{\rho(x')}{\rho(x)} = 1.
\]
Since  
\[
       \lim_{i\to\infty} \frac{\mbox{Vol}(f_i^{-1}(x))}{\mbox{Vol} (M_i)} =
       \frac{\rho_X(x)}{\hbox{Vol}(X)},
\]
we conclude
\[
     \frac{\rho(x')}{\rho(x)}\frac{\rho_X(x)}{\rho_X(x')} =1,
\]
which shows that $\rho_X$ is of class $C^2_*$.

Next consider the general case when $X$ is not a Riemannian manifold.
Since the above argument is local, 
it follows  that $\rho_X$ is of class $C^2_*$ on $X^{reg}$.
This completes the proof of Lemma \ref{lem:chi_X}.
\hfill QED
\medskip

Later on, we will also need the properties of the 
class of the orthonormal frame bundles $FM$
over Riemannian manifolds 
$(M, p, \mu) \in \frak M \frak M_p(n, \Lambda, D)$. The frame bundles $FM$ are equipped with the Riemannian metric $\tilde h$,
inherited from $(M, h)$ and the corresponding  probability measure $\tilde \mu$, see (\ref{1.15.01.2012}).
We denote this class by  $\frak F \frak M \frak M (n, \Lambda, D)$. By Theorem \ref{thm:limitstr} and O'Neill's formula \cite{Oneill},
\beq \nonumber
& &\hbox{dim}(FM)=n+ \hbox{dim}(O(n)), \quad \hbox{diam}(FM) \leq D_F,
\\ 
& &|\textrm{sec}(FM,  \tilde h_{FM})| \leq \Lambda^2_F, \quad
\hbox{for}\, \, FM \in \frak F \frak M \frak M (n, \Lambda, D).
\label{2.15.01.2012}
\eeq
The closure $\overline {\frak F \frak M \frak M }(n, \Lambda, D)$
with respect to the measured Gromov-Hausdorff topology gives rise to $C^3_*$-smooth manifolds $Y$ equipped with $C^2_*$-smooth Riemannian metric $\tilde h_Y$,
which appears in Theorem \ref{thm:limitstr}(2) and in the proof
of Lemma \ref{lem:fuk-den}. 
The $C^3_*$-smooth structure of $Y$ is defined by the limit of the transition maps for harmonic coordinates on a sequence $FM_i$ converging to $Y$, see \cite{A90,AKKLT,GW:lipschitz,Pets:conv}.
Hence $Y$ has a compatible $C^{\infty}$-smooth structure contained in the $C^3_*$-smooth structure.
Analyzing the
proof of Lemma \ref{lem:fuk-den}, we have

\begin{corollary} \label{density_on_Y}
Let $Y \in \overline {\frak F \frak M \frak M }(n, \Lambda, D)$. Then $Y$ is a smooth manifold with $C^2_*$-smooth Riemannian metric $\tilde h_Y$ and strictly
positive density function
$\rho_Y\in C^2_*(Y)$. Moreover, $O(n)$ acts by isometries on $Y$ and $\rho_Y$ 
 is  $O(n)$-invariant.
\end{corollary}

Let  $M_i$ be a sequence in ${\frak M}(n, \Lambda, D)$
converging to  $X \in \overline{{\frak M}}(n, \Lambda, D) $,
and let $FM_i$ denote the orthonormal frame bundle of $M_i$ endowed with
the natural Riemannian metric, which  has uniformly bounded  sectional 
curvature and diameter.
Put
\bequ \label{1.15.01.2012}
       d \tilde\mu_i = \frac{dV_{FM_i}}{\hbox{Vol}(FM_i)},
\eequ
where $dV_{FM_i}$ is the natural volume element on $FM_i$ with
 $O(n)$ acting isometrically   on $FM_i$.
Passing to a subsequence, we may assume that 
$(FM_i,\tilde\mu_i, O(n))$ converges to $(Y,\tilde\mu, O(n))$ in the
equivariant measured Gromov-Hausdorff topology, where $Y$ is a Riemannian manifold and 
$\tilde\mu$ is a probability measure on $Y$ invariant under the 
$O(n)$-action. 
The notion of equivariant measured Gromov-Hausdorff topology is defined 
in a way similar to that of measured Gromov-Hausdorff topology, where
a measurable $\e_i$-approximation maps $\tilde \psi_i:FM_i \to Y$ with
the property $(\tilde \psi_i)_*(\tilde\mu_i) \to \tilde\mu$ is required to satisfy 
conditions similar to  \eqref{29.06.3}.
By Theorem \ref{thm:fib}, there are $\e_i$-regular maps
$\tilde f_i: FM_i \to Y$  and $f_i: M_i \to X$ 
such that $f_i\circ\pi_i = \pi\circ\tilde f_i$.
We may replace $\tilde \psi_i$ by $\tilde f_i$.

The probability measure $\tilde\mu$ on $Y$ can be written as 
\[
       d\tilde \mu = \rho_Y\, \frac{dV_Y}{\hbox{Vol}(Y)},
\]
with a strictly positive $O(n)$-invariant density function $\rho_Y$ by Corollary \ref{density_on_Y}.
It follows that there is a strictly positive function $\overline\rho_X$ on $X$ 
with 
    $\overline\rho_X\circ \pi = \rho_Y$, 
where   $\pi:Y\to X$ is the projection.
Since $d_X(x, x')=d_Y(\pi^{-1}(x),\, \pi^{-1}(x'))$, then $\overline\rho_X$ is Lipschitz.

The projection $\pi_i:FM_i \to M_i$ is a Riemannian submersion 
with totally geodesic fibers isometric to $O(n)$.
Thus,  
\[
   d  \tilde\mu_i = \frac{dV_{FM_i}}{\mbox{Vol}(O(n))\times
     \mbox{Vol}(M_i)}.
\]
Then it follows that $(\pi_i)_*(\tilde \mu_i)=\mu_i$.
For any continuous function $f$ on $X$, by Fubini's theorem, we have
\beq
& &  \int_X f(x)\rho_X(x)\, \frac{dV_X}{\hbox{Vol} (X)}=\int_X f d\mu=
 \lim_{i\to\infty}\int_{M_i}(f\circ  f_i)\,d\mu_{M_i}\\ \nonumber
       & &= \lim_{i\to\infty}\int_{FM_i}(f\circ f_i\circ\pi_i)\,d\tilde\mu_{i}
       = \lim_{i\to\infty}\int_{FM_i}(f\circ \pi\circ \tilde
       f_i)\,d\tilde\mu_{i}\\ \nonumber
& &       = \int_{Y}(f\circ\pi) \,d\tilde\mu = \int_{Y}(f\circ\pi) \,\rho_Y \frac{dV_Y}{\hbox{Vol}(Y)}
\\ \nonumber
       & &= \frac{1}{\mbox{Vol}(Y)} \int_{X}  \left(\int_{\pi^{-1}(x)}  \rho_Y(y) 
         \,d{\mathcal H}^{\ell}_{\pi^{-1}(x)}(y) \right) f(x)  dV_X(x)
	 \\ \nonumber
      & &=  \frac{1}{\mbox{Vol}(Y)} \int_{X}  f(x) \overline\rho_X(x) \mbox{Vol}_{\ell}(\pi^{-1}(x))
      \,dV_X(x), \label{AAAA}
\eeq
where ${\mathcal H}^{\ell}_{\pi^{-1}(x)}$ is 
the $\ell$-dimensional Hausdorff measure of ${\pi^{-1}(x)}$ with
  $\ell=\dim (O(n))$.
Then it follows
\bequ \label{000.15}
\pi_*(\tilde \mu)=\mu
\eequ
 and the function
\bequ \label{eq:rhoo}
 \rho_X(x) =  \overline \rho_X(x) \hbox{Vol}_{\ell}(\pi^{-1}(x))
      \frac{\hbox{Vol}(X)}{\hbox{Vol}(Y)} 
\eequ
is the required density function on $X$ with $d\mu_X = \rho_X  dV_X/ \hbox{Vol}(X)$, as stated in Lemma \ref{lem:fuk-den}(2).

\section{On properties of eigenfunctions}\label{section: Aux results}
In this section, we consider some properties of the eigenpairs $\{\la_j,\phi_j\}_{ j =0}^\infty$ 
of the weighted Laplace operator $\Delta_X $, $(X, p, \mu_X) \in
{\overline {\frak M \frak M}}_p (n, \Lambda, D)$ and $\{\tilde \la_j ,\tilde\phi_j\}_{ j =0}^\infty$ 
for $\Delta_Y $, $(Y, \mu_Y)  \in 
{\overline {\frak F \frak M \frak M}}_p (n, \Lambda, D)$.
The basic properties of these operators are given in Appendix \ref{OTA}.

As the Riemannian metric $h_Y$ and the density function $\rho_Y$ on $Y$ are of class $C^{2}_*$ by Corollary \ref{density_on_Y}, it is natural to work in the $C^3_*$-smooth structure of $Y$  
and
invariantly define
  the Sobolev spaces
$W^{k, q}(Y),\, k\in \{0, 1, 2\},$ $ 1 \leq q < \infty$, and H\"older spaces 
$C^{l, \a}(Y),\, l \in \{0, 1, 2\},\, 0\le \a <1,$ and the Zygmund space $C^{3}_*(Y)$.  

Moreover, due to Theorem \ref{thm:MZ} with $M=Y,\, G=O(n)$, we can define $O(n)$-invariant subspaces $W_O^{k, q}(Y) \subset W^{k, q}(Y)$,$\,C_O^{l, \a}(Y) \subset C^{l, \a}(Y)$
and $C^{3}_{*, O}(Y) \subset C^{3}_*(Y)$, which consist of functions invariant
with respect to $O(n)$-action, that is, functions $f\in W^{k,q}(Y)$ (resp. $C^{l,\alpha}(Y)$, $C^3_*(Y)$) satisfying $f(y)=f(o(y))$ for all $y\in Y$ and $o\in O(n)$.

Introduce the operator
\bequ \label{25.12.2011.2}
{\Bbb P}_O: L^2(Y) \to L^2_O(Y), \quad
\left({\Bbb P}_O u^*\right)(y)= \int_{O(n) y} u^*({\it o} ( y)) \, d_n {\mathcal H}^\ell_{O(n) y},
\eequ
where $d_n{\mathcal H}^\ell_{O(n) y}:= \frac{1}{{\rm Vol}_\ell(O(n) y)}d {\mathcal H}^\ell_{O(n) y}$ and $\ell=\textrm{dim}(O(n)\cdot y)=\dim (O(n))$,
where $d {\mathcal H}^\ell_{O(n) y}$ is the $\ell$-dimensional Hausdorff measure on the orbit $O(n)\cdot y$ of $y$.
In the future, the subindex $_O$ indicates
 $O(n)$-invariance of functions.

\begin{lemma} \label{PO-action} The operator ${\Bbb P}_O$ in (\ref{25.12.2011.2}) 
is an orthogonal projector in $L^2(Y, \mu_Y) $. Moreover,
for  $ l \in \{0, 1, 2\},\, 1 \leq q < \infty$, 
\beq \label{000.16}
& &{\Bbb P}_O: L^q(Y) \to L^q_O(Y),\quad {\Bbb P}_O: W^{2, q}(Y) \to W_O^{2, q}(Y),\,\, \\ \nonumber
& &{\Bbb P}_O: C^{l, 1}(Y) \to C_O^{l, 1}(Y), \quad
{\Bbb P}_O:  C^{3}_{*}(Y) \to C^{3}_{*, O}(Y),
\eeq
are bounded.
In addition, for $ 1\le q <\infty,\,  k\in \{0, 1, 2\}$, $W_O^{k, q}(Y)$, $C^{3}_{*, O}(Y)$ are dense in $L_O^q(Y, \mu_Y)$.
\end{lemma}

\begin{proof}
Let $u^*,\, v^* \in L^2(Y)$. Then,
\bfo
& &\left( {\Bbb P}_O u^*,\, v^* \right)_{L^2(Y)} = 
\int_Y \left( \int_{O(n)y} u^*({\it o}(y)) d_n {\mathcal H}^\ell_{O(n) y} \right) v^*(y) d\mu_Y(y)
\\
\nonumber 
& &
 =\int_Y  \left(\int_{O(n) \tilde y}  v^*({\it o}^{-1}(\tilde y))  d_n {\mathcal H}^\ell_{O(n) \tilde y}
 \right)u^*(\tilde y) d\mu_Y(\tilde y) 
=
\left(  u^*,\, {\Bbb P}_O v^* \right)_{L^2(Y)}.
\efo
Here we have made the substitution $\tilde y= {\it o} y$, and used 
the invariance of $d {\mathcal H}^\ell_{O(n) y}  $ with respect to $O(n)$ and the fact that $\mu_Y$ is $O(n)$-invariant.

Next we prove (\ref{000.16}) for $L^q$. We have
\bfo
& &\|{\Bbb P}_O u^*  \|_{L^{q}(Y)}^q= 
\int_Y \left|\left( \int_{O(n) y} u^*({\it o} ( y)) \, d_n {\mathcal H}^\ell_{O(n) y}\right)\right|^q d\mu_Y(y) \\
& &\leq \int_Y  \int_{O(n) y} |u^*({\it o} ( y))|^q\, d_n {\mathcal H}^\ell_{O(n) y} d\mu_Y(y) =   \|u^*\|_{L^q(Y)}^q,
\efo
where we used H\"older's inequality, invariance of $d_n {\mathcal H}^\ell_{O(n) y} $ and $\mu_Y$ with respect to $O(n)$, and $ \int_{O(n) y} d_n {\mathcal H}^\ell_{O(n) y} 
=1$.

As, by \cite[Thm.\ 9.4.1]{Kry}, 
$
(I+\Delta_Y)^{-1}: L^q(Y) \to 
W^{2, q}(Y),
$
is an isomorphism and $\Bbb P_O$ commutes with $\Delta_Y$, we obtain (\ref{000.16}) for
$W^{2, q}(Y)$. 

Since $O(n)$ acts by isometries, (\ref{000.16}) on $C^{0,1}$ and $C^1_*$ follows from definition (\ref{25.12.2011.2}).
As by \cite[Prop. 2.4.1]{AKKLT} $(I+\Delta_Y)^{-1}: C^1_*(Y) \to 
C^3_*(Y)$
is an isomorphism, we obtain (\ref{000.16}) on $C^3_*(Y)$.

To obtain the density of $C^{3}_{*,O}(Y)$ and therefore, $W_O^{k, q}(Y)$ in $L_O^q(Y)$, 
we approximate $u^* \in L^q_O(Y)$
by $C^3_*$-functions $\tilde u^*_k$ and consider
$
u^*_k ={\Bbb P}_O \tilde u^*_k \in C^3_{*, O}(Y).
$
Since $u^*,u^*_k$ are $O(n)$-invariant,
then $u^*(y)-u^*_k(y)=\int_{O(n)y} \left( u^*-\tilde u^*_k\right) ({\it o} ( y)) \, d_n {\mathcal H}^\ell_{O(n) y}$ and by H\"older's inequality,
\begin{equation*}
\label{25.12.2011.3}
\|u^*-u^*_k\|_{L^q(Y)} \leq \|u^*-\tilde u^*_k\|_{L^q(Y)}.
\end{equation*}
Hence $u^*_k$ approximates $u^*$ in $L^q$.
\end{proof}

Let $(X,\mu_X) \in {\overline {\frak M \frak M}}_p (n, \Lambda, D)$ and $(Y,\mu_Y) \in {\overline {\frak F\frak M \frak M}}_p (n, \Lambda, D)$ such that $X=Y/O(n)$, as constructed in Theorem \ref{thm:limitstr}. Let $\pi:Y\to X$ be the projection.
As $X^{reg}$ is a $C^3_*$-smooth manifold (see \cite{AKKLT}), we  define function spaces
 $W^{k, q}_{loc}(X^{reg})$, $C^{l, \a}(X^{reg})$ and $C^3_*(X^{reg})$.

\begin{lemma} \label{smoothness_of_eigenfunctions_g}
(1). Let $(Y,\mu_Y) \in {\overline {\frak F\frak M \frak M}}_p (n, \Lambda, D)$. Then the eigenfunction 
$
\tilde\phi_j \in C^3_*(Y).
$
\smallskip

(2). Let $(X,\mu_X) \in {\overline {\frak M \frak M}}_p (n, \Lambda, D)$. Then the eigenfunction
\beq \label{1.05.11.11}
\phi_j \in   C^{0, 1}(X) \bigcap
    C^{3}_*(X^{reg}) 
\eeq
\end{lemma}

\begin{proof}
(1)\, Since $Y$ is a $C^3_*$-manifold with $C^2_*$ metric tensor, the $C^3_*$-smoothness of $\tilde\phi_j $ follows from the standard 
interpolation arguments by the interior Schauder estimates, see e.g.
\cite[Thm. 6.17]{GTr}. Since $h_Y$ and $ \rho_Y$ are $C^{2}_*$-smooth,
these arguments can be adjusted to $C^{3}_*$-case by means of 
\cite[Prop.\ 2.4.1]{AKKLT}, see also 
\cite[Thm.\ 14.4.2-3] {Tay} or \cite{Tri}.

\smallskip
 \noindent (2) Denote by $\Delta_Y^O$ the $O(n)$-invariant part of $\Delta_Y$ on $L^2_O(Y)$, namely, $\Delta_Y^O:=\Delta_Y \, \mathbb{P}_O$ where $\mathbb{P}_O:L^2(Y) \to L^2_O(Y)$ is the projection onto the $O(n)$-invariant subspace of $L^2(Y)$.
 By \cite[Lemma 7.1]{Fuk_inv}, see also Appendix \ref{OTA},
 \bequ \label{26.11.2011.1}
 \hbox{spec}(\Delta_Y^O)=\hbox{spec}(\Delta_X), \quad  \phi^O_j= \pi^* \phi_j,
 \eequ
 where $ \la^O_j, \phi^O_j$ are the eigenvalues and eigenfunctions of 
 $\Delta_Y^O$. Then,
 \beq \label{26.11.2011.2}
 |\phi_j(x)-\phi_j(y)|=|\tilde\phi_j^O(x^*)-\tilde\phi_j^O(y^*)| \leq 
 \|\tilde\phi_j^O\|_{C^{0, 1}(Y)}\, d_X(x, y),
 \eeq
 where
  $x, y \in X$ and $x^* \in \pi^{-1}(x),\, y^* \in \pi^{-1}(y),$
 satisfy
 $ d_X(x, y)=  d_Y(x^*, y^*).
$
 The $C^3_*$ inclusion in (\ref{1.05.11.11}) follows
 from (\ref{26.11.2011.1}) and claim (1).
\end{proof}

Denote by $Z_X, \, Z_Y, \, Z_Y^O$ the linear subspaces of finite linear combinations of
the eigenfunctions $\{\phi_j\}_{j=0}^\infty, \,\{\tilde\phi_j\}_{j=0}^\infty$,
and $\{\phi^O_j\}_{j=0}^\infty$, respectively.

\begin{proposition}\label{lem: density}
\begin{enumerate}
\item 
$Z_Y,\, Z_Y^O$ are dense in $W^{k, q}(Y), \,W^{k, q}_O(Y)$,  for any 
$k \in \{0, 1, 2\}$, $1 \leq  q <\infty$, $C^{l,\a}(Y),\, C_O^{l,\a}(Y),\, l\in\{0,1, 2\},\, 0\le \a <1,\,$
and $C^3_*(Y),\, C^3_{*, O}(Y)$, respectively.
\item
$Z_X$ is dense in $W^{k, q}_{loc}(X^{reg})$, $C^{l,\a}(X^{reg})$, $C^3_*(X^{reg})$
and  $C^{0, 1}(X)$.
\end{enumerate}
\end{proposition}

\begin{proof}
{\it 1.} We start with $W^{2, q}$. Using the density of $W^{k, q}$ in $W^{k, q'}$, if $q>q'$, 
we consider only the case
$2 \leq q <\infty$.
As, by \cite[Thm.\ 9.4.1]{Kry}, 
\[
(I+\Delta_Y)^{-1}: L^q(Y) \to 
W^{2, q}(Y),
\]
is an isomorphism, it is enough to show that $Z_Y$ is dense in $L^q$. Clearly, 
this is true for $q=2$ and we will use boot-strap arguments to show this for any $q>2$.

Indeed, assume that $Z_Y$ is dense  in $L^q$ for some $q\geq 2$. Then, as 
$(I+\Delta_Y)^{-1} (Z_Y)=Z_Y$, $Z_Y$ is dense  in $W^{2, q}$. Using Sobolev's
embedding and the fact that $C^3_{*}(Y)$ is dense in $L^{q'}(Y),\, q'<\infty$,
 $Z_Y$ is dense  in $W^{2, q'}$, if 
$
\frac{1}{q'}> \frac{1}{q}-\frac{1}{\hbox{dim}(Y)}.
$
Iterating the arguments, we obtain
the density of $Z_Y$ in $W^{k, q}(Y), \,k=0, 1, 2,\, q \in [2, \infty).$


Observe now that, due to $W^{k, q} \subset C^{0,1}$ for $q> \hbox{dim}(Y)/2$, $Z_Y$ is dense in $ C^{0,1}$. Then the
interior Schauder regularity method, adjusted to Zygmund classes, see e.g. \cite[Prop.\ 2.4.1]{AKKLT} ,
\cite[Thm.\ 14.4.2-3] {Tay} or \cite{Tri}, shows that $(I+\Delta_Y)^{-1}: C^{1}_*(Y)\to C^{3}_*(Y)$ is an isomorphism,
completing the proof of (1) for $Z_Y$.

Turning to the case of $Z_Y^O$, we just note that
\bfo
\left(I+\Delta_Y \right)^{-1}: L^{ q}_O(Y)\to  W^{2,q}_O(Y), \quad
\left(I+\Delta_Y \right)^{-1}: C^1_{*,O}(Y) \to C^3_{*, O}(Y),
\efo
are isomorphisms.
Since $Z_Y^O$ is dense in $L^2_O(Y)$, repeating the 
above arguments gives
 the desired result.

\smallskip
{\it 2.} Note that 
\bequ \label{26.11.2011.3}
C^{0,1}_O(Y)= \pi^*(C^{0,1}(X)),
\eequ
see Appendix \ref{OTA}
or \cite[Section 7]{Fuk_inv}, where the classes $C^1$ can be easily 
substituted by classes
$C^{0, 1}$.  This, together with 
(\ref{26.11.2011.1}), (\ref{26.11.2011.2}) and the density of $Z^O_Y$ in 
$C^{0,1}(Y)$ provides the density 
of $Z_X$ in $C^{0, 1}(X)$.

Next, to prove the density of $Z_X|_{X^{reg}}$  in $W^{k, q}_{loc}(X^{reg}), \, k \in\{0, 1, 2\},\, 1 \le q<\infty$,
or in $C^{l, \a}(X^{reg}),\, l\in\{0,1,2\},\, 0\le \a <1$, it is sufficient to show that 
$Z_X|_{X^{reg}}$ is dense in $C^3_*(X^{reg})$.  To this end,
let $ K$ be a compact set in $X^{reg}$.
To prove that last statement, it is sufficient to show
that any $u \in C^3_*(X^{reg}), \, \hbox{supp}(u) \subset K,$ can be approximated in $C^3_*(K)$ by a sequence of
functions $z_m|_K,\, z_m \in Z_X$.  
To this end, consider
\[
u^* = \pi^*(u) \in C^3_{*, O}(Y).
\]
Since $O(n)$ acts isometrically on $\pi^{-1}(K)$, it follows from Theorem \ref{thm:MZ} that 
\bequ \label{000.12}
\pi_*: C^3_{*,O}(\pi^{-1}(K)) \rightarrow C^3_*(K)
\eequ
is a bounded operator.
 
By part (1), there are  $\{z^O_m\}_{ m=1}^\infty,\, z^O_m \in Z^O_Y$ which approximate
$u^*$ in $C^3_*(Y)$. Thus, $\{z^O_m|_{\pi^{-1}(K)}\}_{ m=1}^\infty\subset Z_X$ approximates $u^*|_{\pi^{-1}(K)}$
in $C^3_{*, O}(\pi^{-1}(K))$ and, by (\ref{000.12}),
$$
\lim_{m\to\infty}\pi_*(z^O_m)|_K = \pi_*(u^*)|_K= u|_K \quad \hbox{in}\,\, C^3_*(K).
\vspace{-5mm}
$$
\end{proof}

\begin{lemma}\label{lem: eigenfunctions injectivity}

(i) The map $\Phi:\,X\to \R^{\N}$ defined by
\bequ
\label{def-Phi}
\Phi(x):=\big(\phi_j(x) \big)_{j=0}^\infty
\eequ
is injective, that is, if $x \neq y$, there is an index 
$j\in \N$ such that $\phi_j(x) \neq \phi_j(y)$.

\smallskip
\noindent (ii)  
For any $x_0 \in X^{reg}$, there is a neighborhood $ U \subset X^{reg}$ of $x_0$
and
indices ${\bf j}=(j_1,\dots, j_d),$ where $ d =\dim(X),$ such that 
\ba
\Phi_{\bf j}:U\to \R^d,\quad \Phi_{\bf j}(x)=\big(\phi_{j_k}(x) \big)_{k=1}^d,\quad x\in U,
\ea
are $C^3_*$-smooth  coordinates  in $U$.
\end{lemma}

\begin{proof}
Suppose that there exists $x\neq y$ such that $\phi_j(x)=\phi_j(y)$ for all $j\in \N$.
Then $f(x)=f(y)$ for all $f\in Z_X$.
Due to the density of $Z_X$ in $C^{0,1}(X)$ (Proposition \ref{lem: density}), this would imply that there is no Lipschitz function taking different values at $x$ and $y$, which is clearly a contradiction e.g., $f(z)=d_X(z,y)/d_X(x,y)$.
The second claim is due to \cite[Lemma 4.2.1]{AKKLT}.
\end{proof}

By Lemma \ref{lem: eigenfunctions injectivity}, we can use the following test to identify
${\bf i}:=(i_1,\dots,i_d)\in \N^d$ such that $\Phi_{\bf i}=(\phi_{i_1}, \dots, \phi_{i_d})$ form a coordinate
system.

\begin{corollary}\label{lem: recognize coordinates}
 Let   $x_0\in X^{reg}$ and,
for ${\bf i}:=(i_1,\dots,i_d)\in \N^d$, 
a neighborhood $ W_{{\bf i}}\subset \R^d$ of $\Phi_{\bf i}(x_0),$
  such that the
function $\phi_l\circ \Phi_{\bf i}^{-1}:W_{\bf i} \to \R$ is $C^{3}_*$-smooth, for any $l \in \{i_1,\cdots,i_d\}$.
Then there is a neighborhood $U$ of $x_0$ such that
$\Phi_{\bf i}:U\to \Phi_{\bf i}(U)$ is a $C^{3}_*$-smooth diffeomorphism.
\end{corollary}

The above results show that the eigenfunctions $\{\phi_j(x)  \}_{j=0}^\infty$ given on 
open set
$\Omega^{reg}\subset X^{reg}$ 
determine the topological and $C^3_*$-differentiable structure of 
$\Omega^{reg}$.
Next we consider the metric tensor $h$ and the density function $\rho$ on $\Omega^{reg}$.

\begin{lemma} \label{recovery}
Let $\Omega^{reg}\subset X^{reg}$ be an open set.
The set $\Omega^{reg}$ and the eigenpairs $\{\la_j,\, \phi_j|_{\Omega^{reg}} \}_{j=0}^\infty$ of $\Delta_{X}$ uniquely determine the metric $h|_{\Omega^{reg}}$ and a function $\normalized \rho$ on $\Omega^{reg}$,
such that $\normalized \rho(x)= \normalized c \rho(x)$,
where
 $\normalized c >0$ is a constant.
\end{lemma}

\begin{proof}
Let $x_0 \in \Omega^{reg}$ and ${\bf x}(x)=(x^1, \dots, x^d),\,{\bf x}(x_0)=0$
be $C^{3}_*$-smooth coordinates 
in a neighborhood $U \subset \Omega^{reg}$ of $x_0$.
Let $\chi(x)$ be a $C^{3}_*$-smooth function with $\hbox{supp}(\chi) \subset \subset U$
and $\chi(x)=1$ in a neighborhood $V$ of 0.
Let 
$$
\chi_i(x)= x^i \chi(x), \; 1\leq i\leq d, \quad
\chi_{j, k}(x)= x^j x^k \chi(x),\;1 \leq j \leq k \leq d.
$$
By Proposition \ref{lem: density}(2), there are  $z^\ell_i,\,z^\ell_{j,k} \in Z_X$ such that, in $C^3_*(U)$,
$$
\chi_i=\lim_{\ell \to \infty} z^\ell_i,\quad 
\chi_{j,k}=\lim_{\ell \to \infty} z^\ell_{j, k}.
$$
Therefore, for sufficiently large $N$ and $y\in V$,
the vectors $\Psi_y(\phi_m)$, $m=0,1,\dots, N,$  span the space $\R^{d(d+3)/2}$, where  
$$
 \Psi_y[f]= \left(\left(\frac {\p f}{\p x^i}(y)\right)_{1\leq i\leq d}, \left(\frac {\p^2 f}{\p x^j\p x^k}(y)\right)_{1\leq j \leq k\leq d}\right).
$$
Consider the equations $\Delta_X \phi_m=\la_m \phi_m$ at an arbitrary point $y\in V$ in coordinates ${\bf x}$, 
\bequ \label{determine-metric}
\Delta_X \phi_m (y)=-h^{jk}(y) \frac{\partial^2 \phi_m}{\partial x^j \partial x^k}(y)- a^j(y) \frac{\partial \phi_m}{\partial x^j}(y)=\lambda_m\phi_m(y), \quad m=0,1,\dots,N,
\eequ
where
$$
a^j=\frac{1}{\sqrt h \, \rho}\,\frac{\partial}{\partial x^k} \left(\sqrt h\, h^{jk} \rho \right)= 
\frac{1}{\sqrt h }\,\frac{\partial}{\partial x^k} \left(\sqrt h\, h^{jk} \right)+ h^{jk} \frac{\partial}{\partial x^k} \log (\rho).
$$
Regarding \eqref{determine-metric} as a system of li equations, the vectors $\Psi_y(\phi_m)$, $m=0,1,\dots, N,$ uniquely determine the coefficients $h^{jk}(y)$ and $a^j(y)$, which determines $h^{jk}(y)$ and, up to a multiplicative constant, $\rho(y)$.
\end{proof}

From now on, we choose and fix some positive function $\normalized \rho$ so that 
\begin{equation} \label{rho-normalized}
\int_{\Omega^{reg}} \normalized \rho(x) dV_X=1,\quad \normalized \rho= \normalized c \rho,
\;\normalized c>0.
\end{equation}

Now we move to the consideration of the pointwise heat data, $PHD$, see (\ref{heat kernel1})-(\ref{intro-PHD}).
We denote by ${\mathcal L}(\la_j)$ the eigenspace of $\Delta_X$,
corresponding to the (multiple) eigenvalue $\la_j$.

\begin{lemma} \label{PHD_to_LSD}
Let $PHD=(H(z_\alpha, z_\beta; t_{\ell}))_{ \alpha, \beta\in \N,\, \ell \in \Z_+}$ be the pointwise heat data
of $\Delta_{X}$ with  $\{z_\alpha\}_{\alpha\in\N}$ dense in an open set $\Omega_0 \subset X$.
Then the data uniquely determine the closure $\overline{\Omega_0}$, the eigenvalues $\la_j$ and, up to an orthogonal transformation in  
${\mathcal L}(\la_j)$, the orthonormalized eigenfunctions $\phi_j|_{\overline{\Omega_0}}$.
\end{lemma}

\begin{proof}
Due to the analyticity of $H(x, y; t)$ with respect to $t,\, \hbox{Re}(t)>0$,
$H_{\alpha, \beta; \ell},\, \ell \in \Z_+$, determine  $H_{\alpha, \beta}(t) \in C(0, \infty)$,
where
$H_{\alpha, \beta}(t)= H(z_\alpha, z_\beta; t).$
As follows from Weyl's asymptotics
for eigenfunctions (Lemma \ref{BBG}) together with the Lipschitz regularity estimate (Lemma \ref{c-beta}), for a given
$X$,
\bequ 
\la_j \geq c\, j^{2/n},\quad \|\phi_j\|_{C^{0,1}(X)} \leq C (1+\lambda_j^2)^s,
\eequ
for some uniform constants $c,s,C$.
This implies that, in the standard
representation,
\bequ \label{000.10}
H(x, y, t)= \sum_{j=0}^\infty e^{-\la_j t} \phi_j(x) \phi_j(y),
\eequ
where the series in the right-hand side converges in $C^{0,1}( X \times X \times (0, \infty))$.
Then the Laplace trasform $\hat H_{\alpha, \beta}(\omega)$ of $H_{\alpha, \beta}(t)$
has a meromorphic continuation to the
whole plane $\C$. 
It has simple poles at $\omega=-\la_j,\, j \in \N$, with corresponding residues,
\bequ \label{04.18.04}
\underset{\omega=-\la_j} {\hbox{Res}}(\hat H_{\alpha, \beta}(\omega))=
\sum_{\la_{j'}=\la_j}  \phi_{j}(z_\alpha) \phi_{j}(z_\beta).
\eequ
Using the results in \cite[Lemma 4.9]{KKL}, these residues
determine, up to an orthogonal li transformation
 in the orthogonal group $O(m_j)$, where $m_j$ is the multiplicity of the eigenvalue $\lambda_j$, 
 the values $\phi_j(z_\alpha)$.  
 
 The map $\Phi: \overline{\Omega_0} \to \R^{\N}$ defined in \eqref{def-Phi} is injective and continuous with respect to 
 the product
 topology in $\R^{\N}$. As $\overline{\Omega_0}$ is compact, $\Phi$ is a homeomorphism onto its image
so that $\phi_{j}(z_{\alpha}),\, j, \alpha \in \N$,
 determine $\Phi(\overline{\Omega_0})$.
\end{proof}

 It follows from the the proof
 of Lemma \ref{PHD_to_LSD} that
 
 \begin{corollary} \label{heat_identification}
Let $X \in \overline{ \frak M \frak M}_p (n,\Lambda,D)$. Assume that $H_X(x, z, t)=H_X(x', z, t)$
 for all $z \in U,\, t \in (a, b) \subset \R_+$, where $U \subset X$ is open. Then $x=x'$.
 \end{corollary}

\begin{remark} \label{regular2}{\rm The pointwise heat data on any open set $\Omega \subset X$ determine
the dimension $d= \hbox{dim}(X)$.
Indeed, $d$ is
the minimal number so that, 
in some open set $\Omega' \subset \Omega$, there are
$d$ eigenfunctions $\phi_{j(1)},\dots, \phi_{j(d)}$, which 
form $C^3_*$-coordinates in $\Omega'.$}
\end{remark}

\section{Continuity of the direct map} \label{Continuity of the direct map}

The principal goal of this section is to show that, in $\overline{{ \frak M \frak M}}_p (n,\Lambda,D)$,
the pointed measured Gromov-Hausdorff convergence implies the convergence
 of the heat kernel associated to the weighted Laplacians.
We define the topology in the set of the heat kernels,
$H_X(\cdot, \cdot, t),\, X \in {\overline{ \frak M \frak M}}_p (n,\Lambda,D)$ in the following way.

\begin{definition} \label{heat-topology}
Let $H_X(x, y, t),\, x, y \in X, t >0,$ and $H_{X'}(x', y', t),\, x', y' \in X', t >0,$
be the heat kernels for spaces 
$(X, p, \mu_X),\, (X', p', \mu_{X'}) \in {\overline{ \frak M \frak M}}_p (n,\Lambda,D)$. 
We say $H_X, \, H_{X'}$ are
$\e$-close, if for every $r>0$, there exist $\e$-nets  
$\{x_i\}_{i=1}^{I(\e)},\, \{x'_i\}_{i=1}^{I(\e)}$ in
$B_X(p,r),\,B_{X'}(p',r)$ satisfying $|d_X(x_i,x_j)-d_{X'}(x'_i,x'_j)|<\e$, and an $\e$-net  $\{t_\ell\}_{\ell=1}^{L(\e)}$ in $(\e, \e^{-1})$, such that
\bequ \label{24a.april.11}
\big|H_X(x_i, x_j, t_\ell)-H_{X'}(x'_i, x'_j, t_\ell)\big|<\e,\;\,\textrm{ for all }\; 1\leq i,j \leq I(\e),\, 1\leq \ell \leq L(\e).
\eequ
\end{definition}

In this section we prove the following result.

\begin{theorem} \label{heat_continuity_pointed} 
For any $\e>0$, there exists $\delta=\delta(\e)>0$ such that, if
$X, X' \in {\overline{\frak M \frak M}}_p (n,\Lambda,D)$ satisfy $d_{pmGH}(X, X') <\delta$, then
their heat kernels $H_X$ and $H_{X'}$ 
are $\e$-close.
\end{theorem}

Note that the choice of $\delta=\delta(\e)$ above is uniform on 
$\overline{{ \frak M \frak M}}_p (n,\Lambda,D)$.
Theorem \ref{heat_continuity_pointed} is a direct consequence of Theorem \ref{heat_continuity} applying to the ball $B_X(p,r)$.
The proof of Theorem \ref{heat_continuity}, based on \cite{Fuk_inv}, is rather long and will occupy the rest of this section.
%
%
%
%


\subsection{Spectral estimates in ${\overline{\frak F\frak M \frak M}}_p$ and ${\overline{\frak M \frak M}}_p$}
\label{eigenfunctions_estimates}

First, we recall the Weyl-type estimate for the eigenvalues. To this end,
recall the counting function, ${{N}}_X(E)$ of $X$, which is the number of the
eigenvalues of $-\Delta_X$, counted with their multiplicities,  that are not larger than $E$.
\begin{lemma} \label{BBG}
There exist constants $c,\,C>0$ such that, for any 
$X \in \overline{{ \frak M \frak M}}_p (n,\Lambda,D)$,
\bequ \label{30.april.11}
\la_{j}^X \geq c \,j^{2/n},\; j=0, 1,  \dots,\quad
 {{N}}_X(E) \leq 1+ C E^{n/2}.
\eequ
\end{lemma}
\begin{proof}
By \cite{BBG}, there exists $c>0$ such that the estimate (\ref{30.april.11}) is valid
for all manifolds in 
${ \frak M \frak M} (n,\Lambda,D)$. 
Since the measured GH-convergence implies the convergence of
eigenvalues, see \cite{Fuk_inv}, the inequality (\ref{30.april.11}) for eigenvalues
remains valid for any space in the closure ${\overline { \frak M \frak M}} (n,\Lambda,D)$. The second inequality in (\ref{30.april.11}) follows from the first.
\end{proof}

We start with uniform spectral estimates in ${\overline{\frak F\frak M \frak M}}_p (n,\Lambda,D)$ and ${\overline{\frak M \frak M}}_p (n,\Lambda,D)$. 
Since these estimates are independent of point $p$, we often omit the sub-index $p$.
Our first goal are the Zygmund-type estimates
of the eigenfunctions on the sequence $FM_k \in {\frak F\frak M \frak M}$ collapsing to $Y\in \overline{{\frak F\frak M \frak M}}$. However, since the
injectivity radii on $FM_k$ are not uniformly bounded from below, we work by lifting the eigenfunctions
to the tangent space $T(FM_k)$ where the injectivity radii are uniformly bounded below, see Section \ref{subsection-basic-property}.
We use the notation $\phi$ for the eigenfunctions on $X\in {\overline{\frak M \frak M}}$, and $\tilde \phi$ for the eigenfunctions on $Y\in {\overline{\frak F\frak M \frak M}}$.


\begin{proposition} \label{alpha-uniform}
There are $C_F=C_F(n, \Lambda, D)>0$ and $s_F=s_F(n, \Lambda, D) >0$ 
such that for any 
 $(FM, \mu_{FM}) \in {\frak F \frak M \frak M}$, and $\tilde\phi$ satisfying 
\bequ \label{19.12.03}
\Delta_Y \tilde \phi = \la \tilde \phi, 
\eequ
we have 
\bequ \label{19.12.02}
\|\tilde \phi\|_{C^{0, 1}(FM)} \leq C_F (1+\la^2)^{s_F/2} \|\tilde \phi\|_{L^2(FM, \mu_{FM})}.
\eequ
\end{proposition}

\begin{proof}
For $y \in FM$, let $\hat B_y(r) \subset T_y(FM)$ be a closed  ball  of radius
$r < \pi/\Lambda_F$, where $\Lambda_F$ is the sectional curvature bound for $FM$, see \eqref{2.15.01.2012}. Let $\tilde h$ denote the metric on $FM$, and let
$\hat  h:=(\exp_y)^* \tilde h$
denote the lifted metric on $\hat B_y(r)$.
Recall that $\hat B_y(r)$ with the lifted metric $\hat h$ has injectivity radii uniformly bounded from below.
Then by \cite{AKKLT}, there exist uniform constants $C_h,r_h>0$ such that there are 
 harmonic coordinates
 $\Psi:\hat  B_y(r_h)\to \R^\ell$ with $\ell=\dim(FM)$, in which
 the metric tensor  satisfies 
 \beq \label{19.12.05}
C_h^{-2}I\leq \Psi_*\hat   h\leq C_h^{2}I,\quad
\|\Psi_*\hat  h\|_{C^2_*(\Psi(\hat  B_y(r_h)))} \leq C_h.
\eeq
In the following, let $R_0=r_h/2 C_h$. Then  the ball $B(R_0) \subset \R^\ell$
satisfies $B(R_0)\subset \Psi(\hat B(r_h/2))$. 
Below, we identify functions on $B(R_0)\subset \mathbb{R}^{\ell}$ and $\Psi^{-1}(B(R_0))$ via $\Psi$, and
use the Lebesgue spaces $L^{p}(B(R_0))$ and 
Sobolev spaces $W^{k,p}(B(R_0))$ which are defined
using the usual Lebesgue measure on $\R^\ell$. We note that as the metric $\hat h$ satisfies (\ref{19.12.05}), the norms of these spaces are equivalent to the norms in the spaces defined using the metric $\hat h$ on $T_y(FM)$.

As $\tilde \phi$ satisfies the eigenvalue equation, the lifted function $\hat \phi:=(\exp_y)^{*} \widetilde{\phi}$ satisfies
\bequ \label{19.12.06}
\Delta_{\hat h} \hat \phi=\la \hat \phi \quad \hbox{in}\,\,  B(R_0).
\eequ
Let us consider radii $R_k=(1-k/(8\ell))R_0$ and $q_k$ given by
$1/q_k=1/2-k/(2\ell)$, $k=0,1,\dots, \ell$.
Then in the harmonic coordinates above, the standard local elliptic regularity estimates, see e.g. \cite[Thm.\ 9.4.1]{Kry}, yield that
\beq \label{05.13.02}
\|\hat \phi\|_{W^{2,q_k}(B(R_k))} 
\leq c_{k, \ell} (1+\lambda^2)^{1/2} \|\hat \phi\|_{W^{2,q_{k-1}}(B(R_{k-1}))},
\eeq
where the constant $c_{k, \ell}>0$
is uniform for $FM \in {\frak F \frak M \frak M}$. 
Combining the estimates (\ref{05.13.02}), for $k=0,1,\dots, \ell,$ together
with the Sobolev embedding theorem, we obtain that
\bequ \label{05.13.01}
 \|\hat \phi\|_{ C^{0,1}(B(R_0/2))} 
\leq C(\ell) (1+\lambda^2)^{s(\ell)/2} \|\hat \phi\|_{L^2(B(R_0))}.
\eequ
By \cite[Lemma 1.2]{K93}, The right-hand side $\|\hat \phi\|_{L^2(B(R_0))}$ of \eqref{05.13.01}
 is controlled by $\|\tilde \phi\|_{L^2(FM, \mu_{FM} )}$, namely
\begin{eqnarray}
\|\hat \phi\|^2_{L^2(\hat B(r), dV_{\hat h})} 
&\leq& 
\frac{V_{\hat h}(\hat B(4r))  }{ V_{\tilde h}(B(y, r)) }\int_{B(y, r)} |\tilde \phi|^2\,  dV_{\tilde h} \\
&\leq&  V_{\hat h}(\hat B(4r)) 
\frac{ V_{\tilde h}(FM)}{ V_{\tilde h}(B(y, r)) }
\int_{FM} |\tilde \phi|^2\, d\mu_{FM} . \nonumber
\end{eqnarray}
Recall that $d\mu_{FM}=dV_{\tilde h}/V_{\tilde h}(FM)$ is the normalized Riemannian volume.
The volume $V_{\hat h}(\hat B(4r))$ is bounded above due to \eqref{19.12.05}, and $V_{\tilde h}(B(y, r))/V_{\tilde h}(FM)$ is bounded below due to the Bishop-Gromov volume comparison theorem. 

At last, observe that for any Lipschitz function $u$ on $FM$, it satisfies that
\bequ
\|u\|_{C^{0,1}(FM)} \leq  \max_{y \in FM}\|(\exp_y)^* u\|_{C^{0,1}(B(R_0/2))},
\eequ
which concludes the proof.
\end{proof}

\begin{lemma} \label{c-beta}
There are $C_F>0$ such that,
for any $(X, \mu_X) \in \overline{{  \frak M \frak M}}$ and $(Y, \mu_Y) \in \overline{{ \frak F \frak M \frak M}}$, the
estimate 
 (\ref{19.12.02}) remains valid for the eigenfunctions $\phi$ and $\tilde \phi$, correspondingly.
\end{lemma}
\begin{proof}
It suffices to consider 
$Y \in \overline{{\frak F \frak M \frak M}} \setminus {{\frak F \frak M \frak M}}$.
Let $(Y,  \mu_Y)$ be the limit of 
$(F M_k, \tilde \mu_k) \in {\frak F \frak M \frak M}$ in the measured GH topology and let
$$
\e_k=d_{pmGH}(Y, F M_k) \to 0, \;\textrm{ as } k\to \infty.
$$
Denote by $\tilde f_k:F M_k \to Y$ a regular fibration of $Y$ which enjoys $\varepsilon$-approximation properties (\ref{GH-approximation}).
Let $\tilde \phi$ be an eigenfunction of $-\Delta_Y$ with $\|\tilde\phi\|_{L^2(Y)}=1$, corresponding to an eigenvalue
$\la <E$.
By \cite{Fuk_inv}, there exists $\delta_k >0$ such that $\lim_{k\to \infty}\delta_k=0$ 
and there are
\begin{equation}
\tilde \phi_k \in \tilde {\mathcal L}_k(\la-\delta_k, \la+\delta_k),\quad \|\tilde \phi_k\|_{L^2(F M_k)}=1, 
\end{equation}
satisfying
\begin{equation} \label{03.09.2011.6}
\|\tilde f^*_k(\tilde \phi)-\tilde \phi_k\|_{L^2(F M_k)} <\delta_k.
\end{equation}
Here and later, for $a<b,\, a, b \notin \hbox{spec}(-\Delta_{F M_k})$, we denote by  
${\tilde {\mathcal L}}_k(a, b)={\tilde {\mathcal L}}_{F M_k}(a, b)$ li combination
of the  eigenfunctions of $-\Delta_{F M_k}$
corresponding to the eigenvalues in the interval $(a, b)$.

Let $\tilde f'_k$ be a (non-necessary continuous) right inverse to $\tilde f_k$, i.e., 
$\tilde f_k \circ \tilde f_k'= \hbox{id}_Y.$
Denote
\bequ
\tilde \phi_k'(y)=\tilde \phi_k(\tilde f_k'(y)), \quad y \in Y.
\eequ
Observe that
\bequ \label{05.09.2011.3}
\|\tilde f_k^*(\tilde \phi_k')-\tilde \phi_k\|_{L^\infty(FM_k)} \leq c (1+E^2)^{s/2}\e_k,
\eequ
where $c=c(E), s=s(E)$ are uniform on ${\frak F\frak M}{\frak M}$. 
Indeed, for any $y_k, y'_k \in \tilde f_k^{-1}(y)$ formula (\ref{GH-approximation})
yields $d_{FM_k}(y_k, y'_k) \leq c \e_k$.
This, together with the uniform Lipshitz continuity (\ref{19.12.02}) of functions $\tilde\phi_{FM_k}$, yields the inequality (\ref{05.09.2011.3}),
which in particular gives
$$
\|\tilde f_k^*(\tilde \phi_k')-\tilde \phi_k\|_{L^2(FM_k)} \leq c(1+E^2)^{s/2} \e_k.
$$
Combining this with (\ref{03.09.2011.6}) yields
\beq \label{03.09.2011.5}
\|\tilde f_k^*(\tilde \phi_k'-\tilde\phi) \|_{L^2(FM_k)} \leq c (1+E^2)^{s/2}(\delta_k+ \e_k).
\eeq

On the other hand, there exists a subsequence, $k=k(p)$,  
and a function $\tilde \phi' \in C^{0,1}(Y)$ such that, for any $y \in Y$,
$$
\lim_{k\to \infty} \tilde \phi'_k(y)=\tilde \phi'(y) .
$$
Indeed, choosing a dense subset $\{y_p \}_{p=1}^\infty \subset Y$ and using
the diagonalization procedure we find $k=k(n),\, k(n) \to \infty$ as $n \to \infty$
such that for all $p\in \Z_+$ there exist limits
$$
\tilde \phi'(y_p)=\lim_{n \to \infty}\tilde\phi'_{k(n)}(y_p).
$$
Using
the uniform Lipschitz bound (\ref{19.12.02}), the estimate $d_{FM_k}(y_k, y'_k) \leq c\e_k$,
for $y_k, y'_k \in \tilde f_k^{-1}(y)$,
 and the definition of the GH distance, 
we see that
\bequ \label{almost-Lip}
|\tilde \phi'_{k(n)}(y)-\tilde \phi'_{k(n)}(y')| \leq 
c (1+E^2)^{s/2}\left(\e_{k(n)}+d_Y(y, y') \right),\quad y, y' \in Y.
\eequ
Thus,
we extend $\tilde \phi'$ from $\{y_p: \ p\in \Z_+\} \subset Y$ to the whole space $Y$ so that 
\bequ \label{03.09.2011.9}
|\tilde \phi'(y)-\tilde \phi'(y')| \leq c (1+E^2)^{s/2} d_Y(y, y'),\quad \hbox{for all } y, y' \in Y.
\eequ
This shows that $\tilde\phi'$ satifies the desired estimate and we need to show $\tilde \phi' =\tilde \phi$.

From \eqref{almost-Lip} and \eqref{03.09.2011.9}, one can show that the pointwise convergence of the subsequence $\tilde \phi'_k$ to $\tilde \phi'$ is uniform, i.e.,
\bequ
\lim_{k\to \infty}\|\tilde \phi'_k -\tilde \phi'\|_{L^{\infty}(Y)} =0,
\eequ
and consequently
$\lim\limits_{k\to \infty}\|\tilde f_{k}^{\ast} (\tilde \phi'_k -\tilde \phi' )\|_{L^{\infty}(FM_k)} =0.$
Hence, together with \eqref{03.09.2011.5}, we have
$$
\lim_{k\to \infty}\|\tilde f_{k}^{\ast} (\tilde \phi' -\tilde \phi )\|_{L^2(FM_k)} =0,
$$
which gives 
$\|\tilde \phi' -\tilde \phi \|_{L^{2}(Y)} =0$
by definition \eqref{measured-def}.
Hence $\tilde \phi' =\tilde \phi$ since they are both continuous.

At last, due to \eqref{4.19.01.2012} and \eqref{26.11.2011.2}, the estimate for $ \overline{ \frak M \frak M}$ follows from that for
$\overline{\frak F \frak M \frak M}$.
\end{proof}




For $a<b,\, a, b \notin \hbox{spec}(-\Delta_{Y})$, we denote by  
${\tilde {\mathcal L}}_Y(a, b)$, for $Y \in \overline{\frak F\frak M \frak M}$,
 the li combination
of the  eigenfunctions of $-\Delta_{Y}$
corresponding to the eigenvalues in the interval $(a, b)$.

\begin{corollary} \label{from_infinity}
Let $\tilde u \in {\widetilde {\mathcal L}}_Y(a, b),\, a, b
\notin \hbox{spec}(-\Delta_Y)$ or $ u \in { \mathcal L}_X(a, b),\, a, b 
\notin \hbox{spec}(-\Delta_X)$. Then 
(\ref{19.12.02}) remains valid for $\tilde u,\, u$ (with $a,b$ instead of $\la$).
\end{corollary}

\subsection{Spectral convergence on $\overline{{\frak F \frak M \frak M}}_p$ and
 $\overline{{ \frak M \frak M}}_p$}

Let $(X_k,\mu_k)\in \overline{{ \frak M \frak M}}_p$ be a sequence converging to $(X,\mu_X)$ in the measured Gromov-Hausdorff topology.
Then as stated in Theorem \ref{thm:fib}, there exist Riemannian manifolds $(Y_k, \tilde \mu_Y^k) \in \overline{{\frak F \frak M \frak M}}_p$ converging to $(Y,\tilde \mu_Y)$ in the measured Gromov-Hausdorff topology, and regular $\e_k$-approximations $\tilde f_k:Y_k\to Y$, $f_k:X_k\to X$, such that $X_k=Y_k/O(n)$, $X=Y/O(n)$, and $\tilde f_k$ are $\e_k$-Riemannian submersions.  
 
For $a<b$, $a,b\notin \hbox{spec}(-\Delta_Y)$, 
we denote by ${\widetilde{\mathcal L}}_Y(a, b)$ the li combinations of the eigenfunctions of $-\Delta_Y$ corresponding to eigenvalues in the interval $(a,b)$, and define similarly ${\mathcal L}_X(a, b)$ as the li combinations of the eigenfunctions of $-\Delta_X$.
Using pullback notations, denote
\beq
&&{\widetilde {\mathcal L}}_k^*(a, b)= \tilde f_k^*\left({\widetilde{\mathcal L}}_Y(a, b) \right) \subset C^{0, 1}(Y_k), \label{LstarY} \\
&&{\mathcal L}_k^*(a, b)= f_k^*\left({\mathcal L}_X(a, b) \right) 
\subset C^{0, 1}(X_k). \label{LstarX}
\eeq
Similarly, denote by ${\widetilde{\mathcal L}}_k(a,b),\mathcal{L}_k(a,b)$ the li combinations of the eigenfunctions of the Laplacian on $Y_k,X_k$, respectively.

Assuming $a, b \notin  \hbox{spec}(-\Delta_{Y})$. In the case of $Y_k=FM_k\in {\frak F \frak M \frak M}_p$, it follows from \cite{Fuk_inv} that for large
$k$, we have $a, b \notin \hbox{spec}(-\Delta_{FM_k})$, 
and
\bequ \label{19.12.08}
d_{L^2(FM_k)}\left({\mathcal B}({\widetilde{\mathcal L}}_k^*(a,b)), \,
{\mathcal B}({\widetilde{\mathcal L}}_k(a,b))  \right) 
\to 0, \quad \hbox{as}\,\,\, k\to \infty,
\eequ
with respect to the normalized Riemannian measure on $FM_k$,
where ${\mathcal B}({\mathcal Z})$ stands for the unit ball in ${\mathcal Z}$, and $d_{L^2(FM_k)}$ is the Hausdorff distance in $L^2(FM_k)$.

\begin{lemma} \label{2023.2.9}
In the case of $Y_k=FM_k\in {\frak F \frak M \frak M}_p$, the convergence in (\ref{19.12.08}) is valid in $C^{\alpha}(FM_k)$, for any $0\leq \alpha<1$.
\end{lemma}

\begin{proof}
We start with the case $\alpha=0$. Assume the convergence does not hold in $C^0(FM_k)$. Then there exists
${\widetilde {\delta}}>0$ such that,
for any $k$, either there
exists $\tilde u_k \in {\mathcal B}({\widetilde{\mathcal L}}_k(a,b))$ such that
\bequ \label{04.09.2011.1}
\|\tilde u_k-\tilde v^*\|_{C(FM_k)}> 
{\widetilde {\delta}},\quad \hbox{for all} \,\,  
\tilde v^*\in  {\mathcal B}({\widetilde{\mathcal L}}_k^*(a,b)),  
\eequ
or there exists $\tilde v^*_k \in  {\mathcal B}({\widetilde{\mathcal L}}_k^*(a,b))$ 
such that 
$$
\|\tilde v_k^*-\tilde u\|_{C(FM_k)}> {\widetilde \delta},
\quad \hbox{for all}\, \,
\tilde u \in {\mathcal B}({\widetilde{\mathcal L}}_k(a,b)).
$$
Without loss of generality we consider the former situation.

On the other hand, by (\ref{19.12.08}), for all $k$ there exists
$\tilde v_k^* \in {\mathcal B}({\widetilde{\mathcal L}}_k^*(a,b))$ 
 such
that
\bequ \label{04.09.2011.2}
\|\tilde u_k-\tilde v_k^*\|_{L^2(FM_k)} \to 0, \quad \textrm{as } k\to \infty.
\eequ
However, by (\ref{04.09.2011.1}), there exists $y_k \in FM_k$ such that
$|\tilde u_k(y_k)-\tilde v_k^*(y_k)|> \tilde{\delta}$. 
Since functions in ${\mathcal B}({\widetilde{\mathcal L}}_k(a,b))$ and ${\mathcal B}({\widetilde{\mathcal L}}_k^*(a,b))$ are uniformly bounded in $C^{0,1}$-norm due to (\ref{19.12.02}) and (\ref{1.22.01.2012}), one can find a neighborhood $\tilde U_k$ of $y_k$ with uniform radius such that $|\tilde u_k(y)-\tilde v_k^*(y)|> \tilde{\delta}/2$ for all $y\in \tilde U_k$.
This implies that $\|\tilde u_k-\tilde v_k^*\|_{L^2(FM_k)}$ is uniformly bounded away from $0$, which is a contradiction to $(\ref{04.09.2011.2})$. This proves the convergence in $C^0$.

Due to (\ref{19.12.02})
and (\ref{1.22.01.2012}),
$$
d_{C^{0,1}(FM_k)}\left({\mathcal B}({\widetilde{\mathcal L}}_k^*(a,b)), \,
{\mathcal B}({\widetilde{\mathcal L}}_k(a,b))  \right)$$ 
is bounded. 
Then standard interpolation gives the convergence in $C^{\alpha}$.
\end{proof}

The following lemma deals with the general case $Y_k \in \overline{{\frak F \frak M \frak M}}_p$.

\begin{lemma} \label{everywhere_C}
\begin{enumerate}
\item
Assume $Y \in \overline{{\frak F \frak M \frak M}}_p,\, a, b \notin \hbox{spec}(-\Delta_Y)$,
and $Y_k \in \overline{{\frak F \frak M \frak M}}_p$ converge to $Y$ in the measured GH-topology.
Then for any $0\leq \a <1$, 
\bequ \label{2.22.01.2012}
d_{C^\a(Y_k)}\left(
{\mathcal B}({\widetilde{\mathcal L}}_k^*(a, b)),\,
{\mathcal B}({\widetilde{\mathcal L}}_k(a, b)) \right) \to 0,\quad \hbox{as}
\,\, k \to \infty.
\eequ 
\item
Assume $X \in \overline{{ \frak M \frak M}}_p,\, a, b \notin \hbox{spec}(-\Delta_X)$
and $X_k \in \overline{{ \frak M \frak M}}_p$ converge to $X$ in the measured GH-topology.
Then for any $0\leq \a <1$, 
\bequ \label{04.09.2011.5}
d_{C^\a(X_k)} \Big( {\mathcal B}({\mathcal L}_k^*(a, b)),\,
{\mathcal B}({\mathcal L}_k(a, b)) \Big) \to 0,\quad \hbox{as}
\,\, k \to \infty.
\eequ 
\end{enumerate}
\end{lemma}

\begin{proof}
{\it 1.} As earlier, we start with the case $\a=0$. Assuming the opposite, there is $\tilde \delta>0$ such that
\bequ \label{3.22.01.2012}
d_{C(Y_k)}\left({\mathcal B}({\widetilde{\mathcal L}}_k^*(a, b)),\,
{\mathcal B}({\widetilde{\mathcal L}}_k(a, b)) \right) >\tilde \delta,
\quad  \hbox{for\,all}\,\, k.
\eequ
Approximate, in the measured GH-topology, $(Y_k, \tilde\mu_Y^k)$ 
by $(FM_k, \tilde \mu_k) \in 
\frak F \frak M \frak M_p$.
Let $\Delta_{kk}$ stand for the Laplacian on $FM_k$. 
Denote by
 ${\widetilde{\mathcal L}}_{kk}(a, b) \subset L^2(FM_k)$  the eigenspace 
of the Laplacian $-\Delta_{kk}$ on $FM_k$ corresponding to the eigenvalues in $(a, b)$, and by
$$
{\widetilde{\mathcal L}}_{kk}^*(a, b) :=
f^*_{kk}\left({\widetilde{\mathcal L}}_k(a, b)\right),
$$
where
$f_{kk}: FM_k \rightarrow Y_k$ is the regular approximation
described in Theorem \ref{thm:fib} with $\e_k$ in (\ref{1.22.01.2012}) changed into
$\e_{kk}$.
Moreover, we can assume that
\bequ \label{5.22.01.2012}
d_{L^2(FM_k)}\left({\mathcal B}({\widetilde{\mathcal L}}_{kk}^*(a, b)),\,
{\mathcal B}({\widetilde{\mathcal L}}_{kk}(a, b)) \right) < \frac1k.
\eequ
Let
$
\tilde f_k: Y_k \to Y
$
be a $\tilde\e_k$-Riemannian submersion, satisfying (\ref{1.22.01.2012}) and
providing the measured GH-convergence of $Y_k$ to $Y$, which is guaranteed
by Theorem \ref{thm:fib}.
We can assume $\tilde f_k$ to be an $O(n)$-map.

Consider
$$
\tilde f_{kk}=  \tilde f_k \circ f_{kk}:\, FM_k \to Y,
$$
which provides a regular $\tilde \e_{kk}$-Riemannian submersion of $FM_k$ to 
$Y,\, \tilde \e_{kk}\to 0,$ as $k \to \infty$. Therefore,
$$
\lim_{k\to \infty}d_{L^2(FM_k)} \left({\mathcal B}({\widetilde{\mathcal L}}_{kk}(a, b)),\,  
{\mathcal B}({\widetilde{\mathcal L}}^*_{kkk}(a, b))\right)= 0,
$$
where ${\widetilde{\mathcal L}}_{kk}(a, b)$ is the eigenspace of $-\Delta_{kk}$ on $FM_k$ corresponding
to the interval $(a, b)$, and
\bequ \label{6.22.01.2012}
{\widetilde{\mathcal L}}^*_{kkk}(a, b) := 
\tilde f^*_{kk}({\widetilde{\mathcal L}}_Y(a, b))=
 f^*_{kk}({\widetilde{\mathcal L}}^*_{k}(a, b)).
\eequ
This, together with (\ref{5.22.01.2012}) implies that
$$
\lim_{k\to \infty}
d_{L^2(FM_k)} \left({\mathcal B}({\widetilde{\mathcal L}}_{kk}^*(a, b)),\, 
{\mathcal B}({\widetilde{\mathcal L}}^*_{kkk}(a, b))\right) = 0.
$$
On the other hand, (\ref{3.22.01.2012}) together with definition 
(\ref{6.22.01.2012}) yields that, for large $k$,
$$
d_{C(FM_k)} \left({\mathcal B}({\widetilde{\mathcal L}}_{kk}^*(a, b)),\, 
{\mathcal B}({\widetilde{\mathcal L}}^*_{kkk}(a, b))\right) > \tilde \delta/2.
$$
Then a similar argument as the proof of Lemma \ref{2023.2.9} shows that the above two inequalities lead to a contradiction.
This proves (\ref{2.22.01.2012}) for $\a=0$.

To obtain the result for any $0\leq \a <1$, we use again the fact that, 
due to  Lemma \ref{c-beta}
and (\ref{1.22.01.2012}), functions in ${\mathcal B}({\widetilde{\mathcal L}}_k(a,b))$ and ${\mathcal B}({\widetilde{\mathcal L}}_k^*(a,b))$ are uniformly bounded in $C^{0,1}(Y_k)$. Then the convergence in $C^{\alpha}$ follows from interpolation.

\smallskip
\noindent {\it 2}. 
Since $\tilde f_k$ is an $O(n)$-map,
(\ref{2.22.01.2012}) remains valid with  $O(n)$-invariant spaces
${\widetilde{\mathcal L}}_{k,O}^*(a, b)$ and $
{\widetilde{\mathcal L}}_{k, O}(a, b)$ instead of
${\widetilde{\mathcal L}}_k^*(a, b),\,
{\widetilde{\mathcal L}}_k(a, b)$. 
Denote $\pi_k: Y_k \to X_k=Y_k/ O(n)$. Then (\ref{2.22.01.2012}) for
the $O(n)$-invariant functions,
 together with (\ref{3.19.01.2012}) and (\ref{7.22.01.2012})
proves (\ref{04.09.2011.5}).
\end{proof}

\begin{remark} \label{point}{\rm 
If $X= \lim_{pmGH} X_k= \{point, 1\}$, then $L^2(X)=\R$ and $\Delta_X=0$. Thus,
the only eigenvalue is $\la_0=0$ with the corresponding eigenfunction $1$.
Due to (\ref{30.april.11}), (\ref{04.09.2011.5}) remains trivially valid with,
when $k$ is sufficiently large,
with
${\widetilde{\mathcal L}}_k(a, b),\, {\widetilde{\mathcal L}}^*_k(a, b)$ consisting
of constant functions,  if $0 \in (a, b)$, and of only $0$-function if $0 \notin (a, b)$.}
\end{remark}

Denote by $f'_k: X \to X_k$ the almost right inverse to $f_k$, i.e., $f_k, f'_k$
satisfy \eqref{GH-approximation} (see e.g. \cite[Lemma 2.5]{Fuk_inv}), and for $x_k\in X_k,\,x\in X$,
\bequ \label{1.25.01.2012}
d_k(x_k,\, f'_k \circ f_k(x_k)),\,
d_X(x,\, f_k \circ f'_k(x)) \to 0, \quad \hbox{as}\,\, k \to \infty.
\eequ
Then using (\ref{04.09.2011.5}) and the uniform $C^{0, 1}$-boundedness of 
${\mathcal B}({\mathcal L}_k(a, b)),\, {\mathcal B}({\mathcal L}(a, b))$, we obtain

\begin{corollary} \label{in_X} 
Let $(X_k, \mu_k), (X,\mu_X) \in \overline{\frak M \frak M}_p$ be such that
$(X_k, \mu_k)$ converges to $(X,\mu_X)$ in the measured GH-topology.
Then there exists $\sigma_k \to 0$, as $k \to \infty$, such that
\bequ \label{4.25.01.2012}
d_{L^\infty(X)} \left({\mathcal B}({\widetilde{\mathcal L'}}_{k}^*(a, b)),\, 
{\mathcal B}({\mathcal L}_X(a, b))\right) < \sigma_k,
\eequ
where
$$
{\widetilde{\mathcal L'}}_{k}^*(a, b)= (f'_k)^* ({\mathcal L}_{k}(a, b)).
$$
Moreover, there is $c>0$, such that if  $d_X(x, x') <\sigma_k$ then for
all $u^*_k \in {\mathcal B}({\widetilde{\mathcal L'}}_{k}^*(a, b))$ we have
\bequ \label{3.25.01.2012}
|u^*_k(x) -u^*_k(x')| <c \sigma_k.
\eequ
\end{corollary}
\begin{proof}
The first claim follows from the uniform $C^{0,1}$-estimate in
${\mathcal L}_X(a, b)$ together with (\ref{04.09.2011.5})
for $\a=0$
and (\ref{1.25.01.2012}).
To prove (\ref{3.25.01.2012}) we also use the uniform $C^{0, 1}$-estimate
in ${\mathcal L}_X(a, b)$
 together with
(\ref{4.25.01.2012}). 
\end{proof}


Our next goal is to obtain continuity of the eigenfunctions  in $L^\infty$-norm, 
with respect to the measured GH distance which is uniform on ${\overline{\frak M \frak M}}_p$.

\begin{definition} \label{compact_subsets}
{\rm Let $(a_\ell, b_\ell),\, \ell=1, \dots, L$, be a finite collection of open intervals,
 and $d_\ell, \ell=1, \dots, L,$ be positive numbers. Denote the set 
 ${(a_\ell, b_\ell, d_\ell)}_{\ell=1}^L$ by ${\mathcal I}$.
Then we define
\ba
{\overline{{ \frak M \frak M}}}_{\mathcal I}= 
\Big\{X \in \overline{{ \frak M \frak M}}_p:\,\,
d\big(\hbox{spec}(X), \{a_\ell\} \big)\geq d_\ell\hbox{ and }
 d\big(\hbox{spec}(X), \{b_\ell\} \big) \geq d_\ell \Big\}.
\ea}
\end{definition}
Note that ${\overline{{ \frak M \frak M}}}_{\mathcal I}$ is closed and thus compact with respect to the measured GH-topology.

\begin{corollary} \label{uniform_C}
Let ${\mathcal I}={(a_\ell, b_\ell, d_\ell)}_{\ell=1}^L$ and ${\mathcal J}={(a_\ell, b_\ell, 
 d'_\ell)}_{\ell=1}^L,\, d_\ell >d'_\ell$. Then,
for any $\e>0$ there is $\sigma =\sigma_{{\mathcal I}, {\mathcal J}}(\e)>0$ such that, if 
$X \in {\overline{{ \frak M \frak M}}}_{\mathcal I}$ and
$d_{mGH}(X, X') < \sigma$, then
 $X' \in {\overline{{ \frak M \frak M}}}_{\mathcal J}$ 
and $X,X'$ satisfy 
$$
{\rm dim}({\mathcal L}_X(a_\ell, b_\ell))={\rm dim}({\mathcal L}_{X'}(a_\ell, b_\ell)):=
n(\ell).
$$
 Moreover, there are measurable maps
$f:X \to X', \, f':X' \to X$, which are $\e$-approximations satisfying \eqref{GH-approximation}
and
(\ref{1.25.01.2012}),
such that
\beq \label{2.25.01.2012}
& &d_{L^\infty(X)}\Big({\mathcal B}(f^*({\mathcal L}_{X'}(a_\ell, b_\ell))),\,
{\mathcal B}({\mathcal L}_{X}(a_\ell, b_\ell))  \Big) <\e, 
\\ \nonumber
& &d_{L^\infty(X')}\Big({\mathcal B}((f')^*({\mathcal L}_{X}(a_\ell, b_\ell))),\,
{\mathcal B}({\mathcal L}_{X'}(a_\ell, b_\ell))  \Big) <\e.
\eeq
In addition,
 for any $\ell=1, \dots, L$, if 
$\phi_{i,\ell},\, i=1, \dots, n(\ell),$ and
$\phi'_{i,\ell},\, i=1, \dots, n(\ell),$
form an
 orthonormal eigenfunction basis
in ${\mathcal L}_X(a_\ell, b_\ell), \, {\mathcal L}_{X'}(a_\ell, b_\ell)$, correspondingly,
then there exist  orthogonal matrices 
\ba
{\Bbb U}_\ell =[u_{ij}]_{i, j=1}^{n(\ell)}\in O(n(\ell)), \quad
{\Bbb U}'_\ell =[u'_{ij}]_{i, j=1}^{n(\ell)}\in O(n(\ell)),
\ea 
such that
$\phi^{*}_{i, \ell}= \sum_{j=1}^{n(\ell)}u_{i j} f^*(\phi'_{j, \ell})$
satisfy
\beq \label{04.09.2011.11}
\| \phi_{i, \ell}-\phi^{*}_{i, \ell}\|_{L^\infty(X)}< \e,\quad i=1, \dots, n(\ell),\,\,
\ell=1, \dots, L.
\eeq
Similar result is valid for $(f')^*(\phi_{i, l})$ if we use ${\Bbb U}'_\ell$.
\end{corollary}

\begin{proof}
By compactness arguments, the first statement of the corollary 
follows immediately from Theorem 0.4(A) in
\cite{Fuk_inv}.
To prove the second statement, in particular, (\ref{2.25.01.2012}), assume 
that there are $X_k,\, \hat X_k$  and $\delta_k \to 0$ satisfying
$d_{mGH}(X_k, \hat X_k) < \delta_k$, but (\ref{2.25.01.2012}) is not valid.
Without loss of generality, we can assume that $X_k, \,\hat X_k \to X$, with respect
to the measured GH-topology, and $f_{k} : X_k \to \hat X_k,\,f'_{k} : \hat X_k \to X_k$
which provide $\delta_k$-approximations, are of the form
$$
f_{k}=\hat f_{kk} \circ f'_{kk}\,\quad f'_{k}= f_{kk} \circ  \hat f'_{kk},
$$
where
$$
f_{kk}: X \to X_k,\;\,f'_{kk}: X_k \to X \quad \hbox{and} \quad
\hat f_{kk}: X \to \hat X_k,\;\,\hat f'_{kk}: \hat X_k \to X,
$$
provide $\delta_k/4$-approximations, in the sense of \eqref{GH-approximation},
of $X, X_k$ and $X, \hat X_k$, correspondingly. 
Note that, if $f_{kk},\, f'_{kk}$ and 
$\hat f_{kk},\, \hat f'_{kk}$ are $\delta_k/4$-almost inverse of each other in the sense
of (\ref{1.25.01.2012}), then $f_k, f'_k$ are $\delta_k$-almost inverse. 

By (\ref{4.25.01.2012}) and the triangle inequality, we have
$$
d_{L^\infty (X)}\left({\mathcal B}(f_{kk}^*({\mathcal L}_{k}(a_\ell, b_\ell))),\,
{\mathcal B}(\hat f_{kk}^*({\widehat{\mathcal L}}_{k}(a_\ell, b_\ell)))  \right) \to 0,
\quad \hbox{as}\,\, k \to \infty,
$$
where ${{\mathcal L}}_{k}(a_\ell, b_\ell)) ,\,{\widehat{\mathcal L}}_{k}(a_\ell, b_\ell)) $
stand for ${{\mathcal L}}_{X_k}(a_\ell, b_\ell) ,\,{\widehat{\mathcal L}}_{{\widehat X}_k}(a_\ell, b_\ell).$
Then, using the facts that $f_k^*=(\hat f_{kk} \circ f'_{kk})^*=(f'_{kk})^*\circ \hat f_{kk}^*$, and $f_{kk},f_{kk}'$ are almost inverse of each other, and the $C^{0,1}$-boundedness of eigenfunctions,
one can show that
$$
d_{L^\infty(X_k)}\left( {\mathcal B}({\mathcal L}_{k}(a_\ell, b_\ell)),\, {\mathcal B}(f^*_k({\widehat{\mathcal L}}_{k}(a_\ell, b_\ell)))
  \right) \to 0\quad \hbox{as}\,\, k \to \infty.
$$

To prove the last statement of the corollary, recall that, due to the measure convergence
in (\ref{eq:meas-approx}) and Lemma \ref{c-beta}, $f^*(\phi'_{j, \ell}),\, j=1,\dots, n(\ell),$ satisfy
\ba
\int_X f^*(\phi'_{j, \ell})\, f^*(\phi'_{i, \ell})d\mu_X \to  \delta_{i j},
\quad \hbox{as}\,\,\, d_{mGH}(X, X') \to 0.
\ea
Together with (\ref{2.25.01.2012}), this implies the existence of 
${\Bbb U}_\ell={\Bbb U}_\ell(X, X'), \, \ell=1,\dots,L,$ such that
for $ i=1,\dots, n(\ell),$
\ba
\|\phi_{i, \ell}- \phi^*_{i, \ell}\|_{L^\infty(X)} \to 0, \quad \hbox{as 
$d_{mGH}(X, X') \to 0$},
\ea
which yields (\ref{04.09.2011.11}).
\end{proof}

\subsection{Heat kernel convergence}
In this subsection, we consider the continuity of the heat kernel $H_X(x, y, t)$, 
with respect to the measured GH-convergence on 
$\overline{{ \frak M \frak M}}_p(n, \Lambda, D)$.
Recall that the heat kernel on $X$ can be written in terms of eigenfunctions 
and eigenvalues as
\bequ \label{21.12.01}
H_X(x, y, t)= \sum_{j=0}^\infty e^{-\la_ j t} \phi_j(x) \phi_j(y), 
\eequ
where the eigenfunctions $ \phi_j,\, j= 0, 1, \dots,$ form an orthonormal basis in $L^2(X, \mu_X)$ and $\la_0=0$ and 
$\phi_0=1$.

\begin{lemma}\label{heat_kernel_estimate}
There exist constants $C >0$ and $s_F>0$, as introduced in (\ref{19.12.02}),
such that, for any $E\geq 1$ and $X \in {\overline{\frak M \frak M}}_p (n,\Lambda,D)$, we have
\bequ \label{21.12.03}
\Big|H_X(x, y, t) -\sum_{\la_j < E} e^{-\la_ j t} \phi_j(x) \phi_j(y) \Big|
\leq   C t^{-(2s_F+\frac{n}{2}+1)} e^{-\frac12 Et},\quad t>0.
\eequ
\end{lemma}
\begin{proof}
 For any $x, y \in X$, using the $C^0$-norm estimate in Lemma \ref{c-beta} and \eqref{30.april.11}, we have
\begin{eqnarray*}
\bigg|\sum_{ \la_j \geq E} e^{-\la_ j t} \phi_j(x) \phi_j (y) \bigg| 
 &\leq&  \sum_{ \la_j \geq  E} e^{-\la_ j t} |\phi_j(x)| |\phi_j(y) | \\
&\leq& C_F^2 \sum_{ \la_j \geq  E} e^{-\la_ j t} (1+\la_j^2)^{s_F} \\
&\leq& C \sum_{\lambda_j \geq  E,\,|\lambda_j - \lambda_k|\geq 1} e^{-\la_j t} \la_j^{2s_F+\frac{n}{2}} ,
\end{eqnarray*}
where $C_F, \, s_F$ are the same constants as in \eqref{19.12.02}.
It is straightforward to verify that, for any $t>0$,
\begin{equation}
e^{-\frac12 \lambda t} \lambda^{2s_F+\frac{n}{2}} \leq C(n,s_F) t^{-(2s_F+\frac{n}{2})},\; \textrm{ for any }\lambda\geq 0.
\end{equation}
Without loss of generality, we can take $E$ to be an integer.
Then,
\begin{eqnarray*}
\bigg|\sum_{ \la_j \geq E} e^{-\la_ j t} \phi_j(x) \phi_j (y) \bigg| 
 &\leq&  C t^{-(2s_F+\frac{n}{2})} \sum_{\lambda_j \geq  E,\,|\lambda_j - \lambda_k|\geq 1} e^{-\frac12 \la_j t}, \\
&\leq&  C t^{-(2s_F+\frac{n}{2})} \sum_{j= E}^{\infty} e^{-\frac12 j t} \\
&=&  C t^{-(2s_F+\frac{n}{2})} \frac{e^{-\frac12 E t}}{1-e^{-\frac12 t}} \leq C t^{-(2s_F+\frac{n}{2}+1)} e^{-\frac12 E t}.
\end{eqnarray*}
\vspace{-2mm}
\end{proof}

As a consequence, we obtain an estimate for the continuity of the heat kernel in a particular space $X \in \overline{{\frak M}{\frak M}}_p (n,\Lambda,D)$.

\begin{corollary} \label{H_continuity}
Let $X \in \overline{{\frak M}{\frak M}}_p (n,\Lambda,D)$ and $0<\hat \e<1$.
Then there exists a constant $C>0$ such that, for any $x, y, x', y' \in X$ and $t,t'>\hat \e$, 
we have
\beq \label{13.09.2011.1}
\big|H_X(x, y, t)-H_X(x', y', t')\big|  < 
C \hat\e^{-(2s_F+\frac{n}{2}+2)} \left(d_X(x, x')+d_X(y, y')+|t-t'| \right) .
\eeq
\end{corollary}
\begin{proof}
For the part of the eigenfunction expansion of the heat kernel with $\lambda_j\geq 1$, one can follow the same argument in Lemma \ref{heat_kernel_estimate} with $E=1$ using the Lipschitz norm estimate in Lemma \ref{c-beta}. The additional order in $\hat\e^{-(2s_F+\frac{n}{2}+2)}$ comes from the derivative of $e^{-\lambda_j t}$ with respect to $t$. The lower part of the eigenfunction expansion with $\lambda_j<1$ is clearly Lipschitz with bounded Lipschitz constant by Lemma \ref{c-beta} and \eqref{30.april.11}.
\end{proof}

We now prove the continuity of the heat kernel with respect to the measured Gromov-Hausdorff topology.

\begin{theorem} \label{heat_continuity} 
For any $\hat\e>0$, there exists $\delta=\delta(\hat\e)>0$ such that, if
$X, X' \in {\overline{\frak M \frak M}}_p(n,\Lambda,D)$ satisfy $d_{mGH}(X, X') <\delta$, then
there exist
$\hat\e$-nets $\{x_i\}_{i=1}^{I(\hat \e)} \subset X$ and  $ \{x'_i\}_{i=1}^{I(\hat\e)} \subset X'$
satisfying
$|d_X(x_i,x_k)-d_{X'}(x'_i,x'_k)|<\hat \e$, 
 and
\bequ \label{45.april.11}
|H_X(x_i, x_k, t)-H_{X'}(x_i', x_k', t)| < \hat\e, \quad \hbox{for}\,\, t\in (\hat\e,\hat\e^{-1}).
\eequ
\end{theorem}

\begin{proof}
By Lemma \ref{heat_kernel_estimate}, for sufficiently small $\hat \e$ and
$E=E(\hat\e)= \hat \e^{\,-2}$, 
 the right-hand
side of  (\ref{21.12.03}) is smaller than $\hat\e/4$ for any $X \in {\overline{\frak M \frak M}}_p$ and any $t >\hat \e$. 

By (\ref{30.april.11}), for any $X\in {\overline{\frak M \frak M}}_p$ we can
choose in the interval $(-1, \hat\e^{\, -2}+4 \hat \e^{\,\beta+n/2+1})$, where $\beta >0$ is to be chosen later, 
 the 
 subintervals $(a_\ell, b_\ell),\, \ell=1, \dots, L=L(X),$ satisfying the following properties:
 \begin{itemize}
 \item[(i)]\quad 
 $\, b_\ell-a_\ell \leq \hat\e^\beta,$ 
 \item[(ii)]\quad 
 $a_{\ell+1} > b_\ell+4 \hat \e^{\,\beta+n+1}$,
 \item[(iii)]\quad 
 $\hbox{spec}(-\Delta_X) \cap (-\infty, \hat \e^{\,-2})\subset \cup_{\ell=1}^L(a_\ell, b_\ell)$,
 \item [(iv)]\quad 
 $d (\hbox{spec}(-\Delta_X),a_\ell), \,d (\hbox{spec}(-\Delta_X), b_\ell) > 
 4 \hat \e^{\,\beta+n+1},
 \, \ell =1, \dots, L$.
 \end{itemize}
 Note that this system of subintervals, associated with $X$, satisfies conditions of Definition
 \ref{compact_subsets} with $d_\ell=4 \hat \e^{\,\beta+n+1}$. Denote it by 
 ${\mathcal I}={\mathcal I}_{\hat \e}(X)$
 and consider ${\overline{{ \frak M \frak M}}}_{{\mathcal I}}$. 
 Choose $d'_\ell=3 \hat \e^{\,\beta+n+1}$
 to define the corresponding ${\mathcal J}(X)$. Then there is $\tau=\tau({\mathcal I}(X), \hat \e, \beta)$,
 such that, if 
 \ba
 X' \in {\mathcal U}_{\tau, X} ({\overline{{ \frak M \frak M}}}_{{\mathcal I}})=
 \{ X' \in \overline{ \frak M \frak M}_p:\, d_{mGH}({\overline{{ \frak M \frak M}}}_{{\mathcal I}},\, X') < \tau\},
 \ea
 then $X' \in {\overline{{ \frak M \frak M}}}_{{\mathcal J}(X)}$.
 Clearly, $\{{\mathcal U}_{\tau, X}:\, X \in {\overline{ \frak M \frak M}}_p\}$
 form an open covering of ${\overline{ \frak M \frak M}}_p$. Choose a finite
 subcovering, ${\mathcal U}_1, \dots, {\mathcal U}_N,\, N=N(\hat \e)$.  Observe that 
 \ba
{\mathcal U}_k \subset {\overline{{ \frak M \frak M}}}_{{\mathcal J}(X)},
\quad \hbox{for\,some}\,\,\, X=X_k.
 \ea 
Then, due to Corollary
\ref{uniform_C}, there
 is $\sigma_k>0$ so that if $d_{mGH}(X, X') <\sigma_k, \, X \in {\mathcal U}_k,$ 
 then equations
 (\ref{2.25.01.2012}), 
(\ref{04.09.2011.11}) are valid with $\e=\hat\e^n$ (and $d'_\ell=2\hat \e^{\,\beta+n+1}$). Since 
$
\cup_{k=1}^N {\mathcal U}_k = {\overline{{ \frak M \frak M}}}_p,
$
taking $\sigma(\hat \e)= \min \sigma_k$, we see that, for any $X, X' \in {\overline{{ \frak M \frak M}}}_p$,
if $d_{mGH}(X, X') <\sigma(\hat \e)$, equations (\ref{2.25.01.2012}), 
(\ref{04.09.2011.11}) are valid with $\e=\hat\e^n$. Moreover, there is $k$ such that $X, X' \in 
{\mathcal U}_k$.

\smallskip
Consider
\beq \label{04.09.2011.12} 
& &\bigg|\sum_{\la_j < c\hat \e^{-2}} e^{-\la_j t} \phi_j(x) \phi_j(y)-
\sum_{\la_j < c\hat \e^{-2}} e^{-\la'_j t} \phi'_j(f(x)) \phi'_j(f(y))\bigg| 
\\ \nonumber 
&=&\bigg|\sum_{\ell=1}^L \bigg(\sum_{\la_j \in (a_\ell, b_\ell)} e^{-\la_j t} \phi_j(x) \phi_j(y)-\sum_{\la'_j \in (a_\ell, b_\ell)} e^{-\la'_j t} \phi'_j(f(x)) \phi'_j(f(y)) 
\bigg)\bigg|,
\eeq
where $f: X \to X', \, f': X' \to X$ are almost inverse $\sigma$-approximations.
We analyze each term on the right side of (\ref{04.09.2011.12}) separately.

First,  observe that for $t \leq \hat\e^{\, -1},\, \la_j, \la'_j \in (a_\ell, b_\ell)$,
\ba
\max\big\{|e^{-\la_j t}-e^{-a_\ell t}|,|e^{-\la_j' t}-e^{-a_\ell t}| \big\}
\leq c(\beta) \hat \e^{\,(\beta-1)} e^{-a_\ell t}.
\ea
Using Lemma \ref{c-beta} together with (\ref{30.april.11}), we see from the previous
equation that replacing the exponents
$e^{-\la_j t}$ and $e^{-\la'_j t}$ by
$e^{-a_\ell t}$ on the right side of (\ref{04.09.2011.12}) gives
rise to an error which can be estimated by $C \hat\e^{\,(\beta -1 -n-2s_F)}$
for any  $t <\hat\e^{\,-1}$.  
Let us next choose $\beta =3+n+2s_F$.
Then there is $\hat\e_1$ such that, for   $\hat \e <\hat \e_1$, 
\ba
& &\bigg|\sum_{\ell=1}^L \sum_{\la_j \in (a_\ell, b_\ell)} \left(e^{-\la_j t}-e^{-a_\ell t} \right) \phi_j(x) \phi_j(y)
\\ \nonumber
& &-\sum_{\ell=1}^L
\sum_{\la'_j \in (a_\ell, b_\ell)} \left(e^{-\la'_j t}-e^{-a_\ell t}\right) \phi'_j(f(x)) \phi'_j(f(y)) \bigg| < \frac{\hat \e}{16}.
\ea
Next, observe that
 \ba
\sum_{i=1}^{n(\ell)} \left(\sum_{j=1}^{n(\ell)} u_{i j}^\ell  f^*(\phi'_{j, \ell})(x)\right) 
\cdot \left(\sum_{k=1}^{n(\ell)} u_{i k}^\ell  f^*(\phi'_{k, \ell})(y)\right)=
\sum_{i=1}^{n(\ell)} f^*(\phi'_{i, \ell})(x) f^*(\phi'_{k, \ell})(y).
\ea
This equality,
together with Lemma \ref{c-beta} and (\ref{GH-approximation}) with $\e=\hat \e^{\,\beta}$
 implies that, that there is $0<\hat \e_2<\hat \e_1$, such that 
 for $\hat \e < \hat \e_2$,
\ba
& &|H_X(x, y, t)-H_{X'}(f(x), f(y), t)| 
\\ \nonumber
& &\leq \frac12 \hat \e +
\sum_{\ell=1}^L \,\sum_{i=1}^{n(\ell)} e^{-a_\ell t} \big|\phi_{i, \ell}(x) \phi_{i, \ell}(y)-
\phi^*_{i, \ell}(x) \phi^*_{i, \ell}(y)\big|.
\ea
Therefore, if equation (\ref{04.09.2011.11}) is valid with $\e=\hat \e^\beta$, 
there is $\hat \e_3<\hat \e_2$ 
such that for $\hat \e <\hat \e_3 $,
\bequ \label{05.09.2011.1}
|H_X(x, y, t)-H_{X'}(f(x), f(y), t)|<\hat \e,  \;\;\textrm{ for }\, t\in (\hat\e,\hat\e^{-1}).
\eequ


At last, choosing $\{x_i\}_{i=1}^{I(\hat \e)} \subset X$ to be an $\hat \e/2$-net in $X$, we see from 
(\ref{GH-approximation}) that, for 
$\delta <\delta(\hat \e)$, the set $\{f(x_i)\}_{i=1}^{I(\hat \e)} \subset X$ is $\hat \e$-dense in $X'$.
Due to the compactness of ${\overline{\frak M \frak M}}_p$, $I(\hat \e)$ can be chosen uniformly
on ${\overline{\frak M \frak M}}_p$.
\end{proof}

\section{From the local spectral data to the metric-measure structure}\label{metric}

In this section, we consider the inverse spectral  problem on ${\overline{\frak M \frak M}}_p$.

\begin{definition}
\label{local_data} 
{\rm Let $(X,\mu) \in {\overline{\frak M \frak M}}$, and let $\Omega \subset X$ be open and
non-empty.
The set $$\big(\Omega,\, (\la_k)_{k=0}^\infty ,\, (\phi_k|_{\Omega})_{k=0}^\infty \big),$$
where $\phi_k$ is orthonormalized in $L^2(X,\mu)$,
is called 
the local spectral data (LSD) for $X$ on $\Omega$.}

{\rm Moreover, suppose that $\big(\Omega,\, (\la_k)_{k=1}^\infty,(\phi_k|_{\Omega})_{k=0}^\infty\big)$
and $\big({\Omega'},\, ( \la'_k)_{k=0}^\infty,( \phi'_k|_{{ \Omega'}})_{k=0}^\infty\big)$ 
are the  LSD of 
 $(X, \mu)$ and $(X', \mu')$, respectively.
We say that these data are equivalent by a map
$\Psi_{ \Omega}: \Omega \to {\Omega'}$, if $\Psi_{ \Omega}$  is a homeomorphism, and 
\ba
\la_k = \la'_k,\quad \Psi^*_{ \Omega}(\phi'_k|_{{\Omega'}})= \phi_k|_{\Omega}, \quad \hbox{for all } k\in \N.
\ea}
\end{definition}
Our aim is to prove  the uniqueness  result.

\begin{theorem} \label{main-map}
Assume that the  LSD of 
 $(X, \mu),\,(X', \mu') \in {\overline{\frak M \frak M}}$ are equivalent by a  homeomorphism
$\Psi_{ \Omega}: \Omega \to {\Omega'}$.
Then there exists a measure-preserving isometry $\Psi:X\to X'$
    such that  $\Psi|_{\Omega}=\Psi_{\Omega}$, i.e.,
\ba
d_{X'}(\Psi(x),\,\Psi(y))=d_X(x, y),\quad x, y \in X,\; \textrm{ and }\,
\mu= \Psi^*(\mu').
\ea
Moreover, $\Psi|_{X^{reg}}: X^{reg} \to (X')^{reg}$ is a $C^3_*$-Riemannian isometry.
\end{theorem}
\begin{corollary} \label{for_ball}
Let $(X, p, \mu),\, (X', p', \mu') \in \overline{\frak M  \frak M}_p$.
Let $z_\a,\, z'_\a,\, \a=0, 1, \cdots,$ be  dense in the ball $B_X(p, r),\,B_{X'}(p', r)$ of some radius $r>0$ with $z_0=p$, $z'_0=p'$, respectively, and $t_\ell,\,
 \ell=1, 2, \cdots$, be dense in $(0, \infty)$. Assume that
 \ba
 H_X(z_\a, z_\beta, t_\ell)=H_{X'}(z'_\a, z'_\beta, t_\ell), \quad \textrm{for all }\,
 \a, \beta\in \N,\; \ell\in \Z_+.
 \ea
 Then  the conclusion of Theorem \ref{main-map} is  valid and, in addition, 
$\Psi(p)=p'.$
\end{corollary}

First, we prove the corollary using Theorem \ref{main-map}.

\smallskip
\noindent {\bf Proof of Corollary \ref{for_ball}}. By Lemma \ref{PHD_to_LSD}, the conditions of the
corollary imply that $\la_j=\la'_j, \, \phi_j(z_\a)=\phi'_j(z'_\a)$ for all $j,\alpha\in \N$, where the last equality is considered modulus
an orthogonal transformation in ${\mathcal L}(\la_j)$ and,
without loss of generality, we assume this transformation to the identity. 
Let us denote $B=B_X(p,r)$ and $B'=B_{X'}(p',r)$.
By taking closure in the product topology of $\R^{\N}$ of the set
$\Phi(z_\a)=\Phi'(z'_\a),\, \a=0,1,\cdots$, where
\bequ \label{1.09.02.2012}
\Phi(x)= \big(\phi_j(x)  \big)_{j=0}^\infty, \quad
\Phi'(x')= \big(\phi'_j(x')  \big)_{j=0}^\infty,
\eequ
 we obtain
the images 
$\Phi(\overline{B})=\Phi'(\overline{B'}) \subset \R^{\N}$.
Since $\overline{B},\, \overline{B'}$ are compact,
by Lemma \ref{lem: eigenfunctions injectivity}, 
there exists a finite set
 $J \subset \N$, such that 
$\Phi_J: B \to \R^{|J|},\,  
\Phi'_J: B' \to \R^{|J|}$, where $\Phi_J(x)=(\phi_{j(1)}(x), \dots, \phi_{j(J)}(x))$ 
and similarly for $\Phi'_J$, 
are homeomorphisms onto their images. 
Then these images coincide and the map
$$
 \Psi_{\overline{B}}:= (\Phi'_J)^{-1} \circ \Phi_J:\, \overline{B} \to \overline{B'}
$$
is a homemorphism, and the LSD are equivalent by $\Psi_{\overline{B}}$.
Let us choose any open set $\Omega\subset \subset B$ such that the image $\Psi_{\overline{B}}(\Omega) \subset \subset B'$.
For example, the set $\Omega$ can be chosen as an open neighborhood of an interior point in $\Psi_{\overline{B}}^{-1}(B_{X'}(p',r/2))\subset \overline{B}$.
Then $\Psi_{\overline{B}}|_{\Omega}$ is a homeomorphism from $\Omega$ to $\Psi_{\overline{B}}(\Omega)$, and the LSD is equivalent on these sets by definition.
Note that the sets $\Omega, \Psi_{\overline{B}}(\Omega)$ are open with respect to the topology of $X,X'$, respectively. 
Hence Theorem \ref{main-map} applies to $\Omega,\Psi_{\overline{B}}(\Omega)$, and we obtain a measure-preserving isometry $\Psi:X\to X'$ such that $\Psi|_{\Omega}=\Psi_{\overline{B}}|_{\Omega}$. 
In particular, $\Psi|_{\{z_{\alpha}\}}$ is an isometry from the dense points $\{z_{\alpha}\}\subset B$ to $\{z'_{\alpha}\}\subset B'$, and thus $\Psi(p)=p'$.
\hfill \textbf{QED}

\medskip
 The proof of Theorem \ref{main-map}
 is rather long and will consist of several parts.

\subsection{Blagovestchenskii identity on ${\overline{\frak M \frak M}}$}

For  $(X, \mu_X) \in {\overline{\frak M \frak M}}_p$,  consider the following 
initial-boundary value problem (IBVP) for the wave equation
\beq
\label{IBVP1}
& &\left( \p_t^2 +\Delta_{X} \right)u(x,t)= H(x,t), \quad \hbox{in } X\times \R_+,
\\ \nonumber
& &  u|_{t=0}=u_0, \quad \partial_t u|_{t=0}=u_1,
\eeq
with $u_0 \in H^1(X),\, u_1 \in L^2(X)$, and $H \in L^2 (\Omega  \times \R_+)$ for an open set $\Omega\subset X$, where 
\ba
L^2(\Omega \times (0, T))=\Big\{H \in L^2(X \times (0, T)): H=0\hbox { outside } \Omega \times (0, T) \Big\}.
\ea 
We denote the  unique solution of $(\ref{IBVP1})$, when $u_0=u_1=0$,  by $u^H(x,t)$.
Sometimes we denote it by $u^H(t)=u^H(\,\cdotp,t)$. Recall that,
\bequ \label{1.06.02.2012}
    u^H(x, t)= \sum_{j=0}^\infty u_j^H(t)  \phi_j(x), \quad
u_j^H(t)=\bra u^H(t),\phi_j\cet_{L^2(X, \mu_X)}.
\eequ

\begin{lemma} \label{Blagovestchenskii}
Assume that we are given the local spectral data (LSD), that is,  $\Omega\subset X^{reg}$, $(\la_j)_{j=0}^\infty,$ 
$(\phi_j|_{\Omega})_{j=0}^\infty$. Let $H \in L^2 (\Omega  \times \R_+)$. 
Then 
\bequ \label{04.10.02}
u^H_j(t)= \normalized{c}^{-1}
\int_0^t \int_{\Omega} \frac{\sin({\sqrt \la_j} (t-\tau) )}{{\sqrt \la_j}}\,H(x,\tau)  \phi_j(x)\,
\normalized \rho(x)  dV_X(x) d\tau,
\eequ
where
the metric $h$ and $\normalized \rho$ on $\Omega$ are determined by  Lemma \ref{recovery}, and $\normalized c>0$ is some (unknown) constant, and $d V_X$ stands for the Riemannian volume element on the regular part $X^{reg}$.
Hence, when we are given LSD on $\Omega$ and  
%
a function
$H\in L^2(\Omega  \times \R_+)$, we can  determine the numbers
$\normalized u_j^{H}(t)=\normalized{c} u_j^{H}(t)$.

%
\end{lemma}

\begin{proof}
Recall that $d\mu_X=\rho \, dV_X$ with the density function $\rho$.
Using (\ref{IBVP1}), the self-adjointness of $\Delta_X$ w.r.t. $d\mu_X$, we obtain
\beq \nonumber
\left( \p_t^2 +\la_j \right) u^H_j(t)&= & \int_X 
 \left( \p_t^2 +\la_j \right) 
u^H (x,t)  \phi_j(x)  d \mu_{X} \\
&=& \int_X 
 \phi_j(x)\left( -\Delta_X +\la_j \right) 
u^H(x,t)   d \mu_{X} + \int_X 
H(x,t) \phi_j(x)  d \mu_{X} \nonumber \\
& =&\int_X 
H(x,t) \phi_j(x) \rho(x) d V_{X} \nonumber \\
\label{eq: integration by parts}
&=&\normalized{c}^{-1} \int_{\Omega} H(x, t)  \phi_j(x)\bar \rho(x) dV_X,
\eeq
where 
$\normalized \rho=\normalized c \rho$
and $\normalized c$ defined in (\ref{rho-normalized}). 
Since $u_j^H(\,\cdotp, 0)=\partial_t u_j^H(\,\cdotp,0)=0$, 
equation (\ref{eq: integration by parts}) implies (\ref{04.10.02}).
\end{proof}
%
%
%

\subsection{Approximate controllability}

\begin{definition}
\label{domain_of_influence}
Let $V \subset X$ be open, $V \neq \emptyset$.  The domain of influence of $V$ in $X$ at time $T$,
 is the open set
$
X(V,T)=\big\{x \in X: d(x, V) <T \big\}.
$
\end{definition}

Below, we use the representation $X=Y/O(n)$ for $Y \in {\overline{\frak F \frak M \frak M}}$, and denote by $\pi:Y\to X$ the map from $y\in Y$ to the
corresponding orbit. 
We have the following generalization
of finite speed propagation of waves.  

\begin{lemma} \label{finite_velocity}
Let $V \subset X$ be open and $T>0$. Assume
$\hbox{supp}(u_0),\,\hbox{supp}(u_1) \subset  V$ and $\hbox{supp}(H) \subset  V \times \R_+$.
Then
$
\hbox{supp}(u(\cdot, T)) \subset X(V, T).
$
\end{lemma}
\begin{proof}
Let $\tilde{u}(y,t)$ be the solution of the wave equation on $Y$,
\beq\label{IVP on Y}
& &\left( \p_t^2 +\Delta_{Y} \right)\tilde{u}(y,t)= \tilde{H}(y,t), \quad \hbox{in } Y\times \R_+,
\\ \nonumber
& &  \tilde{u}(y,t)|_{t=0}=\tilde{u}_0(y), \quad \partial_t \tilde{u}(y,t)|_{t=0}=\tilde{u}_1(y)
\eeq
where $\tilde{u}_0 = u_0\circ \pi$ and $\tilde{u}_1 = u_1\circ \pi$ are functions on $Y$ that are
supported on $\tilde{V}:=\pi^{-1}(V) \subset Y$, 
and $\tilde{H} = H\circ \pi$ is a function on $Y\times \R_+$ supported on $\tilde{V} \times \R_+$.
Since $\rho_Y, \, \tilde{h}_Y \in C^2_*(Y)$ by Corollary \ref{density_on_Y}, a slight modification of the classical result
shows that
$
\hbox{supp}(\tilde{u}(\cdot, T)) \subset Y(\tilde{V},T),
$
where $ Y(\tilde{V},\, T)$ is the domain of influence on $Y$ corresponding to $\tilde{V}\subset Y$. Observe that $X(V,T)=\pi (Y(\tilde{V},T))$.

Since $\tilde{u}_0$, $\tilde{u}_1$ and $\tilde{H}$ are all $O(n)$-invariant,
we see that the solution $ \tilde{u}(y,t)$ is also
$O(n)$-invariant, see Appendix \ref{OTA}. Thus, $u(x,t)=\tilde{u}(y,t)$ for any $y\in \pi^{-1}(x)$,
and hence $\hbox{supp}(u(T))\subset \pi (Y(\tilde{V},T))=X(V,T)$.
\end{proof}

Next, we generalize Tataru's controllability result.
\begin{theorem} \label{Tataru} Let $V\subset X$ be an open set. 
  Then the set $\mathcal{Z}_X=\{u^H(\cdot, T):\, H \in C^3_*(X^{reg} \times (0,T)),\, \hbox{supp}(H) 
  \subset V \times (0, T)\}$  
  is dense in $L^2(X(V, T))$.
\end{theorem}
\begin{proof}
Let $\tilde{V}=\pi^{-1}(V)$ be an open subset of $Y$. 
Since $\rho_Y, \, \tilde{h}_Y \in C^2_*(Y)$ by Corollary \ref{density_on_Y}, the classical Tataru's unique continuation
result, see \cite{Tataru}, remains valid for the wave equation (\ref{IVP on Y}) on $Y$.
Let $\tilde{u}^{\tilde{H}}(y,t)$ be the solution of  (\ref{IVP on Y}) with $\tilde{u}_0=0$ and $\tilde{u}_1=0$. 
 Using the standard duality
arguments, e.g. \cite[Sec.\ 2.5]{KKL}, this implies that
 the set
$\mathcal{Z}_Y=\{\tilde{u}^{\tilde{H}}(\cdot, T):\, \tilde{H} \in C^3_*(Y \times (0,T)),\,
\hbox{supp}(\tilde{H}) 
  \subset \tilde{V} \times (0, T)\}$ is dense in $L^2(Y(\tilde{V}, T))$.
Hence by the proof of Lemma \ref{finite_velocity}, the $O(n)$-invariant part $\mathcal{Z}_Y^O=\{\tilde{u}^{\tilde{H}}(\cdot, T):\, \tilde{H} \in C^3_{*,O}(Y \times (0,T)),\,
\hbox{supp}(\tilde{H}) 
  \subset \tilde{V} \times (0, T)\}$ is dense in $L^2_O(Y(\tilde{V}, T))$, where $L^2_O(Y(\tilde{V}, T))$ denotes the $O(n)$-invariant subspace of $L^2(Y(\tilde{V}, T))$, see Appendix \ref{OTA}.
Together with \eqref{000.12}, this implies that $Z_X=\pi_*(Z_Y)$  is dense in $L^2\big(X(V, T) \big)=\pi_*\big(L^2_O(Y(\tilde{V}, T)) \big)$.
\end{proof}
%
%
%
%

\begin{lemma} \label{basis}
  For any open set $V \subset \Omega \subset X^{reg}$ and $ s>0$, there exist 
  $\{F_k\}_{ k=1}^\infty\subset   L^2(V \times (0, s))$
such that  for $ a_k(x)= u^{F_k}(x, s)$,
 \medskip
 
\hspace{0.5cm} $(*)$\, 
$\{a_k(x)\}_{k=1}^\infty$
 form an  
 orthonormal basis in $L^2(X(V, s), \normalized{c}^{2} \mu_X)$.
 \medskip
 
\noindent Furthermore,   using the LSD, that is, $\Omega,$ $\la_j, \phi_j|_\Omega,\, j=0, 1, \dots$, it is possible to construct sources
 $F_k$ such that
 $(*)$ holds.
\end{lemma}

\begin{proof}
Let $H_k \in C^3_*(V \times (0, s)),$ $ k=1, 2, \dots,$ be dense in $L^2(V \times (0, s))$. Then by Theorem \ref{Tataru}, functions 
$ b_k(x):=u^{H_k}(x, s)$ are dense in $L^2(X(V, s), \mu_X)$.
By Lemma \ref{Blagovestchenskii}, the LSD and functions 
$H_k, H_l$ determine $\normalized u_j^{H_k}(s)$ and  $\normalized u_j^{H_l}(s)$, and thus
the non-negative quadratic form 
by 
\bequ \label{04.11.06}
Q[H_k,\, H_l]:=\sum_{j=0}^\infty \normalized u_j^{H_k}(s)\, \normalized u_j^{H_l}(s)= 
\normalized{c}^2 \bra u^{H_k}(s), \,u^{H_l}(s)\cet _{L^2(X, \mu_X)}.
\eequ
Applying the Gram-Smidt orthonormalization algorithm to the sequence $( H_k)_{k=1}^\infty$ with respect to the
quadratic form $Q[\,\cdotp,\,\cdotp]$, we obtain a sequence  $(F_k)_{k=1}^\infty$ of functions $F_k\in L^2(V \times (0, s))$
such that $Q[F_k, F_l]=\delta_{k l}$.
Hence $a_k=u^{F_k}(\,\cdotp,s)$ are  orthonormal and span a dense set in $L^2(X(V, s), \normalized{c}^{2} \mu_X)$
\end{proof}


By (*), the functions 
$$\Phi_j(x):=\normalized{c}^{-1}\phi_j(x)$$ 
are an orthonormal basis of  $L^2(X, \normalized{c}^2\mu_X)$.
For a measurable set $A\subset X$, we denote 
$$
(P_A v) (x):=   \chi_{A}(x) v(x),
$$
where $\chi_A$ is the indicator function of $A$.
By Lemma \ref{basis}, we have 
$$
(  P_{{X(V, s)}} v) (x)=   \chi_{X(V, s)}(x) v(x) 
 = \sum_{k=1}^\infty \bra v, a_k\cet_{L^2(X,\normalized{c}^2d\mu_X)}\,  a_k(x).
$$
Then, 
\ba
& &\bra P_{{X(V, s)}}\Phi_i,\Phi_j\cet_{L^2(X,\normalized{c}^2d\mu_X)}=\bra  \sum_{k=1}^\infty \bra \Phi_i, a_k\cet_{L^2(X,\normalized{c}^2d\mu_X)} \,a_k ,\Phi_j
\cet_{L^2(X,\normalized{c}^2d\mu_X)}
\\
&& =\sum_{k=1}^\infty \bra \Phi_i, a_k\cet_{L^2(X,\normalized{c}^2d\mu_X)} \,\bra a_k ,\Phi_j\cet_{L^2(X,\normalized{c}^2d\mu_X)}
= \sum_{k=1}^\infty \normalized u_i^{F_k}(s) \normalized u_j^{F_k}(s). 
\ea
%
%
Thus, as we can compute $\normalized u_i^{F_k}(s)$ using
 Lemma \ref{Blagovestchenskii}, we see that if
we are given LSD on $\Omega$,
we can compute the operator $\M_{X(V,s)}$ for any $V\subset \Omega$ corresponding to  the orthogonal projector $P_{X(V,s)}$, defined by
\beq\label{M-map}
   \M_{X(V,s)}: \ell^2 \rightarrow \ell^2,\quad \M_{X(V,s)} = 
   {\mathcal F}_{\normalized{c}}\circ P_{X(V, s)}\circ
     {\mathcal F}^{-1}_{\normalized{c}},
\eeq
where ${\mathcal F}_{\normalized{c}}:L^2(X,\normalized{c}^2\mu_X)\to \ell^2$ is the isometry 
$$
{\mathcal F}:v \mapsto \Big(\int_X v(x) \Phi_j(x) \,\normalized{c}^2d\mu_X \Big)_{j=0}^\infty.
$$
We note that the operator   $\M_{X(V,s)}: \ell^2 \rightarrow \ell^2$ is independent of $\normalized{c}$,
and it can be considered as an infinite matrix whose elements are
\bequ \label{M-element}
\bra P_{{X(V, s)}}\Phi_i,\Phi_j\cet_{L^2(X,\normalized{c}^2d\mu_X)}=\bra P_{{X(V, s)}}\phi_i,\phi_j\cet_{L^2(X,d\mu_X)}=\int_{X(V,s)}\phi_i \phi_j \, d\mu_X.
\eequ

Let  $V_l \subset \Omega$ be open sets and $s_l \geq 0$, $l=1, 2,\dots, 2L$. Denote
\bequ
\label{I1}
   I:=  \bigcap_{l=1}^L \Big( X(V_{2l-1}, s_{2l-1}) \setminus X(V_{2l},
         s_{2l}) \Big). 
 \eequ
Then for given $V_l$  and $s_l$, LSD determine the operator
\bequ\label{MI-matrix}
\M_I={\mathcal F}_{\normalized{c}} \circ P_{I} \circ
     {\mathcal F}^{-1}_{\normalized{c}} ,\quad P_I=  \prod _{l=1}^L P_{X(V_{2l-1}, s_{2l-1})}\big(I-P_{X(V_{2l},
         s_{2l})} \big).
\eequ
The operator $\M_I$ is non-zero if and only if the set $I$  has non-zero measure.



%
%
%
%
%
%
%
%

\begin{lemma} \label{intersection} 
Let $\big(\Omega,(\lambda_j)_{j=0}^\infty,(\phi_j|_\Omega)_{j=0}^\infty \big)$ and 
$\big(\Omega',(\lambda'_j)_{j=0}^\infty,(\phi'_j|_{\Omega'})_{j=0}^\infty \big)$, where $\Omega\subset X^{reg}$, $\Omega'\subset (X')^{reg}$,
be the LSD for $(X, \mu)$ and $(X', \mu')$,
and assume that these data are equivalent by $\Psi_{\Omega}:\Omega\to \Omega'$.
Let $V_l \subset \Omega$, $V_l':=\Psi_{\Omega}(V_l)\subset \Omega'$ be open sets, and $s_l \geq 0$, $l=1,2, \dots, 2L$.
Let $I\subset X$ be given by (\ref{I1}) and define
\bequ
\label{I'}
   I':=  \bigcap_{l=1}^L \left( X'(V_{2l-1}', s_{2l-1}) \setminus X'(V_{2l}',
       s_{2l}) \right)\subset X'.
\eequ
Then  $\mu(I)\neq 0$ if and only if  $\mu'(I')\neq 0$, and
\bequ \label{04.18.09}
 \int_{I} \phi_i(x) \phi_j(x) d\mu=
\int_{I'} \phi'_i(x') \phi'_j(x') d\mu',\quad \forall \, i,j\in \N.
\eequ
\end{lemma}

\begin{proof}
As the LSD on $X$  and $X'$  coincide, the operators $\M_{X(V_{l}, s_{l})}$ and  $\M_{X'(V'_{l}, s_{l})}$ given by \eqref{M-map} coincide. Hence $\M_I$ is non-zero if and only if  $\M_{I'}$ is non-zero, and the integrals \eqref{04.18.09} coincide by definition, see \eqref{M-element}.
\end{proof}

\subsection{Cut locus}

All locally compact length spaces are geodesic
spaces, that is, any pair of points can be joined by a distance-minimizing path.
By \cite{GLP},  the Gromov-Hausdorff limit  of a sequence of 
locally compact length spaces 
is again a  locally compact length space, and 
hence,  $X\in{\overline{\frak M \frak M}}$ is a geodesic space.


\begin{definition} \label{geodesic}{\rm 
We say that a path $\g:[0,\ell]\to X$ is a geodesic on $X$ if it is a locally  distance-minimizing path and is parametrized with the arc length. The curve $\gamma\cap X^{reg}$ 
is a geodesic in the Riemannian sense on $X^{reg}$.

For  $x \in X^{reg}$ and a unit vector $\xi \in S_x(X)$,  
we denote by 
$\g_{x, \xi}(t)$, $t\ge 0$ a geodesic of $X$ starting from $x$ in the direction
$\xi$.
 We denote by $I_{x, \xi}$ the largest interval on
which the geodesic  $\g_{x, \xi}(t)$ can be defined.}
\end{definition}

\begin{definition} \label{cut_locus}
{\rm For  $x\in X^{reg}$ and $\xi\in S_x(X)$,  let
      $i(x,\xi) $
be the supremum of those $t\in I_{x, \xi}$ that
  $d(x,\g_{x,\xi}(t))=t$. 
We call $i(x, \xi)$ the cut locus distance function.
The injectivity radius $i(x)$ at $x$ is defined as
$
      i(x) = \inf\limits_{\xi\in S_x(X)} i(x,\xi).
$
}
\end{definition}



\begin{lemma}\label{lem: distance to Xsing}
For any $x\in X^{reg}$, we have $
i(x)\leq d_X(x,X^{sing}).$
\end{lemma}

\begin{proof} 
Let $z\in X^{sing}$ be a est point in $X^{sing}$ from $x$.
If the claim is not valid, then we have 
$
t_0=d_X(x, z) <i(x)$.
Then there exists a shortest $\g_{x, \xi}(t)$ joining $x$ and $z$ with $\g_{x, \xi}(t_0)=x'$, and the geodesic 
$\g_{x, \xi}([0,t_1])$ is distance-minimizing at least till $t_1=i(x)$. However, due to \cite[Thm. 1.1(A)]{Pet},
this implies that $\g_{x, \xi}([0,t_1)) \subset X^{reg}$, contradicting 
$z=\g_{x, \xi}(t_0) \in X^{sing}$.
\end{proof}

\begin{figure}[h]
  \begin{center}    \includegraphics[width=0.58\linewidth]{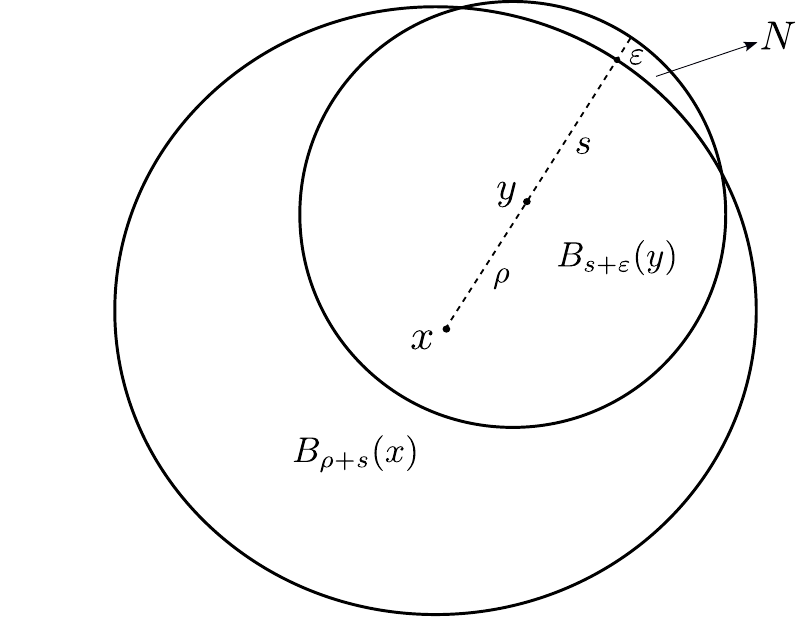}
    \caption{The set $N=N(x,\xi;\rho,s,\varepsilon)$ where $y=\gamma_{x,\xi}(\rho)$.}
    \label{fig_setN}
  \end{center}
\end{figure}

Let $y= \g_{x, \xi}(\rho)$  with $0<\rho < i(x)$ and $s, \e >0$.
We denote
\bequ \label{domains}
N(x, \xi;\, \rho,s, \e)
=B_{s+\e}(y)  \setminus B_{\rho+s}(x),
\eequ
where $B_r(x)$ denotes the open ball of radius $r$ centered at $x$, see Figure \ref{fig_setN}.


\begin{lemma} \label{11}
Let $x\in X^{reg}$ and $0<\rho < i(x)$.  
\begin{enumerate}
 \item[(a)]  If  $\rho+s <\, i(x)$, then
$
      \bigcap_{\e>0} {N(x, \xi; \rho,s, \e)} =\{ \g_{x,\xi}(\rho+s)\}.
$
\item[(b)]  If $\rho+s >\, i(x)$, then there exist
     $\xi\in S_x(X)$ and $\e >0$ such that 
$
      N(x,\xi; \rho,s, \e) = \emptyset.
$
\item[(c)]    The injectivity radius satisfies
$$ \hspace{-1cm}i(x)=\inf_{s>0} \Big\{\rho+s \,:\hbox{ there exist }\xi\in S_xX,\,  \e >0
 \hbox{ such that } \hbox{Vol}\,\big(N(x,\xi; \rho,s, \e)\big)= 0 \Big\}.$$
\end{enumerate}
\end{lemma}

\begin{proof}
$(a)$ For any $\xi\in S_x(X)$ and sufficiently small $\e>0$,
$\g_{x, \xi}([0, \rho+s+\e])\subset X$ is  distance minimizing. This implies 
that 
$
   \g_{x, \xi}([\rho+s, \rho+s+\e)) \subset N(x,\xi;\rho,s, \e).
$
Thus, $\g_{x, \xi}(\rho+s) \in  \bigcap_{\e>0} N(x, \xi; \rho,s, \e)$.

For any $q\in \bigcap_{\e>0} N(x, \xi; \rho,s, \e)$, it satisfies that $d(y,q)<s+\e$ for any $\e>0$ and $d(x,q)\geq \rho +s$. The former yields $d(y,q)\leq s$, which implies that
$$\rho+s\leq d(x,q)\leq d(x,y)+d(y,q)\leq \rho+s.$$
Here we have used the condition that $\rho<i(x)$ so that $d(x,y)=\rho$. The inequality above gives $d(x,q)=d(x,y)+d(y,q)=\rho+s$, 
which means that $q$ lies on the extension of the geodesic $\gamma_{x,\xi}$. Since $\rho+s<i(x)$, we have $q=\gamma_{x,\xi}(\rho+s)$.

\smallskip
(b) Assume the claim is not true: $N(x,\xi,\rho,s,\e)\neq \emptyset$ for any $\xi\in S_x(X)$ and any $\e>0$.
Let $\xi_0\in S_x(X)$ be a unit vector satisfying $\rho+s>i(x,\xi_0)$. We pick a sequence $\e_i\to 0$ such that $N(x,\xi_0,\rho,s,\e_i)\neq \emptyset$ and pick a point $q_i\in N(x,\xi_0,\rho,s,\e_i)$.
Then it satisfies that $d(y,q_i)<s+\e_i$ and $d(x,q_i)\geq \rho +s$. Taking a converging subsequence of $q_i\to q$, we have 
$d(y,q)\leq s$ and $d(x,q)\geq \rho +s$. Repeating the proof of (a) gives $d(x,q)=d(x,y)+d(y,q)=\rho+s$, which means that $q$ lies on the extension of the geodesic $\gamma_{x,\xi_0}(t)$, and the geodesic is minimizing until $t=\rho+s$, contradiction to $\rho+s>i(x,\xi_0)$.

\smallskip
(c) Observe that $N(x,\xi; \rho,s, \e)$ contains an open ball of radius $\e/2$ centered at $\gamma_{x,\xi}(\rho+s+\e/2)$ for sufficiently small $\e>0$ if $\rho+s<i(x)$. Then the claim (c) is a consequence of claims (a) and (b).
\end{proof}

Next, we combine Lemma \ref{intersection} and Lemma \ref{11}.

\begin{lemma} \label{geodesics}
Let $(X, \mu),\, (X', \mu') \in {\overline{\frak M \frak M}}$. Let $\Omega \subset X^{reg},\, \Omega' \subset (X')^{reg}$, and 
$\big(\Omega,\, (\la_k)_{k=0}^\infty,$ $(\phi_k|_{\Omega})_{k=0}^\infty\big)$
and $\big({\Omega'},\, ( \la'_k)_{k=0}^\infty,( \phi'_k|_{{ \Omega'}})_{k=0}^\infty \big)$ be the  LSD of 
 $X, X'$.
Assume that the LSD are equivalent by a homeomorphism
$\Psi_{ \Omega}: \Omega \to {\Omega'}$. 
Then $\Psi_{ \Omega}: \Omega \to {\Omega'}$ is a $C^{3}_*$-diffeomorphism.

Moreover, let  
$x\in \Omega$, $x':=\Psi_{\Omega}(x)$ and
$\exp_{x}:T_xX\to X$, $\exp_{x'}:T_{x'}X'\to X'$ be the exponential maps.
Then $i(x)=i(x')$, and the map $E:B(x,i(x))\to B'(x',i(x'))$ defined by
\bequ
E(z)=\exp_{x'}\Big (d \Psi_{ \Omega}|_x(  \exp_x^{-1}(z))\Big)
\eequ
satisfies $h=E^*(h'),\, \normalized\rho=E^*(\normalized\rho')$, and 
$\phi_j=\phi_j'\circ E$ on $B(x,i(x))$ for all $j\in \N$.
\end{lemma}

\begin{proof}
By Lemma \ref{lem: eigenfunctions injectivity}, 
  any $z\in \Omega$ has a neighborhood $V\subset\Omega $ 
such that there is an index ${\bf j}=(j_1, \dots, j_d)$ with $d=\hbox{dim}(X)$, for which $\Phi_{\bf j}:V\to \R^d$, $\Phi_{\bf j}(x)=(\phi_j(x))_{j\in {\bf j}},$
 defines $C^3_*$-smooth coordinates. Let $z'=\Psi_\Omega(z)$ and
 $V'=\Psi_\Omega(V)$. Then,
 as $\phi_j'(\Psi_{\Omega}(x))=\phi_j(x)$ for all $j$ by assumption,
 it follows from Lemma \ref{lem: recognize coordinates} and Remark \ref{regular2}
 that the 
 map $\Phi'_{\bf j}:V'\to \R^d$ 
defines $C^{3}_*$-smooth coordinates for the same index  ${\bf j}$. Moreover, $\Psi_\Omega=(\Phi'_{\bf j})^{-1}\circ \Phi_{\bf j}$ in $V$.
 Thus, $\Psi_\Omega:\Omega\to \Omega'$ is a $C^{3}_*$-smooth
 diffeomorphism. 

Let $\rho < \min(i(x),i'(x'))$ be small so that $B_{\rho}(x)\subset \Omega$ and $B_{\rho}(x')\subset \Omega' $,  and  let $\xi\in S_{x}(X),\,
\xi' :=d\Psi_{\Omega}|_x(\xi) \in S_{x'}(X')$.
By Lemma \ref{intersection}, for any $s,\e>0$,
$$
\mu(N(x,\xi; \rho,s, \e))=0\quad\hbox{if and only if}\quad \mu'(N'(x',\xi'; \rho,s, \e))=0.
$$
Therefore, due to Lemma \ref{11}, $i(x)=i(x').$ 

Next, we again take  $\xi\in S_{x}(X),\, 
\xi' :=d\Psi_{\Omega}|_x(\xi)$  and $0<s <i(x)-\rho$.
Note that $B(x,i(x))\subset X^{reg}$ and $B'(x',i(x'))\subset (X')^{reg}$ due to Lemma \ref{lem: distance to Xsing}, and thus the exponential maps are defined in the Riemannian sense.
Consider $x_0=\g_{x, \xi}(\rho+s),\, x'_0=\g'_{x', \xi'}(\rho+s)=E(x_0)$.
Then by Lemma \ref{intersection}, we have
$$
\phi_j(x_0)=\lim_{\e\to 0}
\frac{1}{\mu(N)}\int_{N(x, \xi; \rho, s, \e)} \phi_j\phi_0\, d\mu 
=\lim_{\e\to 0}
\frac{1}{\mu'(N')}\int_{N'(x', \xi'; \rho, s, \e)}\phi_j'\phi_0'\,d\mu'
=\phi_j'(x'_0),
$$
where we have used $\phi_0=\phi'_0=1$ since $\phi_j,\phi'_j$ are normalized w.r.t. $\mu,\mu'$.
At last, using Lemma \ref{recovery}, we obtain that $h= E^*(h'),\, 
\normalized \rho=E^*(\normalized \rho')$ on $B(x,i(x))$.
\end{proof}

We are now ready to prove Theorem \ref{main-map}.

\medskip
\noindent {\bf Proof of Theorem \ref{main-map}.}\,\,
Let ${\mathcal O}$ be the set of all pairs, $(\Omega_1,\Omega'_1)$,
of connected open subsets $\Omega_1 \subset X^{reg},\,\Omega'_1 \subset (X')^{reg}$ 
containing $\Omega, \Omega'$ 
such that 
$\Psi_{\Omega}$ extends to a diffeomorphism $\Psi_1:\Omega_1\to \Omega'_1$ 
and satisfies  $\Psi_1^*(h'|_{\Omega_1'})=h|_{\Omega_1}$,
$\Psi_1^*(\normalized \rho'|_{\Omega_1'})=\normalized \rho |_{\Omega_1}$
and
 $\Psi_1^*(\phi'_j|_{\Omega_1'})=\phi_j|_{\Omega_1}$ for every $j\in \N$.
 Note that the map $\Psi_1$ is an isometry by the first condition, and by the last condition, the map 
 $\Psi_1$ has to coincide with $(\Phi')^{-1}\circ \Phi|_{\Omega_1}$, where 
$\Phi:X\to \R^{\N}$ and  $\Phi':X'\to \R^{\N}$ are the maps defined in
(\ref{1.09.02.2012}).

 Let us consider  ${\mathcal O}$ as a partially ordered set, where the order is given
 by  the inclusion $\subset$.
Suppose that $I$ is an ordered index set and $(\Omega_i, \Omega'_i)\in {\mathcal O}$, $i\in I$, are
sets so that $\Omega_i\subset \Omega_l,\, \Omega'_i\subset \Omega'_l$ for $i<l$. Let $\Psi_i:\Omega_i\to \Omega'_i$ 
be a corresponding diffeomorphism. Then it follows
 from Lemma \ref{lem: eigenfunctions injectivity} that $\Psi_l|_{\Omega_i}=\Psi_i$ for $i<l$.
Hence we see that
$\left(\bigcup_{i\in I}\Omega_i, \, \bigcup_{i\in I}\Omega'_i \right) \in {\mathcal O}$. 
Then by Zorn Lemma, there  exists a maximal element  
$(\Omega_{\max},\,\Omega'_{\max})$ in ${\mathcal O}$. 
Let  $\Psi_{\max}:\Omega_{\max}\to \Omega'_{max}$ be a diffeomorphism (an isometry) corresponding 
to $(\Omega_{\max},\,\Omega'_{\max})\in {\mathcal O}$.

\begin{lemma} \label{ass:max} The maximal element satisfies
$(\Omega_{\max},\Omega'_{\max})=\big(X^{reg},(X')^{reg} \big) $.
\end{lemma}

\begin{proof}
Assume that the claim is not valid. Then,
without loss of generality,  we can assume that $X^{reg}\setminus \Omega_{\max}\not=\emptyset$. 
As $X^{reg}$ is connected, we see that there exists  a point $z_0$ in
$\partial \Omega_{\max}\cap X^{reg}$. Let  us choose a  sequence of 
points $y_k \in \Omega_{max}$
such that  $ y_k \to z_0$ as $k\to \infty$ and let 
$y_k'=\Psi_{max}(y_k)\in \Omega'_{max}$. As $X'$ is compact, by choosing a subsequence,
we can assume $y'_k$ converges to some $z'_0\in X'$. 
Note that at this moment we do not know whether $z_0'\in (X')^{reg}$.

Let  $W=B(z_0,t_0)\subset X^{reg}$ be a neighborhood
of $z_0$ such that
$\overline W\subset X^{reg}$.
 By Lemma \ref{injectivity radius}, 
there exists a positive number $i_0>0$ such that
$\hbox{inj}(p) \geq i_0$ for all $p \in \overline {W}$.
Let 
$r_0:=\min\{i_0,t_0\}/2.$

Let us choose (and fix) $k$ such that $z_0\in  B(y_k,r_0)\subset X^{reg}$.
Note that $i(y_k)\geq i_0$.
Since $y_k \in \Omega_{max} \subset X^{reg},\, y'_k \in \Omega'_{max} \subset (X')^{reg}$, 
there exists $0<\rho <\min\{i(y_k), i(y'_k)\}$ such that
  $B:=B(y_k, \rho) \subset \Omega_{max}$ and $B':= B'(y'_k, \rho) \subset \Omega'_{max}$.
Then the LSD  
 $\left(B, (\la_j)_{j=0}^\infty,\, (\phi_j|_B)_{j=0}^\infty \right)$ 
and  $\left(B', (\la'_j)_{j=0}^\infty,\, (\phi'_j|_{B'})_{j=0}^\infty \right)$ 
 are  equivalent by the isometry
$\Psi_{\max}|_{B}$. Using Lemma \ref{geodesics}, 
we can determine  $i(y_k)=i(y'_k)$.
Moreover, as $i(y_k)\geq i_0$,
 we can construct a diffeomorphism $E:B(y_k,r_0)\to B'(y'_k,r_0)$
such that
\bequ \label{04.18.10AA}
\phi_j(x)= \phi'_j(E(x)), \quad \hbox{for\,all}\,\, j \in \N \,\,\hbox{and}\,\,x\in B(y_k,r_0).
\eequ
In particular, $B'(y'_k,r_0)\in (X')^{reg}$.
This implies that
 $E$ coincides with the map $(\Phi')^{-1}\circ \Phi|_{B(y_k,r_0)}$, where $\Phi, \Phi'$
 are defined in (\ref{1.09.02.2012}). 
Note that $\Psi_{\max}$ is also the restriction of 
$(\Phi')^{-1}\circ \Phi$ onto $\Omega_{max}$. 
Since $E: B(y_k,r_0)\to  B'(y'_k,r_0)$ and 
$\Psi_{\max}:\Omega_{\max}\to
\Omega_{\max}'$ 
are surjections, we see that 
$\Psi_{ext}:=(\Phi')^{-1}\circ \Phi|_{\Omega_{\max}\cup B(y_k,r_0)}$
is a bijection from $ \Omega_{ext}=\Omega_{\max}\cup B(y_k,r_0)$ to $\Omega_{ext}'=\Omega'_{\max}\cup B'(y'_k,r_0)$.
Then it follows that $\Psi_{ext}$ is a diffeomorphism.
In particular,
we have $\Psi_{ext}^*(\phi_k')=\phi_k$ on $\Omega_{ext}$, and using Corollary
\ref{recovery}, we see that $\Psi_{ext}^*(h')=h$ and $\Psi_{ext}^*(\normalized \rho')=\normalized \rho$
on $\Omega_{ext}$. Thus, $(\Omega_{ext},\,\Omega'_{ext}) \in {\mathcal O}$.

At last, observe that $B(y_k,r_0)\setminus \overline{\Omega}_{max}\neq \emptyset$ since $z_0\in B(y_k,r_0)\cap \partial \Omega_{max}\neq \emptyset$.
Hence
$\Omega_{ext}\not\subset \overline{\Omega}_{\max},\, \Omega'_{ext}\not\subset \overline{\Omega'}_{\max}$, 
which is in contradiction
with the the assumption that $(\Omega_{max}, \Omega'_{max})$ 
is the maximal element of $ {\mathcal O}$.
Thus, $(\Omega_{max}, \Omega'_{max})=\big(X^{reg},(X')^{reg}\big) $.
\end{proof}

\noindent {\bf End of Proof of Theorem \ref{main-map}.} 
Let $\Omega\subset X$ and $\Omega'\subset X'$ be the open sets in Theorem \ref{main-map}, possibly containing singular points. Since any homeomorphism between $\Omega$ and $\Omega'$ has to preserve the regular part, $\Psi_{\Omega}|_{\Omega\cap X^{reg}}$ is a homeomorphism from $\Omega\cap X^{reg}$ to $\Omega'\cap (X')^{reg}$, and hence the LSD on $\Omega\cap X^{reg},\, \Omega'\cap (X')^{reg}$ are equivalent.
Thus, applying the arguments above yields that $X^{reg}$ and $(X')^{reg}$ are isometric via an isometry $\Psi_{max}$.
Since ${\overline {X^{reg}}}=X,\, {\overline {(X')^{reg}}}=X'$, extending $\Psi_{max}$ by continuity to $X$, we obtain
an isometry $
\Psi: X \to X'$.
By construction, $h= \Psi^*(h')$, $\normalized \rho= \Psi^* (\normalized \rho')$ on $X^{reg}$, and $\Psi|_{\Omega}=\Psi_{\Omega}$.
Since $\mu(X)=\mu'(X')=1$, the density functions satisfy $\rho=\Psi^*(\rho')$ on $X^{reg}$.
Hence as $\mu(X^{sing})=\mu'((X')^{sing})=0$, we see that $\mu=\Psi^*(\mu')$ on $X$.
Moreover, the same considerations as in
Lemma \ref{geodesics}, with $\Omega=X^{reg}$, show that $\Psi|_{X^{reg}}$ 
is $C^3_*$-Riemannian isometry.
\hfill {\bf QED}

\section{Stability of inverse problem} \label{Stability of inverse problem}
In this section we prove the main Theorem \ref{main stability thm}.

\medskip
\noindent {\bf Proof of Theorem \ref{main stability thm}.} To prove the result it is enough to show that, for any $\e>0$, there is
$\delta >0$, such that if the conditions of the theorem are satisfied with this $\delta$, then 
(\ref{2.28.07.10}) is valid with $\e$ instead of $\omega_r(\delta)$.

Assume that, for some $\e>0$, there exist pairs  
$(X_i, p_i, \mu_i),\, (X'_i, p'_i, \mu'_i) \in \overline{\frak M \frak M}_p$,
which satisfy the conditions of the theorem with $i^{-(2s_F+n/2+4)}$ instead of $\delta$, but
\ba
d_{pmGH}(X_i,\, X'_i) \geq 6\e.
\ea
Using regular approximations of $X_i, \, X'_i$ by pointed manifolds from $\frak M \frak M_p$,
there are pairs of pointed manifolds 
$(M_i, p_i, \mu_i),\, (M'_i, p'_i, \mu'_i) \in {\frak M \frak M}_p$
which satisfy the conditions of the theorem with $1/i$ instead of $\delta$, but
\beq \label{3.28.07.10}
d_{pmGH}(M_i,\, M'_i) \geq \e.
\eeq
Indeed, we can find $M_i, M_i'$ such that
\ba
d_{pmGH}(X_i, M_i),\; d_{pmGH}(X'_i, M_i') < \e/2,
\ea
so that \eqref{3.28.07.10} is satisfied.
Moreover, using Theorem \ref{heat_continuity_pointed}, there exist $i^{-(2s_F+n/2+4)}$-nets, 
$\{x^i_\beta\} \subset B(p_{X_i}, r),\, \{y^i_\beta\} \subset B(p_{M_i}, r)$
and $\{x'^i_\beta\} \subset B(p_{X'_i}, r),\, 
\{y'^i_\beta\} \subset B(p_{M'_i}, r)$,
and $\{t_{m}\} \subset (i^{-1}, i)$ such that
\beq \label{13.09.2011.2}
& & |H_{M_i}(y^i_\beta, y^i_{\beta'},  t_m)-
H_{X_i}(x^i_\beta, x^i_{\beta'}, t_m)| <i^{-(2s_F+n/2+4)},
\\ \nonumber
& &
|H_{M'_i}(y'^i_\beta, y'^i_{\beta'}, t_m)-
H_{X'_i}(x'^i_\beta, x'^i_{\beta'}, t_m)| <i^{-(2s_F+n/2+4)}.
\eeq
For any $z^i_\a$ involved in the definition of $d_{pmGH}(X_i, \, X_i')$, take a point 
$x^i_{\beta(\a)}$ which is $i^{-(2s_F+n/2+4)}$ close to $z^i_\a$ and similarly for 
$z'^i_\a$. Rename $y^i_{\beta(\a)},\, y'^i_{\beta(\a)}$ as 
$y^i_\a,\, y'^i_\a,\, \a=1, \dots, A(1/i)$.
By the above inequality, for sufficiently large $i$, the points $y^i_\a,\, y'^i_\a$
form $2/i$-net on $M_i, M'_i$, correspondingly. Similarly, for any $t_\ell$ involved
in the definition of $d_{pmGH}(X_i, \, X_i')$ take a point $t_m$ above such that
$|t_\ell-t_m|<i^{-(2s_F+n/2+4)}$.
Compare  $H_{M_i}(y^i_\a, y^i_{\a'}, t_\ell)$ and 
$H_{M'_i}(y'^i_{\a}, y'^i_{\a'}, t_\ell), \, \a =1, \dots, A(1/i)$.
Then, using (\ref{000.3}) with $H_{X_i},\, H_{X'_i}$ instead of $H, H'$, 
and (\ref{13.09.2011.2}), we see from (\ref{13.09.2011.1}) that
\ba
|H_{M_i}(y^i_{\a}, y^i_{\a'}, t_\ell)-H_{M'_i}(y'^i_{\a}, y'^i_{\a'}, t_\ell)|
< i^{-(2s_F+n/2+4)} + 2C i^{-2} < i^{-1},
\ea
for sufficiently large $i$.

\smallskip
Going to a subsequence if necessary, we can assume that 
\bequ \label{4.28.07.10}
(M_i, p_i, \mu_i) 
\stackrel{f_i}\longrightarrow (X, p, \mu_X),\quad
(M'_i, p'_i, \mu'_i) \stackrel{f'_i}\longrightarrow (X', p', \mu_{X'}),
\eequ
in the sense of the pointed measured GH convergence with $f_i,\,f'_i$ being the
corresponding regular approximations. In particular,
\begin{equation} \label{ball-close}
B(p_i, r) \stackrel{f_i}\longrightarrow B(p, r), \quad B'(p'_i, r)\stackrel{f'_i}\longrightarrow B'(p', r).
\end{equation}
Let us show that 
\beq \label{5.28.07.10}
(X, p,\mu_X) \simeq (X', p', \mu_{X'}),
\eeq
where $\simeq$ stands for the measure-preserving isometry. 
This would imply that
\bfo
d_{pmGH}\left((M_i, p_i, \mu_i),\, (M'_i, p'_i, \mu_i')  \right) \to 0,\quad\hbox{as }i\to \infty,
\efo
contradicting (\ref{3.28.07.10}), which would prove the theorem.

\medskip
The rest of the proof is a proof of (\ref{5.28.07.10}). Denote
\begin{equation}
w^i_\a:=f_i(y^i_\a) \in X,\quad w'^i_\a:=f'_i(y'^i_\a)\in X', \quad i=1, \dots,\,\, \a=1, \dots, A(1/i).
\end{equation}
Let ${\Bbb P}$ be the set of the double sequences 
${\it q}:=\{i(k),  \a(k)  \}_{k=1}^\infty$, with
$\{i(k)\}_{k=1}^\infty$ being a subsequence of $\{i\}_{i=1}^\infty$, 
such that for some point
$w_{{\it q}} \in B(p, r),\,w'_{{\it q}} \in B'(p', r)$, we have
\bequ \label{6.28.07.10}
w^{i(k)}_{\a(k)}  \to w_{{\it q}}\quad \hbox{and}\quad  w'^{i(k)}_{\a(k)}  \to w'_{{\it q}},
\quad\hbox{as }k\to \infty.
\eequ
Here we use $B(p,r),\, B'(p',r)$ to denote the closed ball in $X,X'$, respectively.
Using the estimate (\ref{05.09.2011.1}) in the proof of Theorem \ref{heat_continuity},
it follows from (\ref{000.3}) that 
\bequ \label{7.28.07.10}
H(w_{{\it q}},\,w_{\it s},\,t)=H'(w'_{{\it q}},\, w'_{\it s},\,t),
\quad \textrm{ for all }\,{\it  q},\,  {\it s} \in {\Bbb P},\; t>0.
\eequ
\begin{lemma} \label{density_P} Let $\Bbb P$ be the set of the convergent double sequences
defined above. Then,
\beq \label{05.09.2011.5}
\{w_{{\it q}}:\, {\it q} \in {\Bbb P}  \} = B(p, r),\quad \{w'_{{\it q}}:\, {\it q} \in {\Bbb P}  \} = B'(p', r).
\eeq
\end{lemma}
\begin{proof}
Due to \eqref{ball-close} and the fact that  $\{y^i_\a \}_{\a=1}^{A(1/i)}$ 
form an 
$1/i$-net in $B(p_i, r)$, 
the points $\{w^i_\a \}_{\a=1}^{A(1/i)}$ form a
$\delta(i)$-net in $B(p, r)$, where $\delta(i) \to 0$ as $i \to \infty$. Analogously, the same property is valid for 
$\{w'^i_\a \}_{\a=1}^{A(1/i)}$.

Therefore, for any $w \in B(p, r)$ there exists a sequence $\a(i), \, i=1, \dots,$ such that 
$w^{i}_{\a(i)} \to w$ as $i\to \infty$.
Consider the corresponding points $w'^{i}_{\a(i)}$ in $X'$. Since the closed ball $B'(p', r)$ is 
compact, there is a subsequence
$\big\{i(k),\, \a(k)=\a(i(k))  \big\}_{k=1}^\infty$ such that $w'^{i(k)}_{\a(i(k))}$ converge. 
Then
${\it q} =\{i(k),\, \a(k)\} \in {\Bbb P}$, and
\bfo
w^{i(k)}_{\a(i(k))} \to w_{{\it q}}= w \quad\hbox{and}\quad w'^{i(k)}_{\a(i(k))} \to w'_{{\it q}},
\quad\hbox{as }k\to \infty.
\efo
This proves the first claim in (\ref{05.09.2011.5}). Switching the role of $X$ and $X'$ yields the second claim in  (\ref{05.09.2011.5}).
\end{proof}

\begin{lemma} \label{equivalence_on_P}
Suppose two double sequences ${\it q},\, {\hat{\it q}} \in {\Bbb P}$ satisfy
$w_{{\it q}}=w_{{\hat{\it q}}}$.  Then $w'_{{\it q}}=w'_{{\hat{\it q}}}.$
\end{lemma}
\begin{proof}
From the condition of lemma and (\ref{7.28.07.10}), we see that
\bfo
H(w_{{\it q}},\,w_{\it s},\,t)=H'(w'_{{\it q}},\, w'_{\it s},\,t),\quad
H(w_{{\it q}},\,w_{\it s},\,t)=H'(w'_{{\hat{\it q}}},\, w'_{\it s},\,t),
\efo
for any ${\it s} \in {\Bbb P},\, t >0$. Therefore, by Lemma \ref{density_P},  
\bfo
H'(w'_{{\it q}},\, w',\, t)=H'(w'_{{\hat{\it q}}},\, w',\, t),
\efo
for any $w' \in B'(p', r),\, t>0$.
Then it follows from Corollary \ref{heat_identification} that $w'_{{\it q}}=w'_{{\hat{\it q}}}$.
\end{proof}

Lemma \ref{equivalence_on_P} makes it possible to
 introduce an equivalence relation $\approx$ on ${\Bbb P}$:
 \ba
 q \approx \hat q \quad \hbox{iff} \quad w_q=w_{\hat q}\, \textrm{ and } \, w'_q=w'_{\hat q}.
 \ea
Then
the maps
\bfo
{\mathcal F}: \, {\Bbb P}/\approx\ \longrightarrow B(p, r),\quad {\mathcal F}({\it q})= 
w_{{\it q}},
\\ \nonumber
{\mathcal F}': \, {\Bbb P}/\approx\ \longrightarrow B'(p', r),\quad {\mathcal F}'({\it q})= 
w'_{{\it q}},
\efo
are bijections. Thus, the map $\Phi:= {\mathcal F}' \circ {\mathcal F}^{-1}: B(p, r) \to B'(p', r)$ is a bijection.
Moreover, due to (\ref{7.28.07.10}),
\bequ  \label{3.29.07.10}
H(w, \hat w, t)= H'(\Phi(w), \Phi(\hat w), t),\quad \textrm{for all }\, w, \hat w \in B(p, r),\; t>0.
\eequ
Observe that Corollary \ref{for_ball}  remains valid if, instead of using the dense sequences
$\{z_\a\}_{\a\in \N}$, $\{z'_\a\}_{\a\in \N}$, we use all points $w \in B(p, r)$ with the 
corresponding point $w'$ running over the whole $B'(p', r)$. Thus, 
using Corollary \ref{for_ball}, $\Phi$ can be uniquely extended to
 a measure-preserving isometry $\Phi: (X, p,\mu_X) \to (X', p',\mu_{X'})$.
Theorem \ref{main stability thm} is proven.
\hfill \textbf{QED}

\appendix

\section{Operator-theoretical approach to $\Delta_X$.}
\label{OTA}

In this section we redefine the weighted Laplacian operator $\Delta_X$ extending on Fukaya's \cite[Section 7]{Fuk_inv}. 
Recall Theorem \ref{thm:limitstr} (2) and Corollary \ref{density_on_Y} that, for $(X,\mu_X)\in \overline{{\frak M}{\frak M}}(n,\Lambda,D)$, there exists a $C^2_{\ast}$ Riemannian manifold $Y$ on which $O(n)$ acts isometrically in such a way that $X=Y/O(n)$, and an $O(n)$-invariant probability measure $\mu_Y$ on $Y$ satisfying $\pi_{\ast}(\mu_Y)=\mu_X$, where $\pi:Y\to X$ is the natural projection. 
Moreover, by Theorem \ref{thm:MZ} (with $M=Y$ and $G=O(n)$),
the $O(n)$-action on $Y$ is of class $C^{3}_*$,
\bequ 
F_O: O(n) \times Y \rightarrow Y, \quad 
F_O({\it o}, y)= {\it o}(y), \quad
F_O \in C^3_*(O(n) \times Y).
\label{1.29.12.2011}
\eequ
Using the operator (\ref{25.12.2011.2}), we decompose $L^2(Y,\mu_Y)$ into
$$
L^2(Y, \mu_Y)=L^2_O(Y) \oplus L^2_{\perp}(Y), 
$$
where $L^2_O(Y), \, L^2_{\perp}(Y)$ are the invariant subspaces of 
$\Delta_Y$ so that 
$$
\Delta_Y= \Delta_O \oplus \Delta_\perp.
$$
We denote by $\{\la_j^O,\,
\phi_j^O\}$ and $\{\la_j^\perp,\,\phi_j^\perp \}$ the eigenpairs of $\Delta_O$ and $\Delta_\perp$.

Due to the fact that $\pi_{\ast}(\mu_Y)=\mu_X$ by (\ref{000.15}), the map
$$
\pi_*: L^2_O(Y, \mu_Y) \to L^2(X, \mu_X)
$$
is an isometry. Thus, we can define a self-adjoint operator $A$ in $L^2(X,\mu_X)$ by
\bequ \label{25.12.2011.5}
A u= \pi_* \circ \Delta_O \circ \pi^* u, \quad
 {\mathcal D}(A)= 
\pi_*\left({\mathcal D}(\Delta_O)\right)
=\pi_* \left(W^{2, 2}_O(Y)  \right).
\eequ
On the other hand, modifying \cite[Section 7]{Fuk_inv}, we define
the Dirichlet form
$$
 a_O[u^*]= \int_Y |d u^*(y)|_{h_Y}^2 \,d\mu_Y=
\int_{ \pi^{-1}(X^{reg})} |d u^*(y)|^2 \,d\mu_Y, 
\quad {\mathcal D}(a_O)= C_O^{0,1}(Y),
$$
where we use the fact that 
$\mu_Y(\pi^{-1}(X^{reg}))=1$.
Using Kato's theory of quadratic forms \cite{Kato}, $a_O$ is closable with
${\mathcal D}(\overline{a_O})= W_O^{1, 2}(Y)$ and the associated self-adjoint
operator is $\Delta_O$. Observe that
\bequ
\label{3.19.01.2012}
\pi_*:\,C^{0,1}_O(Y) \to C^{0,1}(X)
\eequ
is an isometry, since
\bequ
\label{7.22.01.2012}
d_X(x, x')=d_Y(\pi^{-1}(x),\, \pi^{-1}(x')), \quad d_X(x, x')=d_Y(y, y'),
\eequ
for some $y \in \pi^{-1}(x),\,\, y' \in \pi^{-1}(x')$.
For $u \in C^{0, 1}(X)$, let us define
\beq \label{17.09.2014.S}
a_X[u]=\int_{X^{reg}} |d u|^2_{h_X}\, d\mu_X=
\int_{\pi^{-1}(X^{reg})} |d \left(\pi^* u \right)|^2_{h_Y}\, d\mu_Y
=
a_O[\pi^* u],
\eeq
Then $a_X$ is closable with 
$
{\mathcal D}(\overline {a_X})= \pi_* (W^{1,2}_O(Y)).
$
This defines a self-adjoint operator, $A'$ in $L^2(X,\mu_X)$. Using the
distribution duality, we can see that in local coordinates of $X^{reg}$, $A'$ is given by 
$$
A' u(x)= -\frac{1}{{\sqrt h_X} \rho_X} \p_i \left({\sqrt h_X} h_X^{i j} \rho_X \p_j u(x) \right).
$$
Thus, ${\mathcal D}(A') \subset W^{2, 2}(X^{reg}, \mu_X)$ and we denote
 $A'$ by $\Delta_X$.
Then it follows from (\ref{17.09.2014.S}) that
$$
 \Delta_X =\pi_* \circ \Delta_Y \circ \pi^* u=A, 
\quad {\mathcal D}(\Delta_X)= 
\pi_*\left(W^{2, 2}_O(Y)\right) \subset W^{2, 2}(X^{reg}).
$$
In particular,
$$\textrm{spec}(\Delta_X)=\textrm{spec}(\Delta_O)\subset \textrm{spec}(\Delta_Y).$$
Therefore, the eigenpairs $\{\la_j, \phi_j\}$ of $\Delta_X$ satisfy
\bequ \label{4.19.01.2012}
\la_j=\la_j^{O}, \quad
\phi_j=\pi_*(\phi_j^O).
\eequ



\begin{remark}{\rm 
If a sequence of manifolds $M_i$ collapse to a point $p$,  $O(n)$ acts transitively on
 $Y$. Thus, $L^2_O(Y)$ consists only of constant functions. Therefore,
$\hbox{spec}(\Delta_O)=\{0\}$. Similarly, $\Delta_X$ is the operator of multiplication
by $0$ in the space of $L^2$-functions on $\{p\}$, i.e. constants. The
results of \cite{Fuk_inv} remain valid for this case.}
\end{remark}


\begin{thebibliography}{99}


%












\bibitem{Alexakis} S. Alexakis, A. Feizmohammadi, L. Oksanen,
\emph{Lorentzian Calderón problem under curvature bounds}.
Invent. Math. {\bf 229} (2022), 87--138.

\bibitem{A90} M. Anderson, \emph{Convergence and rigidity of manifolds under Ricci curvature bounds}, Invent. Math. {\bf 102} (1990), 429--445.

\bibitem {AKKLT} 
M. Anderson, A. Katsuda, Y. Kurylev, M. Lassas, M. Taylor, \emph{Boundary regularity for the Ricci equation, geometric convergence, and Gelfand's inverse boundary problem}, Invent. Math. {\bf 158} (2004), 261--321. 


\bibitem{AG}
V. Arnaiz, C. Guillarmou, \emph{Stability estimates in inverse problems for the Schr\"odinger and wave equations with trapping}, Rev. Mat. Iberoam. \textbf{39} (2023), 495--538.


\bibitem {AP}
K. Astala, L. P\"aiv\"arinta, \emph{Calder\'on's inverse conductivity problem in the plane}, Annals
of Math. {\bf 163} (2006), 265--299.

\bibitem{Bel} M. Belishev,  \emph{An approach to multidimensional inverse problems for the wave equation}, Dokl. Akad. Nauk SSSR {\bf 297} (1987), 524--527.

\bibitem{BK} M.~Belishev,~Y.~Kurylev,~\emph{To the reconstruction of a Riemannian manifold via its spectral data (BC-method)},~Comm. PDE.~\textbf{17}~(1992),~767--804.

\bibitem{BD}
M. Bellassoued, D. Dos Santos Ferreira, \emph{Stability estimates for the anisotropic wave equation from the Dirichlet-to-Neumann map}, Inverse Probl. Imag.~\textbf{5}~(2011), 745--773.

\bibitem{BMR} J. Bemelmans, Min-Oo, E. Ruh, \emph{Smoothing Riemannian metrics}, Math. Z. {\bf 188} (1984), 69--74.


\bibitem{BBG} P. Berard, G. Besson, S. Gallot, \emph{Embedding Riemannian manifolds 
by their heat kernel}, Geom. Funct. Anal. {\bf 4} (1994), no. 4, 373--398. 

\bibitem{BL} J. Bergh, J. L\"ofstr\"om, \emph{Interpolation spaces: an introduction}, Grundlehren der Mathem Wissenschaften \textbf{223}, Springer-Verlag, 1976.



\bibitem{Bill:conv}
P. Billingsley, \emph{Convergence of probability measures}, 2nd edition, Wiley Series in Probability and Statistics, John Wiley \& Sons, Inc., 1999.
 
 \bibitem{BKL3}
R.~Bosi,~Y.~Kurylev,~M.~Lassas,~\emph{Reconstruction and stability in Gel'fand's inverse interior spectral problem}, Anal. PDE. \textbf{15} (2022), 273--326.

\bibitem{BBI} D. Burago, Y. Burago, S. Ivanov, \emph{A course in metric geometry}. Graduate Studies in Math. \textbf{33}, AMS, 2001.


\bibitem{BGP} Y. Burago, M. Gromov, G. Perelman, \emph{A.D. Alexandrov spaces with curvature bounded below},  Russian Math. Surveys  {\bf 47}  (1992),  1--58.

\bibitem{BILL}
D. Burago, S. Ivanov, M. Lassas, J. Lu, \emph{Quantitative stability of Gel'fand's inverse boundary problem}, to appear in Anal. PDE., arXiv:2012.04435v4.

\bibitem{Caday}
P. Caday, M. de Hoop, V. Katsnelson,  G. Uhlmann, \emph{Scattering control for the wave equation with unknown wave speed,}  ARMA {\bf 231}  (2019), 409--464. 

\bibitem{CH} E. Calabi, P. Hartman, \emph{On the smoothness of isometries}, Duke Math. J. {\bf 37} (1970), 741--750.




\bibitem{Ch:finite}  J. Cheeger, \emph{Finiteness theorems for Riemannian manifolds}, Amer. J. Math. {\bf 92}  (1970),  61--74.
















\bibitem{Coifman2}
R. Coifman, S. Lafon, \emph{Diffusion maps}, Appl. Comput. Harmon. Anal. \textbf{21} (2006), no. 1, 5--30.


\bibitem{Coifman1} 
R. Coifman, S. Lafon, A. Lee, M. Maggioni, B. Nadler, F. Warner, S. Zucker, \emph{Geometric diffusions as a tool for harmonic analysis and structure definition of data: Diffusion maps}, PNAS \textbf{102} (2005), 7426--7431. 



\bibitem{FIKLN} C.~Fefferman,~S.~Ivanov, Y.~Kurylev,~M.~Lassas,~H. Naranayan,
~\emph{Reconstruction and interpolation of manifolds I: The geometric Whitney problem}, Found. Comp. Math.
\textbf{20} (2020), 1035--1133.

\bibitem{FILLN} C.~Fefferman,~S.~Ivanov,~M.~Lassas,~J.~Lu,~H. Naranayan,
~\emph{Reconstruction and interpolation of manifolds II: Inverse problems for Riemannian manifolds with partial distance data}, to appear in Amer. J. Math., arXiv:2111.14528.

\bibitem{FMN}
C. Fefferman, S. Mitter,  N. Narayanan, \emph{Testing the manifold hypothesis}, J. Amer. Math. Soc. \textbf{29} (2016), no. 4, 983--1049.


\bibitem {Fuk_JDG2}
K. Fukaya, \emph{Collapsing Riemannian manifold to ones of lower
dimension}, J. Differ. Geom. {\bf 25} (1987), 139--156.

\bibitem {Fuk_inv}
K. Fukaya, \emph{Collapsing of Riemannian manifolds and eigenvalues
 of Laplace operator},  Invent. Math.  {\bf 87}  (1987),  no. 3, 517--547.

\bibitem {Fuk_JDG}
K. Fukaya, \emph{A boundary of the set of the Riemannian manifolds with
 bounded curvatures and diameters}, J. Differ. Geom.  {\bf 28}  (1988),  no. 1, 1--21.
 
\bibitem {Fuk:II}
K. Fukaya, \emph{Collapsing Riemannian manifolds to ones with lower
dimension II}, J. Math. Soc. Japan {\bf 41} (1989), 333--356.

\bibitem {Fuk:haus}
K. Fukaya, \emph{Hausdorff convergence of Riemannian manifolds and its applications}, Recent Topics in Differential and Analytic Geometry, Adv. Stud. Pure Math. {\bf 18} (1990), 143--238.






\bibitem{FY:fundgp}
K. Fukaya, T. Yamaguchi, \emph{The fundamental groups of almost nonnegatively curved  manifolds}, Ann. of Math. {\bf 136} (1992), 253--333.
    
\bibitem{Gel} I.~Gel'fand,~\emph{Some aspects of functional analysis and algebra},~Proc. ICM.~\textbf{1}~(1954),~253--277.

\bibitem{GTr} D. Gilbarg, N. Trudinger, \emph{Elliptic partial differential equations of second order}, Grundlehr. der Mathem. Wissensch. {\bf 224}, Springer, 1983.
    


\bibitem{GW:lipschitz}
R. Greene, H. Wu, \emph{Lipschitz convergence of Riemannian manifolds}, Pacific J. Math. {\bf 131}  (1988), 119--141.


\bibitem{GKLU1} A. Greenleaf, Y. Kurylev, M. Lassas, G. Uhlmann, \emph{Full-wave invisibility of
active devices at all frequencies}, Comm.  Math. Phys. {\bf 275} (2007), 749--789.



\bibitem{GLU}  
A. Greenleaf, M. Lassas, G. Uhlmann, \emph{On nonuniqueness for Calder\'on's inverse problem}, Math. Res. Lett. {\bf 10} (2003), 685--693.


\bibitem{GLP}
M. Gromov, \emph{Structures m\'etriques pour les vari\'et\'es riemanniennes}, Edited by J. Lafontaine and P. Pansu, Textes Math\'ematiques {\bf 1}, CEDIC, Paris, 1981.

\bibitem{Gui}
C. Guillarmou, \emph{Lens rigidity for manifolds with hyperbolic trapped sets.} J. Amer. Math. Soc. {\bf 30} (2017),  561–599. 

\bibitem{H50}
P. Hartman, \emph{On the local uniqueness of geodesics}, Amer. J. Math. \textbf{72} (1950), 723--730.
 


\bibitem{Kas:measII} A. Kasue, 
\emph{Measured Hausdorff convergence of Riemannian manifolds and Laplace
operators II},   Lecture
Notes in Pure Appl. Math. {\bf 143} (1992),  97--111. 

\bibitem{K93}
A. Kasue, \emph{Measured Hausdorff convergence of Riemannian manifolds with Laplace operators}, Osaka J. Math. \textbf{30} (1993), 613--651.


 
\bibitem{KKL} A.~Katchalov,~Y.~Kurylev,~M.~Lassas,~\emph{Inverse boundary spectral problems},  Pure and Applied Mathematics, \textbf{123}, CRC-press,~2001.


\bibitem{KKLM}  A.~Katchalov,~Y.~Kurylev,~M.~Lassas, N. Mandache,~\emph{Equivalence of time-domain inverse problems and boundary spectral problems},~Inverse Probl.~\textbf{20}~(2004),~419--436.

\bibitem{Kato} T. Kato, \emph{Perturbation theory for linear operators}, Springer, 1976.

\bibitem{KSU}
C. Kenig, J. Sj\"ostrand, G. Uhlmann, \emph{The Calder{\'o}n problem with partial data},  Ann. of Math. {\bf 165} (2007),  no. 2, 567--591.

\bibitem{KrKuLa}
K. Krupchyk, Y. Kurylev, M. Lassas, \emph{Inverse spectral problems on a closed manifold}, Journal de Mathematique Pures et Appliquees {\bf 90} (2008), 42--59. 

\bibitem{Kry} N. Krylov, \emph{Lectures on elliptic and parabolic equations in Sobolev spaces}, Graduate Studies in Math. \textbf{96}, AMS, 2008.


\bibitem{KOP} Y.~Kurylev,~L.~Oksanen,~G.~Paternain, \emph{Inverse problems for the connection Laplacian.} {J. Differ. Geom.} {\bf 110} (2018), 457--494.


\bibitem{LLY}
M. Lassas, J. Lu, T. Yamaguchi, \emph{Inverse spectral problems for collapsing manifolds II: Quantitative stability of reconstruction for orbifolds}, in preparation.


\bibitem{LO}
M.~Lassas, L.~Oksanen,
\emph{Inverse problem for the Riemannian wave equation with Dirichlet
  data and Neumann data on disjoint sets}, Duke Math. J. \textbf{163} (2014), 1071--1103.


\bibitem{LL}
C. Laurent, M. Leautaud, 
\emph{Quantitative unique continuation for operators with partially analytic coefficients. Application to approximate control for waves.} 
J. Eur. Math. Soc. {\bf 21} (2019),  957–1069.

\bibitem{MZ} D. Montgomery, L. Zippin, \emph{Topological transformation groups}, R.Kruger, 1974.


\bibitem{Na2}
A. Nachman, \emph{Global uniqueness for a two-dimensional inverse boundary value problem}, Ann. of Math. {\bf 143} (1996), no. 1, 71--96.


\bibitem{Oneill}  B. O'Neill, \emph{Semi-Riemannian geometry with applications to relativity}, Acad. Press, New York-Boston-London, 1983.


\bibitem{Pets:conv}
S. Peters, \emph{Convergence of Riemannian manifolds}, Compos. Math. {\bf  62}  (1987),   3--16. 



\bibitem{Pet}
A. Petrunin, \emph{Parallel transportation for Alexandrov space with curvature bounded below}, Geom. Funct. Anal. {\bf 8}  (1998),  no. 1, 123--148.

\bibitem{Shf} S. Z. Shefel, \emph{Smoothness of a conformal map of Riemannian spaces}, Sibirsk. Mat. Zh. {\bf 23} (1982), 153--159.





\bibitem{Singer}
A. Singer, \emph{Mathematics for Cryo-Electron Microscopy}, Proceedings of ICM 2018, 3995–4014, 2019.


\bibitem{SU} P.~Stefanov,~G.~Uhlmann,~\emph{Stability estimates for the hyperbolic Dirichlet-to-Neumann map in anisotropic media},~J. Funct. Anal.~\textbf{154}~(1998),~330--357.

\bibitem{SU2}
P. Stefanov, G. Uhlmann, \emph{Stable determination of generic simple metrics from the hyperbolic
Dirichlet-to-Neumann map}, Int. Math. Res. Not. \textbf{17} (2005), 1047--1061.

\bibitem{SUV1}
P. Stefanov, G. Uhlmann, 
A. Vasy, \emph{Local and global boundary rigidity and the geodesic X-ray transform in the normal gauge.} Ann. of Math. {\bf  194} (2021),  1–95.

\bibitem{SUV2}
P. Stefanov, G. Uhlmann, 
A. Vasy, \emph{Boundary rigidity with partial data.} J. Amer. Math. Soc. {\bf 29} (2016),  299–332. 



\bibitem{SyU}
J. Sylvester, G. Uhlmann, \emph{A global uniqueness theorem for an inverse boundary value problem}, Ann. of Math. {\bf 125} (1987), no. 1, 153--169.

\bibitem{Tataru} D. Tataru, \emph{Unique continuation for solutions to PDE's; between H\"ormander's theorem and Holmgren's theorem},  Comm. PDE. \textbf{20}  (1995), 855--884. 

\bibitem{PDE-Taylor} M. Taylor, \emph{Partial differential equations I}, 
 2nd Ed., Springer,  New York, 2011.

\bibitem{Tay} M. Taylor, \emph{Partial differential equations III},  2nd Ed., Springer, New York, 2011.



\bibitem{ISOMAP}  
J. Tenenbaum, V. de Silva, J. Langford, \emph{A global geometric framework for nonlinear dimensionality reduction}, Science \textbf{290} (2000), 2319--2323.


\bibitem{Thurston} W. Thurston, \emph{The geometry and topology of three-manifolds}, Lecture notes, Princeton University, 1980.


\bibitem{Tri} H. Triebel, \emph{Interpolation theory, function spaces, differential operators}, second edition. Johann Ambrosius Barth, Heidelberg, 1995.


\bibitem{WM}
X. Wang, J. Marron, \emph{A scale-based approach to finding effective dimensionality in
manifold learning}, Electron. J. Stat. \textbf{2} (2008), 127–148.

\bibitem{UW}
G. Uhlmann, A.  Vasy, \emph{The inverse problem for the local geodesic ray transform.} Invent. Math. {\bf 205} (2016),  83–120.




\end{thebibliography}
\end{document}